\documentclass[reqno,11pt]{amsart}
\usepackage[utf8x]{inputenc}
\usepackage[unicode]{hyperref}
\usepackage{mathrsfs,bbm}
\usepackage{graphicx}

\usepackage[english]{babel}
\usepackage{verbatim}
\usepackage{supertabular}
\usepackage{stmaryrd}
\usepackage{longtable}
\usepackage{varioref}
\usepackage{calc}
\usepackage{ifthen}
\usepackage{indentfirst}
\usepackage{fancyhdr}
\usepackage{enumerate}
\usepackage{float}
\usepackage{xcolor}
\usepackage{multicol}
\usepackage[BCOR4mm,DIV10]{typearea}
\usepackage[refpage]{nomencl/nomencl}
\usepackage{amssymb}
\usepackage{amsxtra}
\usepackage{amsmath} 
\usepackage{amsfonts}
\usepackage{etoolbox}
\usepackage{xfrac}
\usepackage{nicefrac}

\usepackage[normalem]{ulem}

\usepackage{todonotes}
\usepackage{ucs}
\usepackage[T1]{fontenc}
\usepackage{here}
\usepackage{tikz}
\usepackage{subcaption} 
\usepackage{mathtools}
\usepackage{subcaption}
\usepackage[capitalise,noabbrev]{cleveref}
\usetikzlibrary{quotes,angles}
\RequirePackage{enumitem}
\RequirePackage{hyperref}
\RequirePackage{cleveref}
\RequirePackage{ifthen}

\makeatletter
\newcommand{\Item}[1][]{
	\ifthenelse{\equal{#1}{}}{%
		\item%
	}{%
		\refstepcounter{Item}%
		\csname @item\endcsname[#1]%
		\protected@edef\@currentlabel{#1}%
		\protected@edef\cref@currentlabel{[enumi][][]{#1}}%
	}%
}
\makeatother

\newcommand{\R}{\mathbb{R}}

\newcommand{\C}{\mathbb{C}}
\newcommand{\Z}{\mathbb{Z}}
\newcommand{\N}{\mathbb{N}}

\newcommand{\Sc}{\R / \Z}
\newcommand{\Ydc}{(\Sc)^2 \setminus Y_\delta(\gamma)}
\newcommand{\Ydch}{(\Sc\times [-\frac{1}{2},\frac{1}{2}]) \setminus \hat{Y}_\delta(\gamma)}

\newcommand\ANG{\mathop{\mbox{$<\!\!\!)$}}\nolimits}

\newcommand{\bleq}[1][ ]{\stackrel{\mathmakebox[\widthof{(1.1)}]{#1}}\leq}

\newcommand{\g}{\gamma}

\DeclareRobustCommand{\rchi}{{\mathpalette\irchi\relax}}
\newcommand{\irchi}[2]{\raisebox{\depth}{$#1\chi$}} 

\DeclareMathOperator{\st}{\textit{St}}

\usepackage{import}
\usepackage{xifthen}
\usepackage{pdfpages}
\usepackage{transparent}

\numberwithin{equation}{section}
\allowdisplaybreaks

\newtheorem{theorem}{Theorem}[section]
\newtheorem{proposition}[theorem]{Proposition}
\newtheorem{lemma}[theorem]{Lemma}
\newtheorem{corollary}[theorem]{Corollary}
\newtheorem{definition}[theorem]{Definition}

\newtheorem*{theorem*}{Theorem}
\hypersetup{
breaklinks=true,
colorlinks=true,
linkcolor=blue,
citecolor=blue,
urlcolor=blue,
}

\DeclareMathOperator{\dist}{dist}

 \renewcommand\S{{\mathbb S}}

\newcommand{\omitted}[1]{}

\newcommand{\Fo}{\,\,\,\text{for }\,\,}
\newcommand{\Foa}{\,\,\,\text{for all }\,\,}

\newcommand{\AND}{\,\,\,\text{and }\,\,}

\newcommand{\heikodetail}[1]{}

\def\g{\gamma}
\def\G{\Gamma}

\def\dist{{\rm{dist}\,}}

\def\TP{{\rm{TP}}}

\def\rtp{{r_\textnormal{tp}}}

\def\HC{\mathcal{C}}
\def\HW{\mathcal{W}}

\def\HF{\mathcal{F}}
\def\HE{\mathcal{E}}

\def\mvint_#1{\mathchoice
          {\mathop{\vrule width 6pt height 3 pt depth -2.5pt
            \kern -8pt \intop}\nolimits_{\kern -3pt #1}}%
             {\mathop{\vrule width 5pt height 3 pt depth -2.6pt
                        \kern -6pt \intop}\nolimits_{#1}}%
              {\mathop{\vrule width 5pt height 3 pt depth -2.6pt
                   \kern -6pt \intop}\nolimits_{#1}}%
              {\mathop{\vrule width 5pt height 3 pt depth -2.6pt
                      \kern -6pt \intop}\nolimits_{#1}}}

\usepackage[left,mathlines,pagewise]{lineno}
\newcommand*\patchAmsMathEnvironmentForLineno[1]{%
\expandafter\let\csname old#1\expandafter\endcsname\csname #1\endcsname
\expandafter\let\csname oldend#1\expandafter\endcsname\csname end#1\endcsname
\renewenvironment{#1}%
{\linenomath\csname old#1\endcsname}%
{\csname oldend#1\endcsname\endlinenomath}}%
\newcommand*\patchBothAmsMathEnvironmentsForLineno[1]{%
\patchAmsMathEnvironmentForLineno{#1}%
\patchAmsMathEnvironmentForLineno{#1*}}%
\AtBeginDocument{%
\patchBothAmsMathEnvironmentsForLineno{equation}%
\patchBothAmsMathEnvironmentsForLineno{align}%
\patchBothAmsMathEnvironmentsForLineno{flalign}%
\patchBothAmsMathEnvironmentsForLineno{alignat}%
\patchBothAmsMathEnvironmentsForLineno{gather}%
\patchBothAmsMathEnvironmentsForLineno{multline}%
}

\title[Optimal planar immersions]{Optimal
planar immersions of prescribed winding number and 
Arnold invariants}
\author{Anna Lagemann}
\address[A.~Lagemann]{
\newline%
RWTH Aachen University,\newline%
Institut f\"ur Mathematik,\newline%
Templergraben 55,
52062 Aachen,
Germany}
\email{lagemann@eddy.rwth-aachen.de}%
\author{Heiko von der Mosel}
\address[H.~von~der~Mosel]{
\newline%
RWTH Aachen University,\newline%
Institut f\"ur Mathematik,\newline%
Templergraben 55,
52062 Aachen,
Germany}
\email{heiko@instmath.rwth-aachen.de}


\date{\today}

\DeclareRobustCommand{\SkipTocEntry}[5]{}
\begin{document}

\begin{abstract}
Vladimir Arnold defined three invariants for generic planar
immersions, 
i.e.\ planar
curves whose self-intersections are all
transverse double points. We use a variational approach to study these 
invariants by investigating a suitably truncated knot energy,
the tangent-point energy. 
We prove existence of energy minimizers for each 
truncation parameter $\delta >0$
in a class of immersions with prescribed winding number
and Arnold invariants, 
and 
establish Gamma convergence of the truncated tangent-point
energies to  a limiting renormalized tangent-point energy
as $\delta\to 0$.
Moreover, we show that any sequence of minimizers subconverges in $C^1$,
and the
corresponding limit curve has the same topological
invariants, self-intersects exclusively 
at  right angles,  and
minimizes the renormalized tangent-point
energy among all curves with right self-intersection angles. In addition,
the limit curve is an almost-minimizer for
all of the original truncated
tangent-point energies as long as the truncation parameter $\delta$
is sufficiently small. Therefore, this limit curve serves as 
an ``optimal'' curve  in the class of generic planar 
immersions with prescribed winding number and Arnold invariants. 
\end{abstract}

\maketitle
\section{Introduction}
Planar closed immersed curves $\g:\R/\Z\to\R^2$
can be classified in terms of their 
winding number $W(\g)$ which equals the degree of the map
$s\mapsto\g'(s)/|\g'(s)|$. Roughly speaking, the winding number
counts the number of turns of the tangent vector while traveling 
once along the curve. The Whitney-Graustein theorem 
\cite[Theorem 1]{whitney}\footnote{For an alternative contact-geometric
proof of the Whitney-Graustein theorem see H. Geiges \cite{geiges_2009}.}
 states: 
\emph{Two planar curves are regularly homotopic if and only if they have the 
same winding number.} V.I. 
Arnold proposed in \cite{Arnold,Arnold_book} 
a much deeper plan of study and introduced three invariants $J^+$, $J^-$ 
and $\st$, which are locally constant on \emph{generic} 
immersions, i.e.\ on 
curves where all self-intersections are transverse double points. 
It turns out that the combinatorics of generic loops is potentially
as rich as that of knots in $3$-space.

Our main motivation  is to explore with variational
tools -- similarly 
in spirit as in
geometric knot theory -- the different connected components
of generic planar immersions with prescribed Arnold invariants.
One may ask, e.g.,   to 
what extent such components of
generic curves 
are potential wells for some energy? 
One typically expects 
energy
minimizers and critical points to be particularly
interesting  
configurations in the given topological class. 
Among other things, one hopes for nice geometric properties, so
that these
minimizers can serve as optimal curves in the class of 
immersions with prescribed invariants. 

In geometric knot theory one minimizes so-called \emph{knot
energies} to find optimal shapes in given knot classes. This is
possible since knot energies blow up along sequences of curves that
converge to a limit curve with self-intersections. Immersed planar curves
with winding number different from $\pm 1$, however, must have
self-intersections, so we need to desingularize the knot energy near
double points to obtain finite energy values.  

Since Arnold's invariants
are well-defined on $C^1$-immersions we choose here the 
\emph{tangent-point energy} $\TP_q$ for some $q>2$ whose natural
energy space is a fractional Sobolev space that embeds compactly
into $C^1$.
This energy is well-suited for our variational
approach to study the space
of generic immersions, since
it is uniquely minimized by the (round) circle among all
embedded closed space curves, and in every prescribed knot class there
is a minimizing knot. Additional critical knots were found by 
symmetric criticality. Moreover,  
long-time existence for
the (Banach-)gradient flow was established,  and subconvergence of
its solution to critical points holds for
variants of the tangent-point energy whose energy space is a Hilbert space. In that case,
the Palais-Smale condition holds, 
which opens up the possibility of Ljusternik-Schnirelman theory. 
In addition, the tangent-point energy induces
a complete Riemannian metric on the open subset of embeddings in that Hilbert space,
together with the existence of  distance-minimizing geodesics connecting
any given pair of knots in the same knot class;
see Section \ref{sec:1.2} for more details and references.

\begin{figure}
  \begin{minipage}{\textwidth}
  \centering
\begin{subfigure}{0.27\textwidth}
    \includegraphics[trim={0 2cm 0 0}, clip, width=\textwidth]{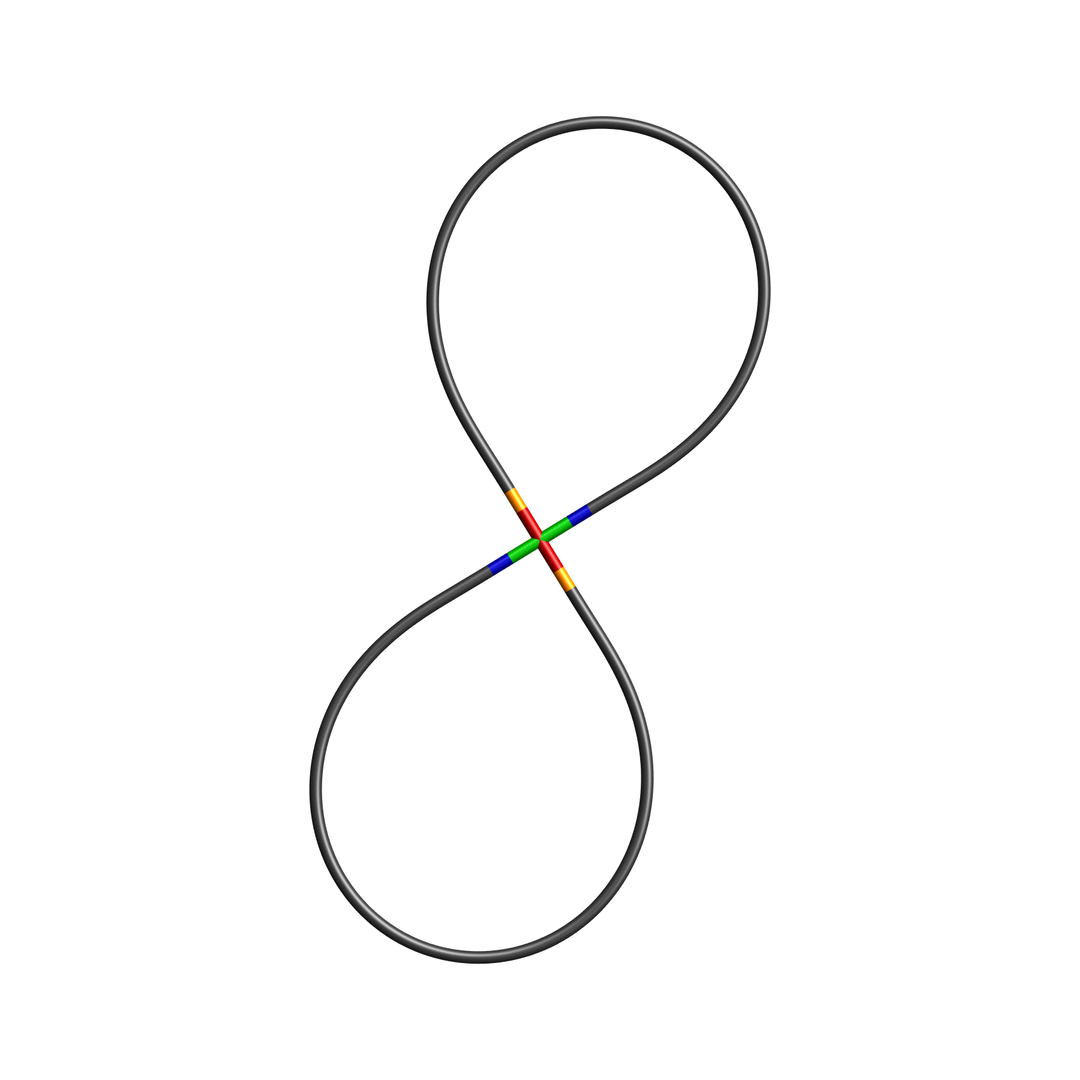}
    \caption{$J^+=0$, $J^-=-1$, $\st=0$, $W=0$.}
\end{subfigure}
\hspace{0.4cm}
\begin{subfigure}{0.27\textwidth}
    \includegraphics[trim={0 2cm 0 0}, clip, width=\textwidth]{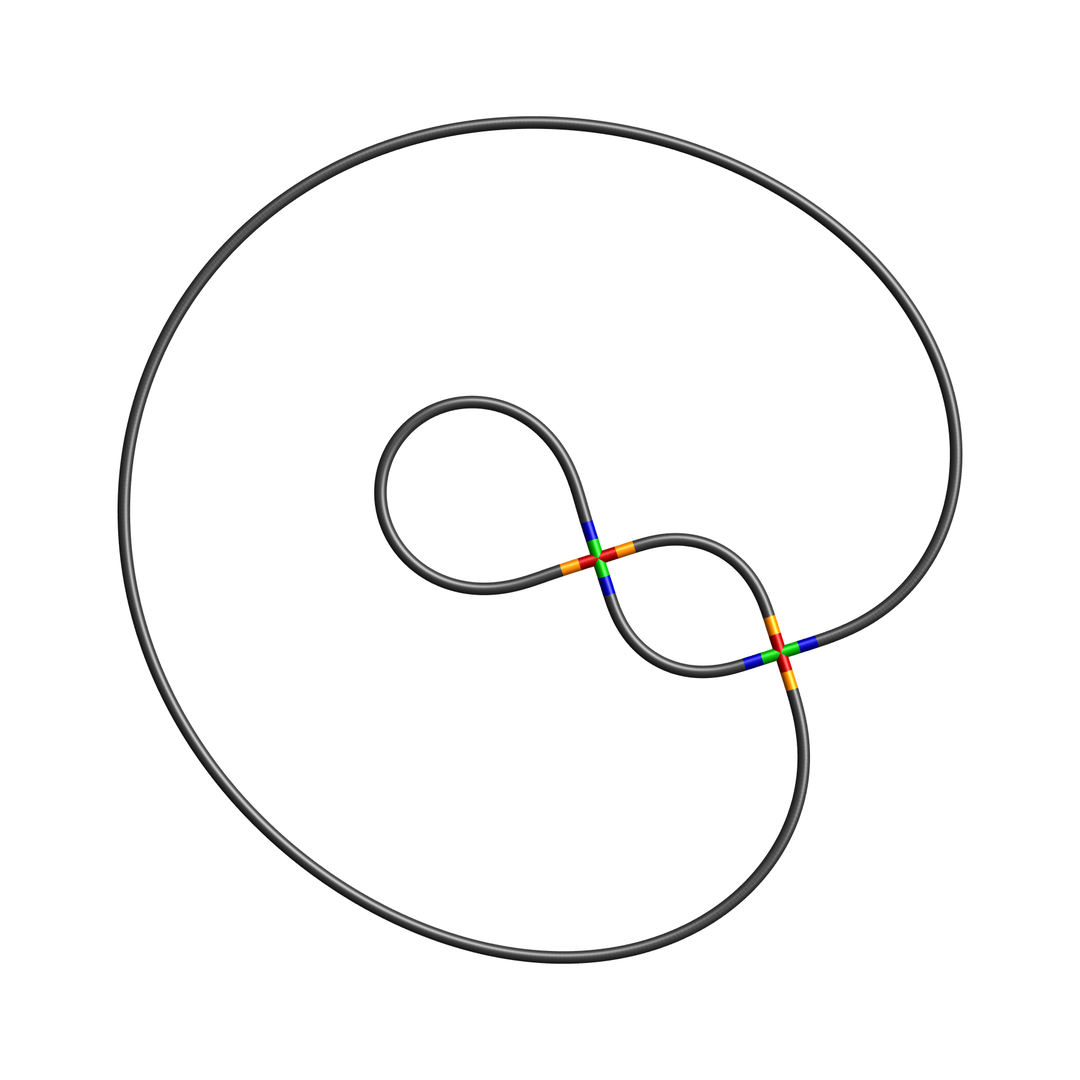}
    \caption{$J^+=0$, $J^-=-2$, $\st=0$, $W=1$.}
\end{subfigure}
\hspace{0.4cm}
\begin{subfigure}{0.27\textwidth}
  \includegraphics[trim={0 2cm 0 0}, clip, width=\textwidth]{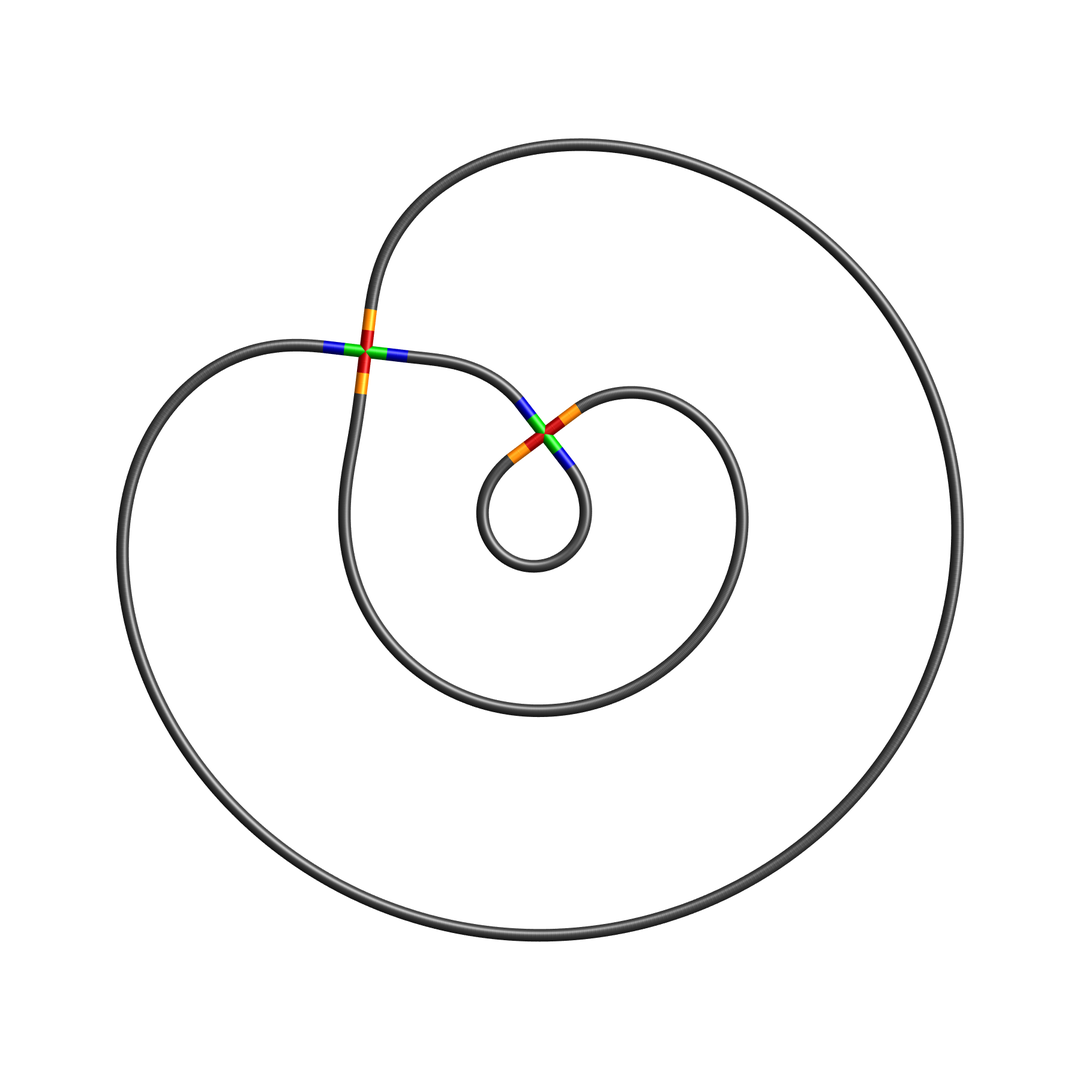}
  \caption{$J^+=-6$, $J^-=-8$, $\st=3$, $W=3$.}
\end{subfigure}
\hfill
\begin{subfigure}{0.27\textwidth}
    \includegraphics[trim={0 2cm 0 0}, clip, width=\textwidth]{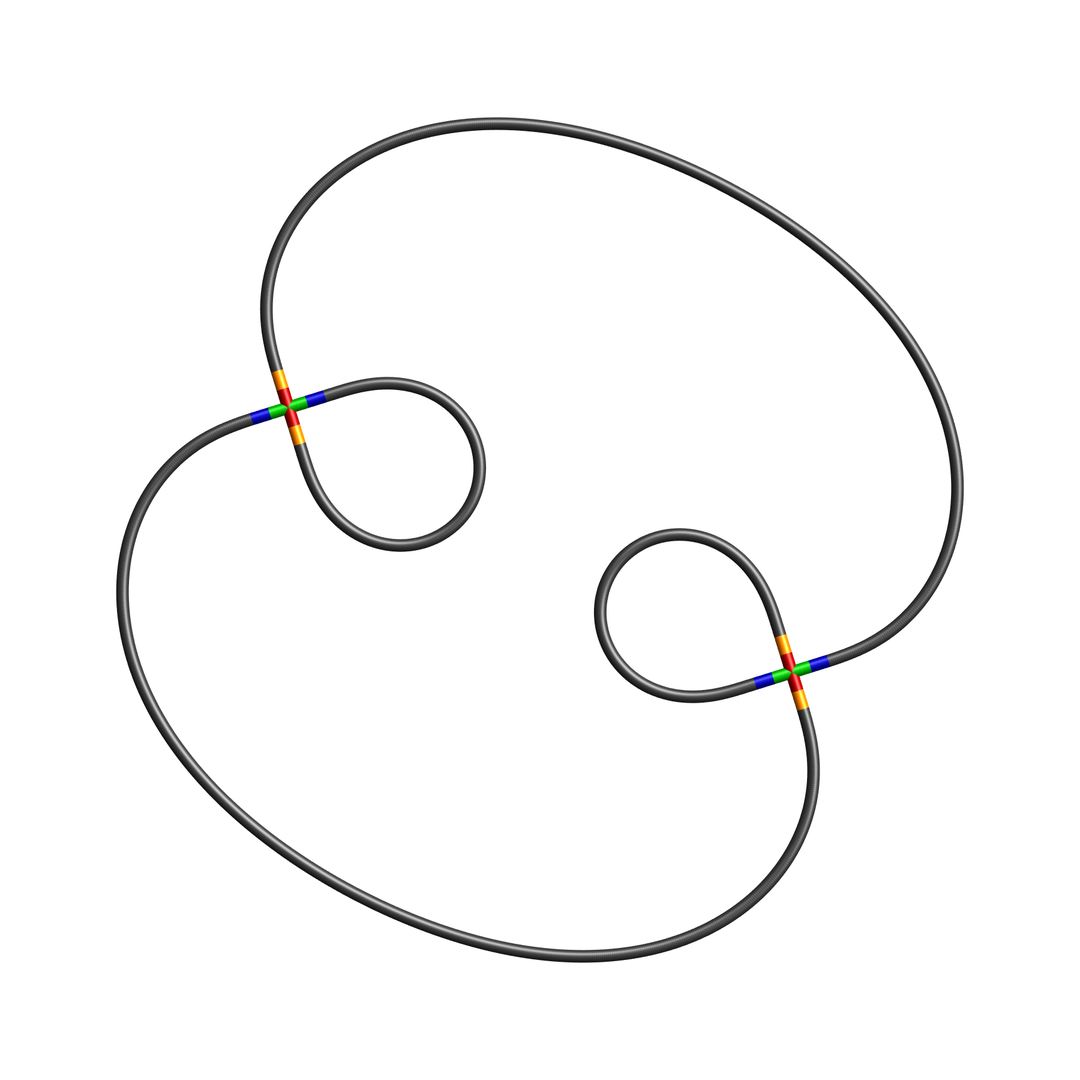}
    \caption{$J^+=-4$, $J^-=-6$, $\st=2$, $W=3$.}
\end{subfigure}
\hspace{0.4cm}
\begin{subfigure}{0.27\textwidth}
  \includegraphics[trim={0 2cm 0 0}, clip, width=\textwidth]{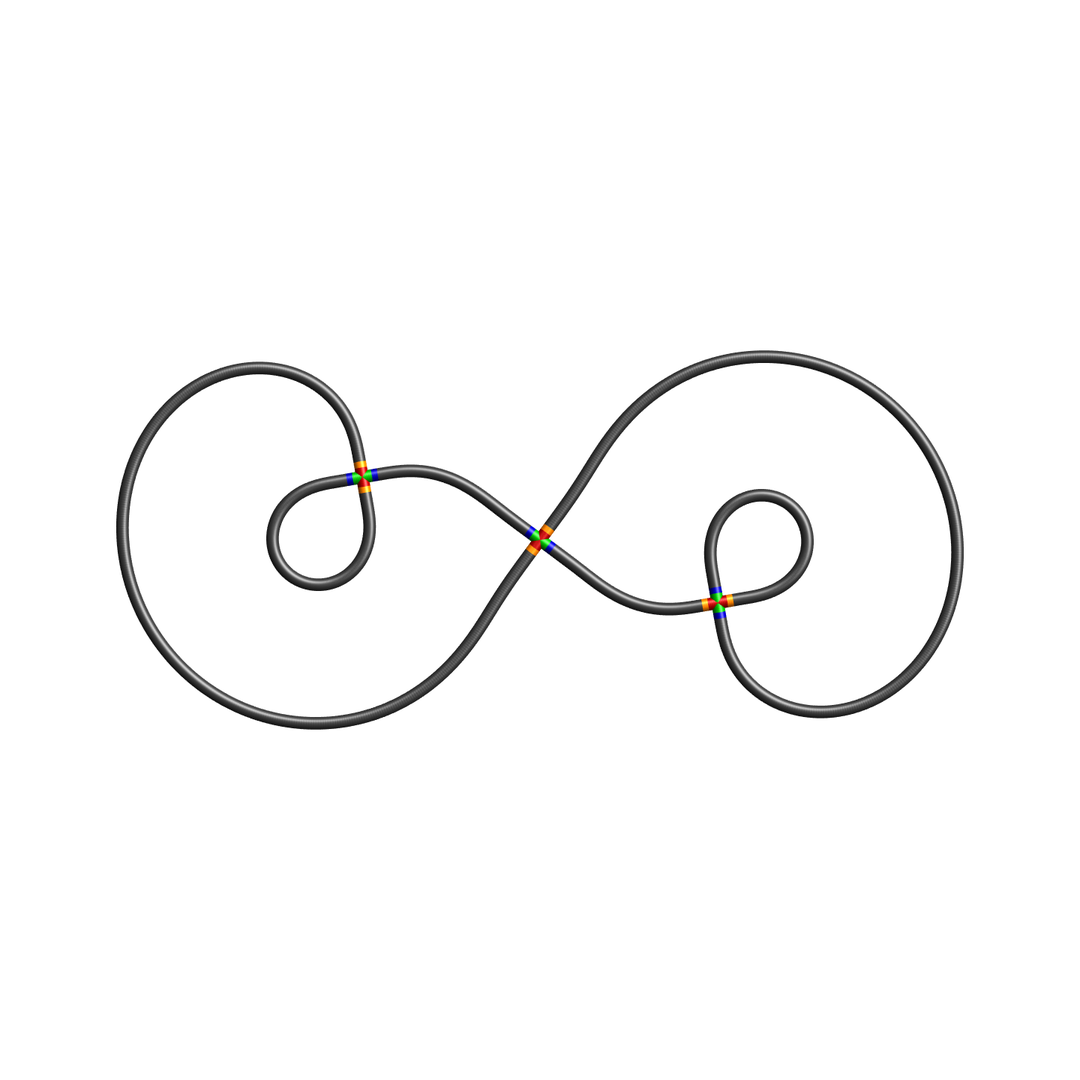}
  \caption{$J^+=-4$, $J^-=-7$, $\st=2$, $W=0$.}
\end{subfigure}
\hspace{0.4cm}
\begin{subfigure}{0.27\textwidth}
  \includegraphics[trim={0 2cm 0 0}, clip, width=\textwidth]{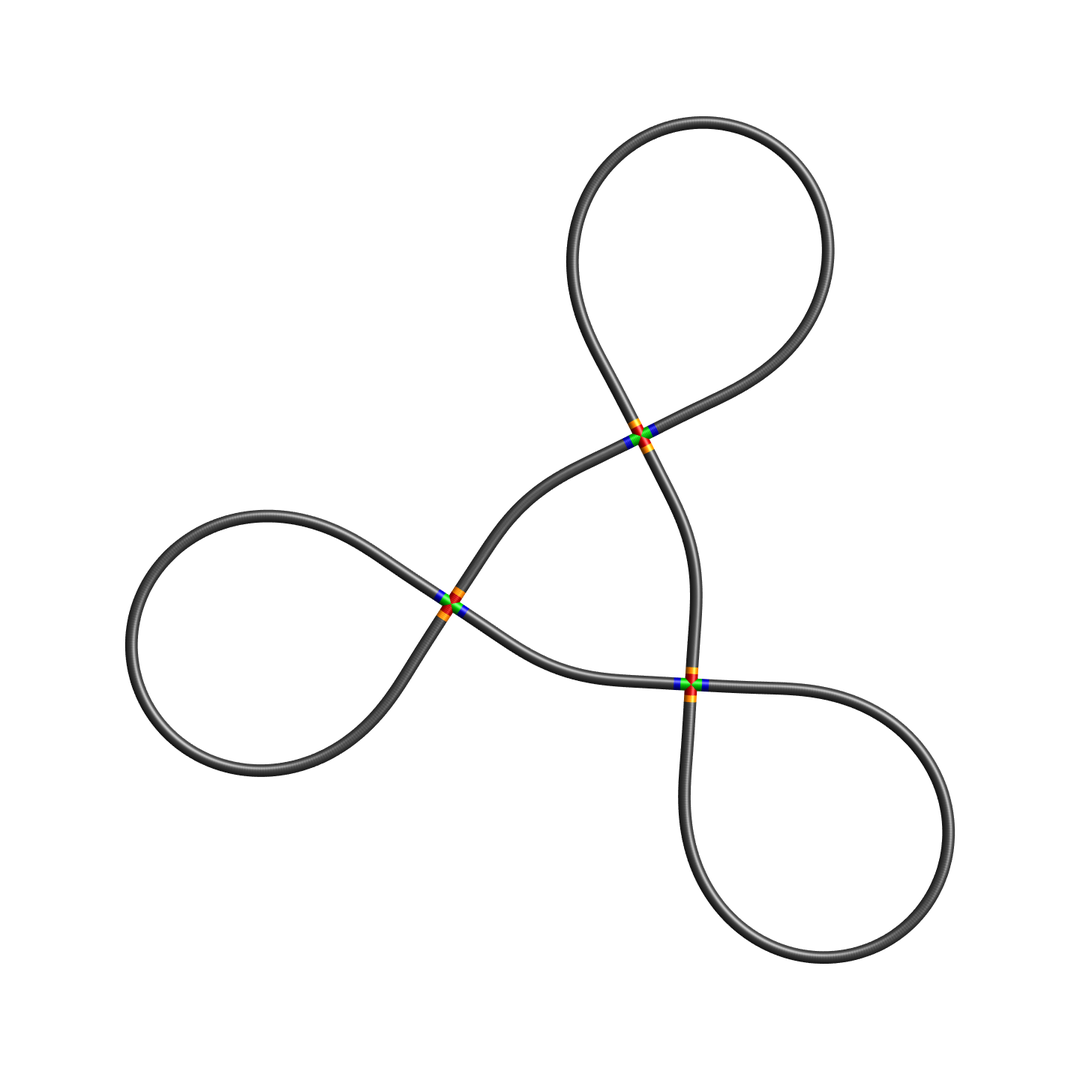}
  \caption{$J^+=0$, $J^-=-3$, $\st=0$, $W=2$.}
  
\end{subfigure}
\hfill
\begin{subfigure}{0.27\textwidth}
  \includegraphics[trim={0 4cm 0 0}, clip, width=\textwidth]{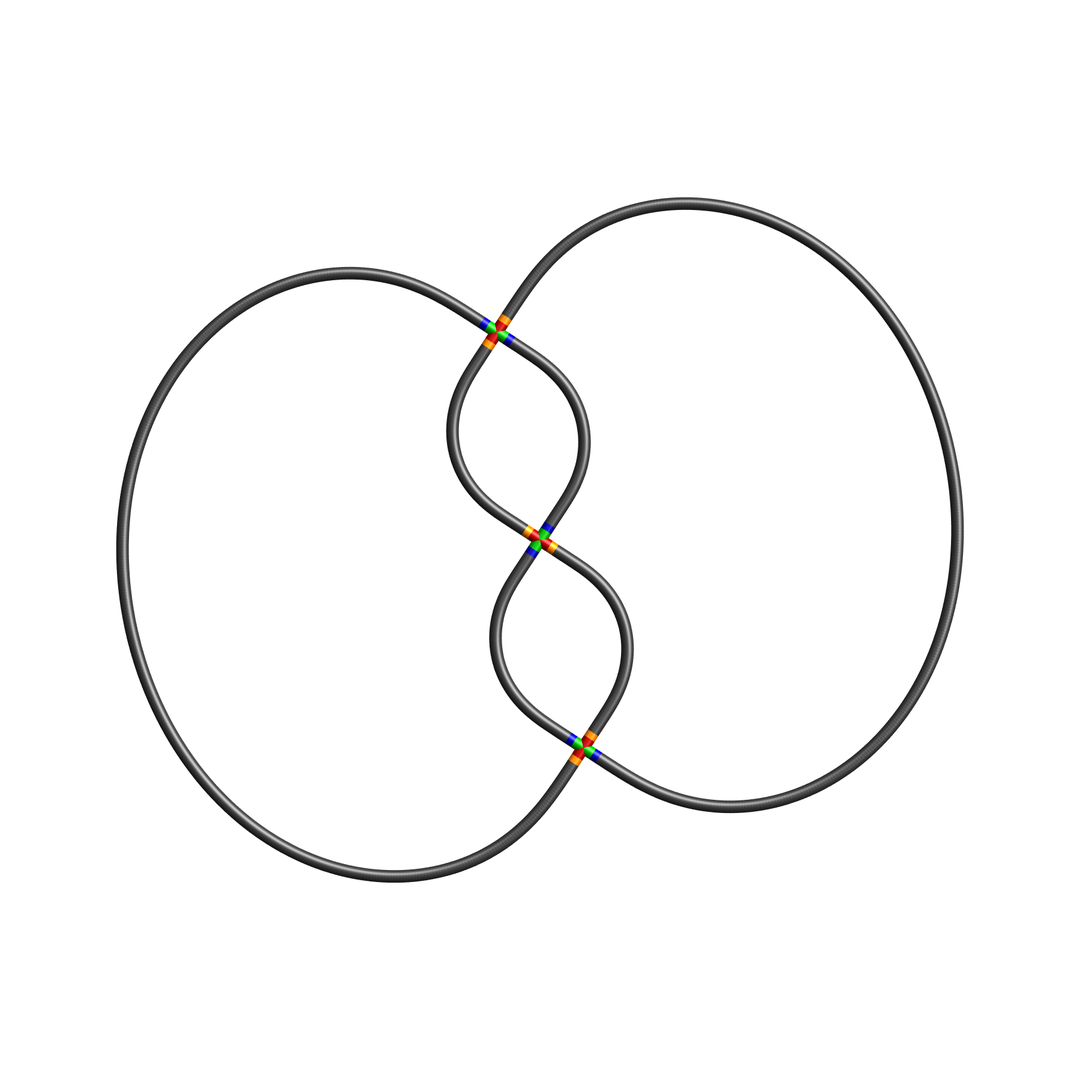}
  \caption{$J^+=2$, $J^-=-1$, $\st=0$, $W=0$.}
\end{subfigure}
\hspace{0.4cm}
\begin{subfigure}{0.27\textwidth}
    \includegraphics[trim={0 4cm 0 0}, clip, width=\textwidth]{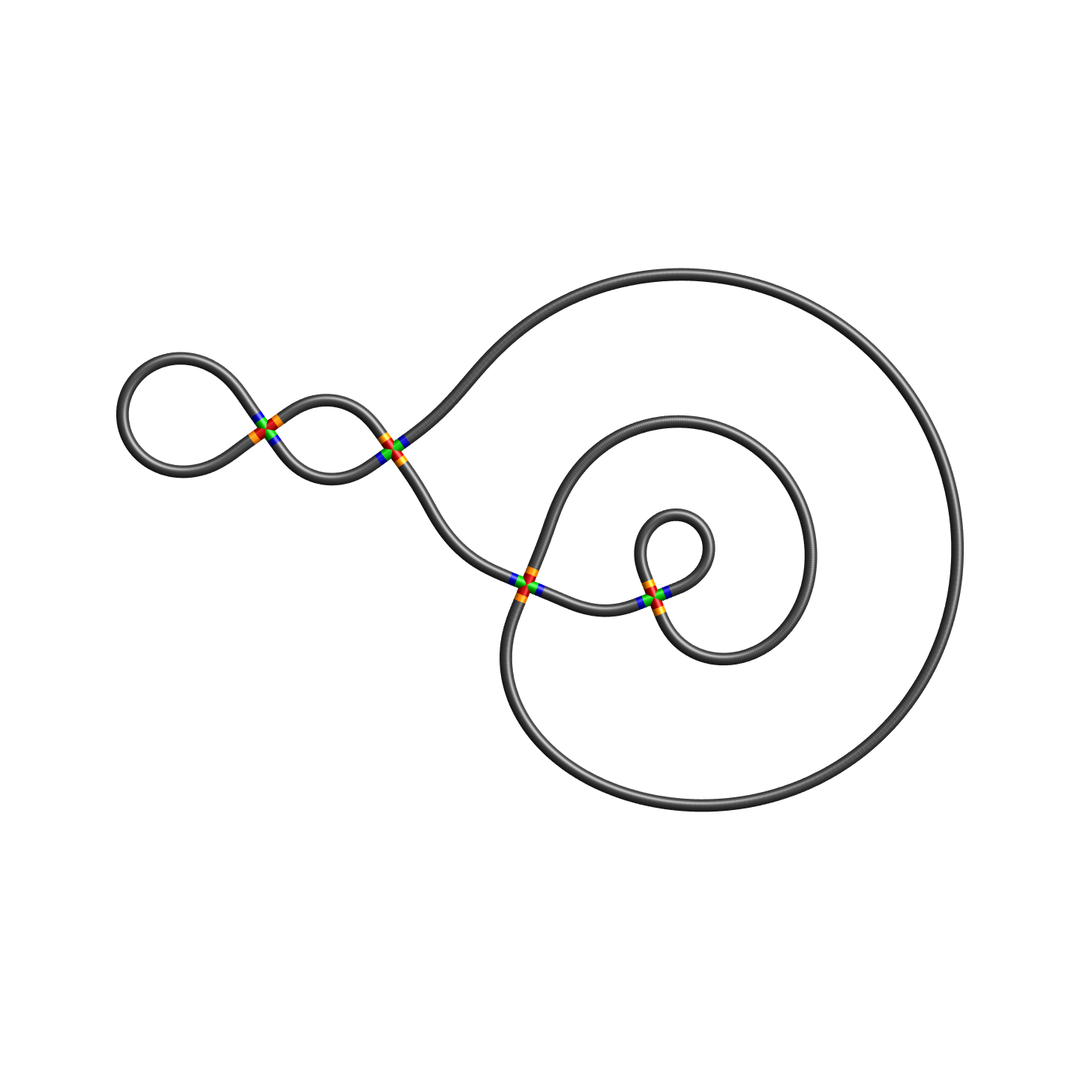}
    \caption{$J^+=-6$, $J^-=-10$, $\st=3$, $W=3$.}
\end{subfigure}
\hspace{0.4cm}
\begin{subfigure}{0.27\textwidth}
  \includegraphics[trim={0 4cm 0 0}, clip, width=\textwidth]{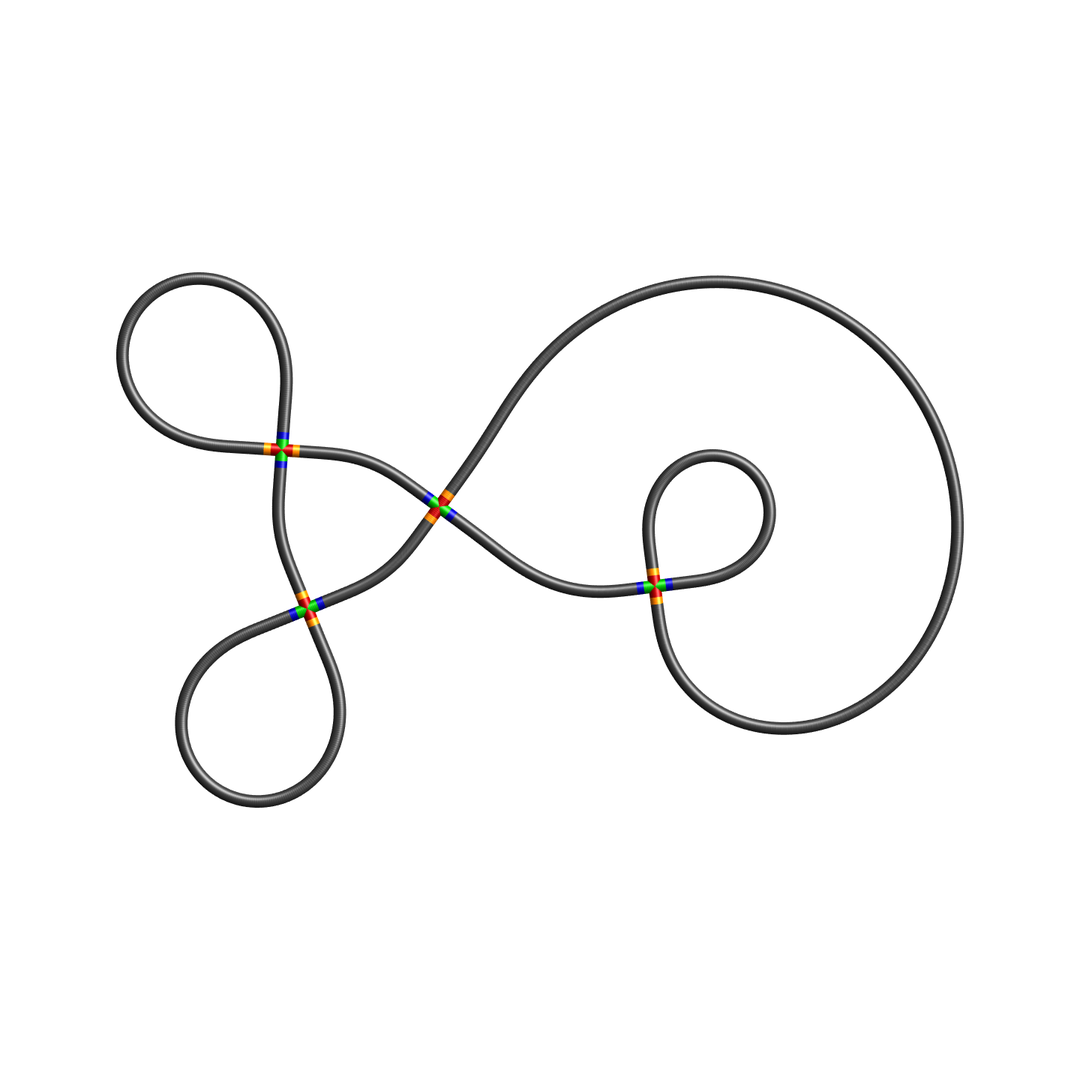}
  \caption{$J^+=-2$, $J^-=-6$, $\st=1$, $W=3$.}
\end{subfigure}
\caption{Examples of numerically computed\protect\footnotemark
minimizers of the energy $\TP_{q,\delta}$ 
in  admissibility 
classes $\HF(\eta,j_\pm,s,\omega)$ of curves with different 
prescribed  Arnold invariants
$J^+=j_+,$ $J^-=j_-$, $\st=s$,
and winding numbers $W=\omega$. All these minimizers seem to
self-intersect exclusively in 
right-angles.
The competing curves in 
$\HF(\eta,j_\pm,s,\omega)$ are affine linear within arclength $\eta$ around
every self-intersection; see Definition \ref{admissibile_curves}. 
The ratio of the truncation parameter $\delta$
(green or red) to $\eta$ (blue or yellow)
is approximately $1/2$ in our computations.
\label{fig:9pictures}
}
\end{minipage}
\end{figure}

To treat immersions with double points we
cut out small $\delta$-neighborhoods
of the curves'  self-intersections to
obtain the \emph{truncated} 
and therefore desingularized variant $\TP_{q,\delta}$ 
of the
tangent-point energy. We prove for each $\delta>0$
the existence of $\TP_{q,\delta}$-minimizers in
a class of generic planar immersions with fixed winding number $W$
and prescribed Arnold invariants $J^+,$ $J^-$, and $\st$.
Rescaling the truncated energies with a suitable power of the
truncation parameter $\delta$ adapted to the energy blow-up
near self-intersections, and then
sending $\delta$ to zero,
we establish Gamma convergence of the rescaled energies to a limiting
energy $R_q$, the \emph{renormalized tangent-point energy},
that measures only the intersection angles between the tangent lines
at the
double points.  The limiting process $\delta\to 0$ also yields
a $C^1$-convergent subsequence of minimizers to a limit curve with 
the same topological invariants and with exclusively right angles at every 
self-intersection. Moreover, this limit curve minimizes the 
renormalized
tangent-point  energy $R_q$ among all
generic immersions of the given topological class
with right self-intersection angles, and it 
almost-minimizes in addition
the original truncated energies $\TP_{q,\delta}$ 
in the full topological class for all sufficiently small $\delta$.
Because of these minimizing properties these  limiting immersions may
be regarded as  optimal configurations in their topological classes.
The precise mathematical statements are formulated in Section
\ref{sec:1.3}, but for  a first impression of the shapes of
minimizing
configurations for various prescribed values of winding number and
Arnold invariants we refer the reader to the numerically computed 
minimizers depicted in 
Figure \ref{fig:9pictures}.

\subsection{Arnold invariants}\label{sec:1.1}
\begin{figure}[h]
  \centering
  \includegraphics[scale=0.85]{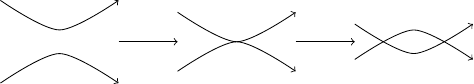}
  \caption{Homotopy through a positive direct self-tangency.}
  \label{positive_crossing}
\end{figure}

If a regular homotopy between two generic planar immersions
leaves the class of generic immersions,
different 
degeneracies could happen: intermediate
curves with self-tangencies or 
triple points, or more complicated self-intersections.
Arnold \cite{Arnold,Arnold_book}
defined the \emph{discriminant $\Delta$} as the 
set of all non-generic immersions. 
It turns out that the discriminant 
contains a ``good'' part, which is a submanifold of codimension one in 
the space of all planar
$C^1$-immersions. 
This good part consists of three different types of 
curves that are generic apart from exactly one degenerate 
self-intersection. The degeneration can manifest itself either as a 
direct self-tangency of multiplicity two, where the tangent vectors 
point in the same direction, an inverse self-tangency of
multiplicity two, 
where the tangent vectors point in opposite directions, or a transverse 
triple point. We denote these submanifolds by $\Delta^d$, $\Delta^i$ and 
$\Delta^t$, respectively. To any transverse intersection of a path in 
the space of immersions with the good part $\Delta^d\cup
\Delta^i\cup\Delta^t$ of the discriminant one can 
then assign a well-defined sign, e.g.,\ an intersection with $\Delta^d$ 
or with $\Delta^i$ is positive if the number of double points 
increases 
along the homotopy, see \cref{positive_crossing}. 

\footnotetext{We thank Henrik Schumacher
who taught us how to use his code developed
for knot energies on embedded space curves 
\cite{yu-etal_2021b}, and how
to adapt it to the present situation  of  planar immersed 
curves.}

\begin{figure}[h]
    \centering
    \def\svgwidth{0.75\columnwidth}
    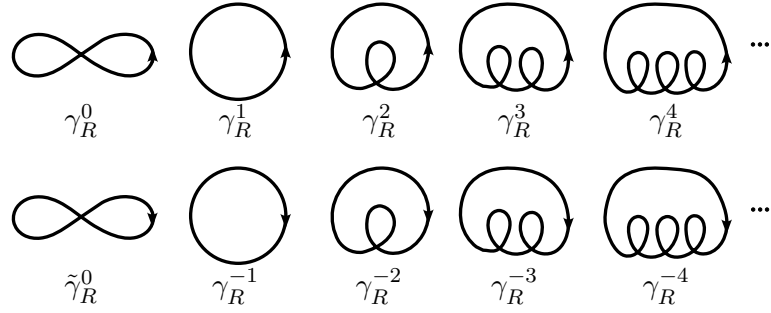

    \caption{Arnold's representative curves $\g_R^\omega$
    for various winding numbers 
    $\omega \in \Z$ (\cite[cf.\ Figure 5]{Arnold}).}
    \label{representative_curves}
\end{figure}

The invariants are then defined as follows: For any winding number, 
fix a representative curve $\g^\omega_R$; 
see \cref{representative_curves}.

Arnold prescribes 
values of the three invariants on the representative curves as follows:
\begin{alignat*}{2}
    &J^+(\gamma_R^0)=J^+(\tilde{\gamma}_R^0)=0, \quad &&J^
    +(\gamma_R^{\omega})=-2\left(|\omega|-1  \right) 
    \text{ for } \omega \in \Z \setminus \{0\}, \\ 
    &J^-(\gamma_R^0)=J^-(\tilde{\gamma}_R^0)=-1, 
    \quad &&J^-(\gamma_R^{\omega})
    =-3\left(|\omega|-1  \right) 
    \text{ for } \omega \in \Z \setminus \{0\}, \\ 
    &\st(\gamma_R^0)=\st\left(\tilde{\gamma}_R^0  \right)=0, 
    \quad &&\st(\gamma_R^{\omega})
    =\left|\omega\right|-1 \text{ for } \omega 
    \in \Z \setminus \{0\}. 
\end{alignat*}

\begin{figure}[h]
  \centering
  \includegraphics[scale=1.5]{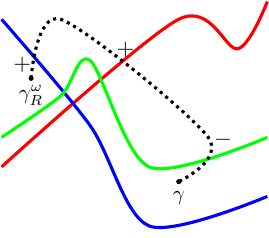}
  \caption{The space of planar $C^1$-immersions with fixed winding
  number equal to $\omega$. 
  The three colored lines without their mutual
  intersections correspond to the 
  good part $\Delta^d\cup\Delta^i\cup\Delta^t$
  of the discriminant $\Delta$. $\gamma^\omega_R$ is the 
  representative 
  curve with winding number $\omega$,
  and the dotted line represents a generic path connecting 
  $\gamma^\omega_R$ with a generic immersion $\gamma$
 of equal winding number.}
  \label{bild_def_AI}
\end{figure}

By the Whitney-Graustein theorem one can find a homotopy between any 
generic immersion and the corresponding representative curve with the 
same winding number. In fact, that path can be chosen in such a way that 
it intersects the discriminant $\Delta$
only in its good part $\Delta^d\cup\Delta^i\cup\Delta^t$
and only in 
finitely many points. Such a path will be called a \emph{generic
path}. 
For an illustration, see \cref{bild_def_AI}.
Arnold 
then establishes rules on how the invariants change at intersections with the discriminant according to the given sign as follows:

\begin{enumerate}
  \item The invariant $J^+$ increases by $2$ under a positive crossing of $\Delta^d$ and remains unchanged at crossings with $\Delta^i$ and $\Delta^t$. 
  \item The invariant $J^-$ increases by $-2$ under a positive crossing of $\Delta^i$ and remains unchanged at crossings with $\Delta^d$ and $\Delta^t$. 
  \item The invariant $\st$ increases by $1$ under a positive crossing of $\Delta^t$ and remains unchanged at crossings with $\Delta^d$ and $\Delta^i$. 
\end{enumerate}
By construction, the following equality holds:
\begin{align}
\label{connection_AI_number_self_intersections}
  J^+(\gamma)-J^-(\gamma)=\text{number of self-intersections of }\gamma. 
\end{align}
To prove existence and well-definedness of the Arnold invariants, one 
needs to show the existence of such a generic path intersecting the 
discriminant only in the codimension one part and at finitely many 
points. Furthermore, one has to verify that the definition of the 
invariants is independent of the chosen path. A detailed proof can be 
found in \cite[Chapter 5]{lagemann_diss}, which fills some gaps in the 
original proof of Arnold. 
Explicit formulas to compute the Arnold invariants as well as
sharp upper and
lower bounds on them can be found, e.g., in
\cite{shumakovich_1995,viro_1996,scheibling_1998}.
Generalizations to planar immersions of  the real line or to closed
spherical or real algebraic curves were treated in
\cite{tabachnikov_1996,viro_1996,arakawa-tetsuya_1999},
whereas the original Arnold invariants reappear as coefficients
in the Taylor expansions of curvature integral formulas
\cite{viro_1996,lanzat-polyak_2013,ito_2023}. They also play a role
in the analysis of periodic orbits of the restricted three-body problem
\cite{cieliebak-etal_2017,kim-kim_2020,cieliebak-etal_2023}, and there
is a connection to Arnold-type invariants for so-called (weak) flat
knot types of curves on surfaces
studied in the context of curve-shortening
or Finsler
geodesic flows or, more generally, in Hamiltonian
dynamics \cite{angenent_2005,hryniewicz-salamao_2013}.

\begin{figure}[h]
    \centering
    \def\svgwidth{0.35\columnwidth}
    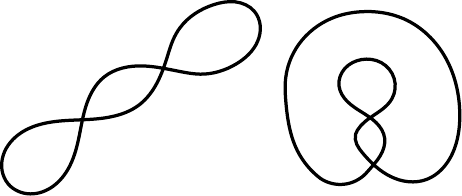

    \caption{Two curves with Arnold invariants
    $J^+=0$, $J^-=-2$, $\st=0$ and winding number $W=1$}
    \label{different_isotopy_classes_same_AI}
    \end{figure}

Because of various relations between the Arnold invariants and 
bounds on them
 (see, e.g., 
\cite[Theorem p. 43]{Arnold} and Arnold's conjectures 
\cite[pp. 60-61]{Arnold}
proven in 
\cite{shumakovich_1995,viro_1996}) it is clear that there
are combinations of values for the Arnold invariants that cannot be
realized by any generic
immersion. But even if there exist, say, two generic immersions
with 
identical winding number and
Arnold invariants it may happen that these curves cannot be deformed into
each other by regular homotopy without leaving the class
of generic immersions during that deformation;
see Figure \ref{different_isotopy_classes_same_AI}.
In other words, the set of all
generic immersions $\g$ with fixed Arnold invariants
$J^+(\g)=j_+,$ $J^-(\g)=j_-$, and $\st(\g)=s$ and winding number
$W(\g)=\omega$
splits into several open
path-connected components with respect to regular
homotopy. In
the sequel we call any such component a 
\emph{$(j_\pm,s,\omega)$-compartment
$\HC(j_\pm,s,\omega)$}, and we are going to investigate such compartments
with energy methods.

\subsection{Desingularized knot energies}\label{sec:1.2}

A variational approach to investigate knotted space 
curves is to 
consider knot energies as proposed by Fukuhara \cite{Fukuhara}, 
and the 
first examples of such energies were introduced by O'Hara in 
\cite{ohara}. Typically, knot energies model self-repulsive behavior 
and therefore provide infinite energy barriers between different knot 
classes since any curve with self-intersections has infinite energy. 
However, planar curves exhibit self-intersections unless they can be 
deformed into the circle. Hence, the value of any knot energy would be 
infinite on curves with self-intersections, which is why a modification 
of these energies is needed. Dunning investigated in 
\cite{dunning_t,dunning_paper} minimizers of a truncated version of the 
M\"obius energy\footnote{The Möbius energy was originally introduced 
by O'Hara in \cite{ohara} and it owes its name to 
its invariance under
M\"obius transformations \cite[Theorem 2.1]{freedman-etal_1994}.} in 
the class of figure-eight shaped curves. Later on, Kube \cite{dennis} 
closed several 
gaps in the work of Dunning and identified a
renormalized version of the 
M\"obius energy 
to characterize an optimally
immersed figure-eight curve in terms 
of its single intersection angle. 

We study a modified version of the so-called 
\emph{tangent-point energy}
\begin{align}\label{eq:TP}
    \TP_q(\gamma)
    &:=\iint\limits_{(\Sc)^2} 
    \frac{|\gamma'(s)||\gamma'(t)|}{\left(r_{\text{tp}}[\gamma](\gamma(t),\gamma(s))\right)^q}
ds dt
    =\iint\limits_{(\Sc)^2} \left(\frac{2\; 
    \dist(\ell(t),\gamma(s))}{|\gamma(t)-\gamma(s)|^2}\right)^q 
    |\gamma'(s)||\gamma'(t)| ds dt
\end{align}
for $q>2$ and regular
closed curves $\gamma \in C^1(\Sc,\R^d)$, $d\ge 2$. 
This energy was  first investigated analytically
by P. Strzelecki and H. 
von der Mosel in \cite{vdm_tangent}. 
Here, $r_{\text{tp}}[\gamma]\left(\gamma(t),\gamma(s)\right)$ is the 
radius of the unique circle passing through $\gamma(s)$ and $\gamma(t)$ 
that 
is tangent to the curve at $\gamma(t)$, and this radius
may be expressed as in \eqref{eq:TP} using
the affine tangent line $\ell(t):=\g(t)+\R\g'(t)$
of $\g$ at the curve point $\g(t)$.
S. Blatt  characterized 
the corresponding energy space and showed that a curve parametrized by 
arclength has finite tangent-point energy if and only if it 
is embedded and lies in the 
fractional Sobolev\footnote{For the definition and a brief account
on the basic properties of one-dimensional periodic fractional
Sobolev spaces it suffices to consult, e.g., 
\cite[Appendix A]{knappmann-etal_2022}.}
space $W^{2-\frac{1}{q},q}(\Sc,\R^d)$; 
see \cite[Theorem 1.1]{blatt_energy_space} or 
\cite[Theorem 1.1]{blatt_reiter_regularity}. 
This characterization 
opened the way towards critical point theory, either via
symmetric criticality (see \cite{gilsbach_2018} for $\TP_q$ alone,
or \cite{gilsbach-etal_2023} where the tangent-point energy is
added as a self-avoidance term to the bending energy), or by
means of (Banach-)gradient flows in \cite{matt-etal_2023}. For variants
of the tangent-point energy whose underlying energy space is Hilbert,
very recently the Palais-Smale condition was shown \cite{freches-etal_2025},
which leads to long-time existence and subconvergence to critical points of the
gradient flow. In \cite{doehrer-etal_2025} this variant of the tangent-point energy
was used to construct a complete Riemannian metric on knot space such that
every pair of knots in a knot class can be connected by a distance-minimizing
geodesic.

To prevent the blow-up of the
energy on curves \emph{with} self-intersections, 
we modify the integration domain and ``cut out'' the points in the 
preimage near the self-intersections. 
To make this more precise, let 
$d_\gamma(x,y)$ denote the intrinsic distance between the points 
$\gamma(x)$ and $\gamma(y)$ along the curve. 
Suppose that $\gamma \in  C^1\left(\Sc,\R^2\right)$ is
 a regular 
curve with finitely many self-intersections, i.e.\ the set
\begin{align}
  \label{def_S_gamma}
  S(\gamma):=\left\{(u,v) \in (\Sc)^2 : 0\le u<v<1 
  \text{ and } \gamma(u)=\gamma(v)  \right\}
\end{align}
is finite. Thus, $S(\gamma)=\{(u_1,v_1),\ldots,(u_n,v_n)\}$ 
for some $n \in \N$. Define furthermore
\begin{align}
  \label{def_T_gamma}
T(\gamma):=\left\{ u_1, v_1, \dots, u_n,v_n \right\}
\AND\Lambda(\gamma)
:=
\min \left\{d_\gamma(a,b) \; | \; a,b \in T(\gamma), a \neq b \right\},
\end{align}
so that $\Lambda(\gamma)
$ is the shortest intrinsic distance along the curve between any two 
self-intersections. For $\delta \in \left(0,\Lambda(\gamma)\right)$ 
we define the set
\begin{align*}
Y_\delta(\gamma):={}&\big\{(x,y)\in (\Sc)^2 \; : \; 
\bigl(d_\gamma(x,u_i)<\delta \text{ and } d_\gamma(y,v_i)<\delta 
\text{ for some } i=1,\dots,n \bigr) \\ 
&\text{ or } \bigl(d_\gamma(x,v_i)<\delta \text{ and } 
d_\gamma(y,u_i)<\delta \text{ for some } i=1,\dots,n \bigr) \big\}. 
\end{align*}
Sometimes, it is more suitable to consider tuples 
$(x,x+w) \in (\Sc)^2$ for $x \in \Sc$ and 
$w \in [ -\frac{1}{2},\frac{1}{2}]$, as this yields 
$|(x+w)-x |_{\Sc}=|w |_{\Sc}=|w |$. Here, $|\cdot|_{\Sc}$ denotes the periodic distance on $\Sc$ defined as 
\begin{align*}
|x-y |_{\Sc}:=\min_{k \in \Z} |x-y+k|. 
\end{align*}
Then we set
\begin{align}
  \hat{Y}_\delta(\gamma):={}&\big\{ (x,w)\in \Sc\times
  [-\tfrac{1}{2},\tfrac{1}{2}] \; : \; 
  (x,x+w) \in Y_\delta(\gamma) \big\}. \label{def_Yhatdelta}
  \end{align}
Define for $0 < \delta < \frac{\Lambda(\gamma)}2$ the \emph{truncated 
tangent-point energy} as 
\begin{align*}
\textstyle
\TP_{q,\delta}(\gamma):=\iint\limits_{\Ydc} 
\frac{
|\gamma'(s)||\gamma'(t)|}{
\left(r_{\text{tp}}[\gamma](\gamma(t),\gamma(s))\right)^q}ds dt
=\iint\limits_{\Sc\times [-\frac12,\frac12]\setminus\hat{Y}_\delta(\g)}
\frac{
|\gamma'(x)||\gamma'(x+w)|}{
\left(r_{\text{tp}}[\gamma](\gamma(x),\gamma(x+w))\right)^q}dwdx
.
\end{align*}
We mostly restrict $\TP_{q,\delta}$ to the following  
class of admissible curves with fixed winding number and 
Arnold invariants, that are affine linear near their self-intersections. 
\begin{definition}
  \label{admissibile_curves}
  For $\eta \in (0,1)$ and  $j_+,j_-,s,\omega \in \Z$, let
 $\HF(\eta,j_\pm,s,\omega)$ be the set of all arclength
 parametrized generic immersions $\g\in C^1(\R/\Z,\R^2)$
 with  $J^+(\gamma)=j_+$, $J^-(\gamma)=j_-$,  $\st(\gamma)=s$,  
$W(\gamma)=\omega$, and $\Lambda(\g)\ge 2\eta$, 
 such that $\gamma|_{B_\eta(u)}$ is
  affine linear for all $u\in T(\g)$.
\end{definition}
Restricting to arclength parametrizations (and therefore to
unit length curves) does not change or reduce the topological
information, so there is no loss of generality in doing so.
Notice that even if a suitable combination of integers
$j_+,$ $j_-,$ $s$, $\omega$ permits a non-empty 
$(j_\pm,s,\omega)$-compartment
$\HC(j_\pm,s,\omega)$ of generic immersions,
the class $\HF(\eta,j_\pm,s,\omega)$ is
empty if $2(j_+-j_-)\eta >1$ by means of
\eqref{connection_AI_number_self_intersections},
since every self-intersection consumes
$2\eta$ of the curve's unit length. We will show in 
Theorem \ref{thm:admissible-curves}, however,
that in any compartment $\HC(j_\pm,s,\omega)$
we find for sufficiently small $\eta>0$ a generic
$C^{1,1}$-immersion contained in $\HF(\eta,j_\pm,s,\omega)$.
By definition, the sets $\HF(\eta,j_\pm,s,\omega)$ are nested
with respect to the parameter $\eta$, i.e., 
\begin{equation}\label{eq:nested-F}
\HF(\eta_2,j_\pm,s,\omega)\subset\HF(\eta_1,j_\pm,s,\omega)\quad
\Foa 0<\eta_1\le\eta_2.
\end{equation}

\subsection{Main results}\label{sec:1.3}

We now state the main results of this paper, 
starting with the  existence of admissible $C^{1,1}$-curves in  
$\HF(\eta,j_\pm,s,\omega)$.

\begin{theorem}[Existence of admissible curves]
\label{thm:admissible-curves}
For every  $(j_\pm,s,\omega)$-com\-part\-ment 
$\mathcal{C}=\HC(j_\pm,s,\omega)$ 
there exists $\eta_0=\eta_0(\HC)>0$ such that
$\HF(\eta,j_\pm,s,\omega)\cap\HC\cap C^{1,1}(\R/\Z,\R^2)\not=\emptyset$ 
for all $\eta\in (0,\eta_0]$.
In addition, there is some $\eta_1=\eta_1(\HC)\in (0,\eta_0(\HC))$ such
that for every $\eta\in (0,\eta_1]$ there is a 
curve in $\HF(\eta,j_\pm,s,\omega)\cap\HC$ of class $C^{1,1}(\R/\Z,\R^2)$
 for which 
the intersection angle at every self-intersection equals $\frac{\pi}2$.
\end{theorem}
From now on  and throughout the paper
we use the positive number $\eta_0(\HC)$ for a  
$(j_\pm,s,\omega)$-compartment $\HC:=\HC(j_\pm,s,\omega)$ to formulate
our results on the non-empty admissible classes  $\HF(\eta,j_\pm,s,\omega)$
for $\eta\in (0,\eta_0]$.
\begin{theorem}[Existence of minimizers]
  \label{ex_min}
  For all  com\-part\-ments $\HC=\HC(j_\pm,s,\omega)$ and any
  $\eta\in (0,\eta_0(\HC)],$ $\delta\in (0,\frac{\eta}2]$ and $q>2$ there exists
  an immersion $\gamma_\delta^\eta \in \HF(\eta,j_\pm,s,\omega)
  \cap W^{2-\frac1{q},q}(\R/\Z,\R^2)\cap\HC$ 
  such that
  \begin{align*}
  \TP_{q,\delta}\left(\gamma_\delta^\eta \right)\leq \TP_{q,\delta}(\gamma) \text{ for all } \gamma \in \HF(\eta,j_\pm,s,\omega)\cap\HC.
  \end{align*}
   \end{theorem}
  For small truncation parameters $\delta$ close to the threshold
  $\frac{\eta}2$ we have computed $\TP_{q,\delta}$-minimizing 
  configurations numerically for various sets of prescribed winding
  numbers and Arnold invariants; see Figure \ref{fig:9pictures}.  As
  a predominant feature we observe in these numerical minimizers
  right angles at every self-intersection. 
 To support this numerical evidence analytically we address the 
 question whether it is possible to send the 
  truncation parameter $\delta$ to zero so that
  the truncated energies $\TP_{q,\delta}$ see an increasingly
  larger portion of a curve $\g$. But the energy values
 $\TP_{q,\delta}(\g)$ would tend to infinity if $\g$ self-intersects.
To take this blow-up into account, we scale the energy by
  the correct blow-up rate $\delta^{q-2}$, which will then
  allow for the limiting process. It turns out, that 
  these scaled versions of the energies 
indeed Gamma-converge towards the
  \emph{renormalized tangent-point energy $R_q:=
  \lim_{\delta\to 0}\delta^{q-2}\TP_{q,\delta}$},
  a functional which  depends only
  on the intersection angles of the curves.

  \begin{theorem}[Gamma convergence as $\delta\to 0$]
    \label{gamma_convergence}
    For all   com\-part\-ments $\HC=\HC(j_\pm, s,\omega)$ and 
    any $\eta\in (0,\eta_0(\HC)]$ and $q>2$ one has
        \begin{align*}
\delta^{q-2} \TP_{q,\delta} \overset{\Gamma}{\longrightarrow} R_q 
\quad
       \text{as }\delta \to 0 
       \quad \text{on }\big(\HF(\eta,j_\pm,s,\omega)\cap \HC\cap
 W^{2-\frac{1}{q},q}\left(\R/\Z,\R^2\right),||\cdot ||_{C^1} \big).
    \end{align*}
The renormalized tangent-point energy $R_q$ depends only on 
    the intersection angles of the curves and is minimized if and only if every intersection angle is a right angle.
    \end{theorem}

Minimizing the renormalized energy $R_q$ itself would lead to minimizing curves
of arbitrarily complicated shapes as long as all self-intersection angles equal $\frac{\pi}2.$ In search of representatives with an optimal shape we prove in addition
to Gamma convergence,  
the convergence of a sequence of  $\TP_{q,\delta}$-minimizers 
$\g_\delta^{\eta}$ to a limit curve $\Gamma^\eta$ as $\delta\to 0$. This limit
$\Gamma^\eta$ has exclusively right self-intersection angles and
can be seen as an optimal curve in the class $\HF(\eta,j_\pm,s,\omega)\cap\HC$, because
 $\Gamma^\eta$ is an \emph{almost-minimizer} of all truncated energies
    for sufficiently small truncation parameters $\delta$. 
    To derive the necessary a priori estimates independent of $\delta$ we use
    comparison curves whose intersection angles equal $\frac{\pi}2$. That is why
    we use the threshold parameter $\eta_1(\HC)$ from Theorem 
    \ref{thm:admissible-curves} instead of 
    $\eta_0(\HC)$ in the formulation of the following results.

  \begin{theorem}[Limit immersion is almost-minimizer]
    \label{limit_delta}
    Let $\HC=\HC(j_\pm,s,\omega)$ be a  com\-part\-ment, $q>2$ and
    $\eta\in (0,\eta_1(\HC)]$. 
    Then every sequence $(\g_\delta^\eta)_\delta
    \subset\HF(\eta,j_\pm,s,\omega)\cap\HC$ 
    of $\TP_{q,\delta}$-minimizers  subconverges (after suitable 
    translations) 
    in $C^1$ to some curve $\G^\eta\in\HF(\eta,j_\pm,s,\omega)
    \cap\HC\cap W^{2-\frac1{q},q}(\R/\Z,\R ^2)$ as $\delta\to 0$. 
    Moreover, 
    all intersection angles of $\G^\eta$ equal $\frac{\pi}2,$
    and for any $\epsilon >0$ there is $\hat{\delta}=\hat{\delta}(\epsilon)>0$
    such that
    \begin{equation}\label{eq:almost-minimal}
    \inf_{\HF(\eta,j_\pm,s,\omega)\cap\HC}\TP_{q,\delta}(\cdot)
    \le\TP_{q,\delta}(\G^\eta)<\inf_{\HF(\eta,j_\pm,s,\omega)\cap\HC}
    \TP_{q,\delta}(\cdot) +\epsilon\quad\Foa \delta\in (0,\hat{\delta}).
    \end{equation}
     \end{theorem}

   Notice that the minimal energy values on the left-hand side of
   \eqref{eq:almost-minimal} blow up as $\delta\to 0$, but the
   energy values $\TP_{q,\delta}(\G^\eta)$ remain in the  fixed
   $\epsilon$-neighborhood of these minimal energies as $\delta\to 0$.
   The proof of \cref{limit_delta} reveals even  more energetic control: 
   If $(\g_{\delta_k}^\eta)_k$ is a subsequence 
   of (suitably translated) $\TP_{q,\delta_k}$-minimizers,
   converging in $C^1$ to $\G^\eta$ for $\delta_k\to 0$ as
   $k\to\infty$, then
   \begin{equation}\label{eq:limit-energy-value}
R_q(\G^\eta)=\lim_{k\to\infty}\delta_k^{q-2}
\TP_{q,\delta_k}(\g_{\delta_k}^\eta).
   \end{equation}

Finally,
it turns out that this optimal curve is a true minimizer for
\emph{all truncated energies} among all admissible curves whose intersection angles
are right angles.
\begin{corollary}[Optimal immersion minimizes among curves with 
intersection angles $\frac{\pi}2$]\label{cor:optimal-among-right-angels}
Under the assumptions of Theorem \ref{limit_delta} one has 
\begin{equation}\label{eq:minimal-right-angles}
\textstyle
\TP_{q,\delta}(\G^\eta)\le\TP_{q,\delta}(\g)\quad\Foa \delta\in (0,\frac{\eta}2],\, \g\in\HF(\eta,j_\pm,s,\omega)\cap\HC,
\end{equation}
if all intersection angles of $\g$ equal $\frac{\pi}2.$
\end{corollary}

{\bf Remarks.}\,
1.\,
At this point it is not clear if the limit curves $\G^\eta$ also minimize the truncated
tangent-point energy in the full admissibility class $\HF(\eta,j_\pm,s,\omega)\cap\HC$.
By Arzela-Ascoli the $\G^\eta$ subconverge uniformly to a Lipschitz continuous
limit curve $\G^0$ as $\eta\to 0$, but we do not have any further information about
$\G^0$, in spite of the fact that the a priori energy bound on the minimizers $\g_\delta^\eta$
in \cref{energy_uniform_bound} in Section \ref{sending_delta_to_zero}
does \emph{not} depend on $\eta$. In addition,
the map $\eta\mapsto\TP_{q,\theta}(\G^\eta)$ is non-decreasing on the interval
$(2\theta,\eta_1(\HC))$; see \cref{cor:eta-monoton}, but uniform convergence
does not suffice to prove lower semi-continuity of $\TP_{q,\theta}$. We do expect,
however, that right intersection angles and straight segments near the 
self-intersections
will remain prevalent in the limit $\eta\to 0$ since this local geometry contributes
least to the energy. Therefore, we believe that the numerically computed minimizers
in Figure \ref{fig:9pictures} are fairly close to the true minimizing 
immersions even 
without the restriction of being locally affine linear near the self-intersections.

2.\,
There is an interesting connection to Legendrian knots in $\R ^3$ equipped
with the standard contact structure; see, e.g., \cite[Chapter 3]{geiges_2008}. The
Lagrangian projection, i.e., the projection onto the $xy$-plane of a Legendrian knot
yields a planar immersed curve. On the other hand, any planar immersion that encloses
zero area can be lifted to a closed Legendrian knot that is unique up to translation
in the $z$-direction. An example of such a planar immersion is a point-symmetric
figure-eight shaped curve with winding number $W=0$ and Arnold invariants
$J^+=0$, $J^-=-1$, and $\st=0$. In \cite[Chapter 7]{lagemann_diss} it was shown
by means of Palais's principle of symmetric criticality
that such curves arise as critical points of the energy $\TP_{q,\delta}$. A natural 
question to ask is whether lifting these symmetric critical points yields special
Legendrian knots. At this point it is open, however, if \emph{all} figure-eight shaped
minimizers $\g_\delta^\eta$  obtained in \cref{ex_min}
have this point symmetry.

The paper is structured as follows. In Section \ref{admissible_curves}
we construct admissible curves in the class $\HF(\eta,j_\pm,s,\omega)$ with
sufficient regularity which will prove  \cref{thm:admissible-curves}. 
Section \ref{def_delta_energy} is devoted to the analysis of the truncated
tangent-point energy, its regularizing effects (\cref{regularity})
and its uniform control on bilipschitz
constants; see \cref{bilipschitz}.  The proof of the existence result 
(\cref{ex_min}) in Section \ref{chapter_existence_minimizers} rests on 
lower semi-continuity (\cref{lower_sc}) and 
compactness  established in \cref{compactness}. To investigate the limit $\delta\to 0$
we introduce in \cref{def:annular} the annular truncated tangent-point
energy for which we prove  in \cref{gamma_c_min} a
 crucial representation that leads to the 
purely angle-dependent Gamma limit $R_q$ proving \cref{gamma_convergence}.  
The results on the limit curve $\G^\eta$ require a priori estimates on 
the minimal energies $\TP_{q,\delta}(\g_\delta^\eta)$ established
in \cref{energy_uniform_bound}. The appendix provides the details for
a specific reparametrization of the immersions with local graph patches
near self-intersections (\cref{lem:bauser}), and the analysis of the angle-dependent
term 
in the annular energy of \cref{def:annular} in \cref{angle_energy_monotone}.
Throughout the paper we use the notation $B_r(x)$ for open subintervals
$(x-r,x+r)\subset\R/\Z$.

\section{\texorpdfstring{Existence of admissible curves}{Existence of
admissible curves}}\label{admissible_curves}
We devote this preparatory section to prove the existence of admissible curves in a given
 compartment $\mathcal{C}(j_\pm,s,\omega)$.

{\it Proof of Theorem \ref{thm:admissible-curves}.}\,
The idea is to deform a suitably reparametrized curve in 
$\HC(j_\pm,s,\omega)$ locally near every of its finitely many transverse 
double points to a pair of graphs which intersect in straight line segments. 
To obtain  right angles at these self-intersections we rotate the inner quarter
of one segment appropriately with a smooth cut-off to leave the rest of the curve
unchanged. 
Here are the details presented in four steps.

{\it Step 1.}\,
Any curve in the open compartment
$\HC=\HC(j_\pm,s,\omega)$ can be smoothened, e.g. 
by convolution, since the winding number and Arnold invariants
are locally constant. In addition, such a smooth representative
can be rescaled to unit length and 
reparametrized to arclength to 
obtain a generic immersion $\g\in C^\infty(\R/\Z,\R^2)
\cap\HC$ with winding number $W(\g)=\omega$ 
and Arnold invariants $J^+(\g)=j_+$,
$J^-(\g)=j_-$, $\st(\g)=s$.  
So, $|\g'|
\equiv 1$ on $\R/\Z$, and $\g$ has exactly $n:=j_+-j_-$ transverse
self-intersections, each of multiplicity two, and no other self-intersections. In other words,
there are mutually distinct arclength parameters $s_i,t_i\in\R/\Z$ with $\g(s_i)=
\g(t_i)$ for $i=1,\ldots,n$.  Abbreviate the unit tangents at these intersection points
by $S_i:=\g'(s_i)$ and $T_i:=\g'(t_i)$ and define
\begin{equation}\label{eq:d}
d:=\min_{i\in\{1,\ldots,n\}}\big\{ |S_i-T_i|  \big\}.
\end{equation}
Note that $d\in (0,2)$ by transversality.

By virtue of Lemma \ref{lem:bauser} for $k:=2$
 we can reparametrize $\g$ preserving the
orientation to obtain 
$\tilde{\g}\in C^2(\R/\Z,\R^2)$ 
which locally near every self-intersection coincides with a
$C^2$-graph over the affine tangent plane 
$\g(s_i)+\R S_i$
and $\g(t_i)+\R T_i$, respectively. More precisely, $\tilde{\g}(\R/\Z)=
\g(\R/\Z)$, $|\tilde{\g}'|\ge \frac14$ on $\R/\Z$, and there is some $r_0>0$ such that
\begin{equation}\label{eq:tilde-gamma}
\tilde{\g}(x)=\begin{cases}
\g(x) & \Fo x\not\in\bigcup_{i=1}^n B_{2r_0}(s_i)\cup B_{2r_0}(t_i)\\
\g(s_i)+(x-s_i)S_i+u_i(x-s_i)S_i^\perp & \Fo
x\in B_{r_0}(s_i)\\
\g(t_i)+(x-t_i)S_i+v_i(x-t_i)T_i^\perp & \Fo
x\in B_{r_0}(t_i).
\end{cases}
\end{equation}
Here, the 
subintervals 
$B_{3r_0}(s_i), B_{3r_0}(t_i)\subset\R/\Z$, $i=1,\ldots,n$, are mutually 
disjoint, and the intrinsic distance $d_{\tilde{\g}}(x,s_i)\ge r_0$ for all $x\not\in B_{r_0}(s_i)$, so that
\begin{equation}\label{eq:Lambda-issue}
\Lambda(\tilde{\g})\ge 2r_0,
\end{equation}
and 
the graph functions $u_i,v_i\in C^2(\R)$  
satisfy $0=u_i(0)=v_i(0)=u_i'(0)=v_i'(0)$ for all $i=1,\ldots,n$.
In addition, we have used in \eqref{eq:tilde-gamma}
the notation $p^\perp:=(-p_2,p_1)^T$ for a vector
$p=(p_1,p_2)^T\in\R^2.$
By continuity we can choose $r_1\in (0,r_0)$ so small that
\begin{equation}\label{eq:u_i-bound}
\max\big\{\|u_i'\|_{C^0(B_{r_1}(0))}, \|v_i'\|_{C^0(B_{r_1}(0))}\big\} 
< \frac{d}{40}\quad\Foa
i=1,\ldots,n.
\end{equation}

\begin{figure}[b]
\centering
\begin{subfigure}{0.45\textwidth}
\includegraphics[width=\textwidth,angle=0,scale=0.8]{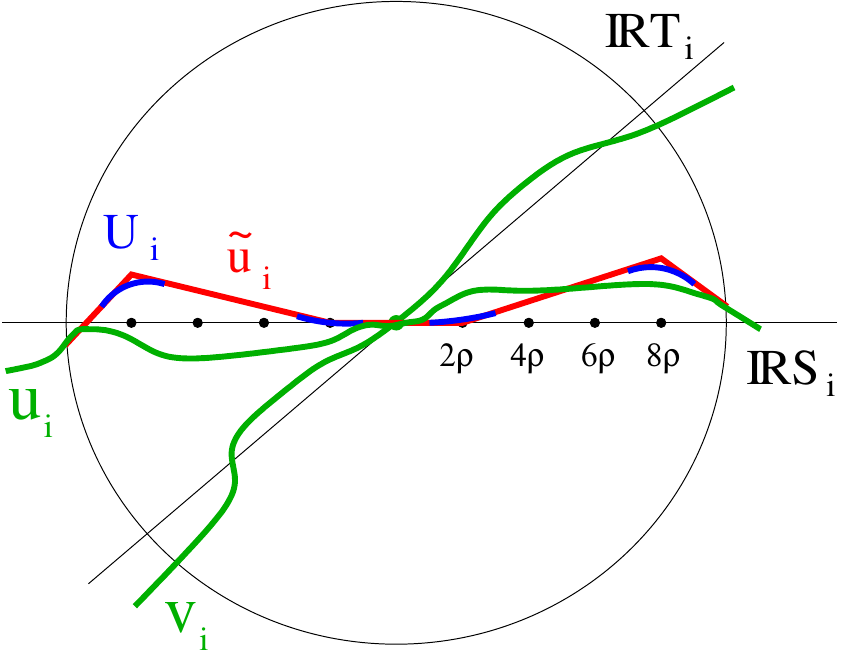}

\caption{
Step 2: The auxiliary piecewise linear function $\tilde{u}_i$ (in red)
interpolates the
point-tangent data of $u_i$ at $\partial B_{10\rho}(0)$ and
$\partial B_{2\rho}(0)$. Smoothing out the corners of $\tilde{u}_i$
with circular arcs (blue) yields the $C^{1,1}$-comparison function 
$U_i$.}
\end{subfigure}
\hspace{1cm}
\begin{subfigure}{0.45\textwidth}
\includegraphics[width=\textwidth,angle=0,scale=0.8]{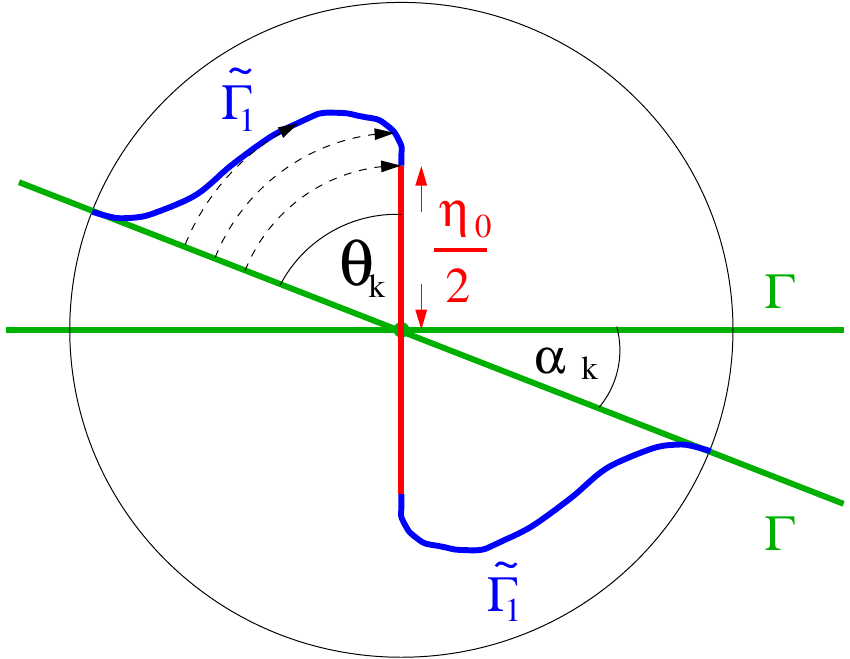}
  
\caption{
Step 4: Rotating one of the two intersecting linear
segments to a perpendicular position (in 
red) and smoothly cutting off the rotational angle (blue) away from
the self-intersection yields
a comparison curve $\tilde{\G}_1$ whose intersection angles 
equal $\frac{\pi}2$.
}
\end{subfigure}

\caption{\label{fig:construct-comparison-curve-step}}
\end{figure}

{\it Step 2.}\,
The quantity
\begin{equation}\label{eq:D}
D:=\min_{i\in\{1,\ldots,n\}}\min_{y\in(\R/\Z)\setminus [B_{r_1}(s_i)\cup B_{r_1}(t_i)]}
\big\{ |\g(s_i)-\tilde{\g}(y)|\big\} 
\end{equation}
is strictly positive, since the reparametrization $\tilde{\g}$ has exactly the same
self-intersection points as $\g$, 
$$
\tilde{\g}(s_i)=\g(s_i)=\g(t_i)=\tilde{\g}(t_i)\quad\Fo i=1,\ldots, n,
$$
and not any other self-intersections.
Now choose a smaller radius 
\begin{equation}\label{eq:rho}
\textstyle\rho\in \big(0, \min\big\{\frac{r_1}{10}, \frac{D}{40} \big\}\big),
\end{equation}
and define for each $i=1,\ldots,n$ the piecewise linear auxiliary function
$\tilde{u}_i:B_{r_1}(0)\to\R$ as 
$$
\tilde{u}_i(x):=\begin{cases}
0 & \Fo x\in \overline{B_{2\rho}(0)}\\
u_i(-10\rho)+(x+10\rho)u_i'(-10\rho) & \Fo x\in [-10\rho,-8\rho]\\
u_i(10\rho)+(x-10\rho)u_i'(10\rho) & \Fo x\in [8\rho,10\rho]\\
\textnormal{linearly interpolating} & \textnormal{\,\,\,on $B_{8\rho}(0)\setminus
\overline{B_{2\rho}(0)}$}\\
u_i(x) & \Fo x\in B_{r_1}(0)\setminus B_{10\rho}(0);
\end{cases}\notag
$$
see Figure \ref{fig:construct-comparison-curve-step}(a).

The slope of $
\tilde{u}_i$ on $[-10\rho,-8\rho]$ and on $[8\rho,10\rho]$ is
bounded by the norm $\|u_i'\|_{C^0(B_{r_1}(0))},$ 
whereas on $B_{8\rho}(0)\setminus
B_{2\rho}(0)$ it is controlled by
$$
\textstyle\left|\frac{\tilde{u}_i(\pm 8\rho)}{6\rho}\right|=\left|
\frac{u_i(\pm 10\rho)\mp 2\rho u_i'(\pm 10\rho)}{6\rho}\right|
\le \frac43\|u_i'\|_{C^0(B_{r_1}(0))}+\frac13\|u_i'\|_{C^0(B_{r_1}(0))}<2\|u_i'\|_{C^0(B_{r_1}(0))}.
$$
Replacing the graph of $\tilde{u}_i$ on the intervals $[-9\rho,-7\rho]$, 
$[-3\rho,-\rho]$, $[\rho,3\rho]$, and on $[7\rho,9\rho]$
by the respective $C^1$-interpolating circular arcs (see Figure 
\ref{fig:construct-comparison-curve-step}(a))
we obtain a piecewise linear and circular function $U_i\in 
C^{1,1}(B_{r_1}(0))$ 
vanishing
on $\overline{B_\rho(0)}$ and coinciding with $u_i$ on $B_{r_1}(0)\setminus
\overline{B_{10\rho}(0)},$ and satisfying
\begin{equation}\label{eq:U-i}
\|U_i'\|_{C^0(B_{r_1}(0))} < 2 \|u_i'\|_{C^0(B_{r_1}(0))}\quad\Fo i=1,\ldots,n.
\end{equation}
Analogously, we obtain a piecewise linear and circular function $V_i\in
C^{1,1}(B_{r_1}(0))$ vanishing on $\overline{B_\rho(0)}$ 
and coinciding with $v_i$ on $B_{r_1}(0)\setminus
\overline{B_{10\rho}(0)},$ such that
\begin{equation}\label{eq:V-i}
\|V_i'\|_{C^0(B_{r_1}(0))} < 2 \|v_i'\|_{C^0(B_{r_1}(0))}\quad\Fo i=1,\ldots,n.
\end{equation}

{\it Step 3.}\,
Consider the convex combinations $u_{i,\lambda}:=(1-\lambda)u_i+\lambda U_i$
on $B_{r_1}(0)$ and $v_{i,\lambda}:=(1-\lambda)v_i+\lambda V_i$
on $B_{r_1}(0)$ for $\lambda\in [0,1]$, $i=1,\ldots,n,$ 
satisfying 
\begin{equation}\label{eq:u-i-lambda-est}
\max\big\{\|u_{i,\lambda}'\|_{C^0(B_{r_1}(0))}, \|v_{i,\lambda}'\|_{C^0(B_{r_1}(0))}  \big\} < 3 \|u_i'\|_{C^0(B_{r_1}(0))} < \frac{3d}{40}
\end{equation}
by means of \eqref{eq:U-i} and \eqref{eq:V-i} in combination with  \eqref{eq:u_i-bound}.
Now define the $1$-parameter family of curves 
$$
\tilde{\g}_\lambda(x):=\begin{cases}
\tilde{\g}(x) & \Fo x\not\in\bigcup_{i=1}^n [B_{r_1}(s_i)\cup B_{r_1}(t_i)]\\
\g(s_i)+(x-s_i)S_i+u_{i,\lambda}(x-s_i)S_i^\perp & \Fo x\in B_{r_1}(s_i)\\
\g(s_i)+(x-t_i)T_i+v_{i,\lambda}(x-t_i)T_i^\perp & \Fo x\in B_{r_1}(t_i).
\end{cases}
$$
{\it Claim: For each $\lambda\in [0,1]$ the curve $\tilde{\g}_\lambda$ has exactly
the same self-intersection points as $\g$,
\begin{equation}\label{eq:intersection2}
\tilde{\g}_\lambda(s_i)=\g(s_i)=\g(t_i)=\tilde{\g}_\lambda(t_i)\quad\Fo i=1,\ldots,n,\notag
\end{equation}
and no other self-intersections.}

Indeed, with \eqref{eq:u-i-lambda-est} and \eqref{eq:rho} we obtain for
any $\lambda\in [0,1]$ and $x\in \overline{B_{10\rho}(s_i)}$
\begin{align*}
|\tilde{\g}_\lambda(x)-\g(s_i)|^2 & = (x-s_i)^2+u_{i,\lambda}^2(x-s_i)\le
(1+\|u_{i,\lambda}'\|_{C^0(B_{r_1}(0))}^2) (10\rho)^2\\
&\overset{\eqref{eq:u-i-lambda-est}}{\le} (1+\frac{9d^2}{1600}) 100\rho^2 <
200\rho^2\overset{\eqref{eq:rho}}{<} 200 \frac{D^2}{40^2}<\frac{D^2}4,
\end{align*}
which implies for every $\lambda\in [0,1]$ by definition of the quantity $D$ in \eqref{eq:D}
\begin{align}\label{eq:dist-global}
|\tilde{\g}_\lambda(x)-\tilde{\g}_\lambda(y)| > \frac{D}2\Foa
x\in \overline{B_{10\rho}(x_i)},\, y\not\in B_{r_1}(s_i)\cup B_{r_1}(t_i).
\notag
\end{align}
For $x\in \overline{B_{10\rho}(s_i)}$ and $y\in B_{r_1}(s_i)\setminus\{x\}$ one has
\begin{equation}\label{eq:dist-same-graph}
|\tilde{\g}_\lambda(x)-\tilde{\g}_\lambda(y)|^2=
(x-y)^2+(u_{i,\lambda}(x-s_i)-u_{i,\lambda}(y-s_i))^2\ge (x-y)^2>0.\notag
\end{equation}
Finally, for $x\in \overline{B_{10\rho}(s_i)}$ and $y\in B_{r_1}(t_i)$ we first estimate the distance 
between the respective tangential projections of $\tilde{\g}_\lambda(x)$ and
$\tilde{\g}_\lambda(y)$, i.e.,
\begin{align}
f_i(x,y) &:= |\g(s_i)+(x-s_i)S_i-[\g(s_i)+(y-t_i)T_i]|\notag\\
&\ge
\max\big\{  \dist((x-s_i)S_i,\R T_i),\dist((y-t_i)T_i,\R S_i) \big\}\notag \\
&= \max\big\{  |x-s_i|,|y-t_i| \big\}\sin\ANG (S_i,T_i).\label{eq:tangential-dist}
\end{align}
With $\sin\ANG(S_i,T_i)>\frac12 \sin\ANG(S_i,T_i)=\frac12 |S_i-T_i|\ge\frac{d}2 >0$
we find from \eqref{eq:tangential-dist} 
\begin{equation}\label{eq:tan-dist}
f_i(x,y)\ge\max\big\{ |x-s_i,|y-t_i| \big\}\cdot\frac{d}2\quad\Foa x\in
\overline{B_{10\rho}(s_i)},\, y\in B_{r_1}(t_i),\,i=1,\ldots,n.
\end{equation} 
On the other hand, by \eqref{eq:u-i-lambda-est}
$$
|\tilde{\g}_\lambda(x)-(\g(s_i)+(x-s_i)S_i)|^2=u_{i,\lambda}^2(x-s_i)\le\frac{9d^2}{1600}
|x-s_i|^2
$$
and, likewise,
$$
|\tilde{\g}_\lambda(y)-(\g(s_i)+(y-t_i)T_i)|^2=v_{i,\lambda}^2(y-t_i)\le\frac{9d^2}{1600}
|y-t_i|^2,
$$
so that we obtain with \eqref{eq:tan-dist}
\begin{equation}\label{eq:dist-semilocal}
|\tilde{\g}_\lambda(x)-\tilde{\g}_\lambda(y)|\ge
\max\big\{  |x-s_i|,|y-t_i| \big\}\cdot\frac{7}{20}d\quad\Fo
x\in\overline{B_{10\rho}(s_i)},\, y\in B_{r_1}(t_i), \notag
\end{equation}
which proves the claim.

In addition, we observe that the map $\lambda\mapsto H(\cdot,\lambda):=\tilde{\g}_\lambda(\cdot)$ is of class $C^0([0,1],C^1(\R/\Z,\R^2))$,
and
$$
\tilde{\g}_\lambda'(x)=\begin{cases}
\tilde{\g}'(x) & \Fo x\not\in\bigcup_{i=1}^n [B_{r_1}(s_i)\cup B_{r_1}(t_i)]\\
S_i+u_{i,\lambda}'(x-s_i)S_i^\perp & \Fo x\in B_{r_1}(s_i)\\
T_i+v_{i,\lambda}'(x-t_i)T_i^\perp & \Fo x\in B_{r_1}(t_i),
\end{cases}
$$
so that $|\tilde{\g}_\lambda'(x)|\ge \frac14 $ for all $x\in\R/\Z.$ 
Consequently, $H$ is a regular homotopy between $H(\cdot,0)=\tilde{\g}(\cdot)$
and $H(\cdot,1)=\tilde{\g}_1\in\HC $ where
$\tilde{\g}_1|_{B_\rho(s_i)}$ and
$\tilde{\g}_1|_{B_\rho(t_i)}$ are affine linear,
by definition of $u_{i,1}=U_i$
and $v_{i,1}=V_i$ for $i=1,\ldots,n,$ 
so that $\Lambda(\tilde{\g}_1)>20\rho$ by
means of \eqref{eq:Lambda-issue}
and \eqref{eq:rho}.
Rescaling $\tilde{\g}_1$ by its length $L:=\mathscr{L}(\tilde{\g}_1)$
and then reparametrizing $\tilde{\g}_1/L$
to arclength yields the desired curve $\G\in \HC\cap 
C^{1,1}(\R/\Z,\R^2)$ 
which is
linear within positive arclength $\eta_0(\HC)=\eta_0:=
\frac{\rho}{L}
$ around each 
self-intersection parameter, satisfying $\Lambda(\G)>20\eta_0>2\eta_0$,
as required in the admissibility class $\HF(\eta_0,j_\pm,s,\omega)
\subset\HF(\eta,j_\pm,s,\omega)$ for all $\eta\in (0,\eta_0]$ 
by
means of \eqref{eq:nested-F}.

{\it Step 4.}\,
Let $\G(\sigma_k)=\G(\tau_k)$ for mutually distinct arclength parameters
$\sigma_k,\tau_k\in\R/\Z$, $k=1,\ldots,n$ be the $n=j_+-j_-$ double points of the 
curve $\G\in\HC $ obtained in the previous step. Notice that our construction
implies that $\G$ intersects the $\eta_0$-ball around each self-intersection
exactly in the two self-intersecting straight line segments, i.e.,
\begin{equation}\label{eq:eta-ball-intersection}
\G(\R/\Z)\cap B_{\eta_0}(\G(\sigma_k))=\G(B_{\eta_0}(\sigma_k))\cup
\G(B_{\eta_0}(\tau_k))\quad\Fo k=1,\ldots,n.
\end{equation}

and our choice of $\rho$ in \eqref{eq:rho}.
Denote the intersection angles of the tangent lines
by $\alpha_k:=\ANG (\R\G'(\sigma_k),\R\G'(\tau_k))\in (0,\frac{\pi}2]$
and set $\theta_k:=\frac{\pi}2-\alpha_k\in [0,\frac{\pi}2)$ for $k=1,\ldots,n.$ 
Identifying $(p_1,p_2)^T\in\R^2$ with points $p_1+ip_2$ in
the complex plane $\C$ we can rotate any point $\xi\in\C\simeq\R^2$ by
the angle $\epsilon_k\theta_k$ by multiplying $e^{i\epsilon_k  \theta_k}\in\C$ with $\xi$.
Here,
the sign $\epsilon_k\in\{1,-1\}$
is determined in such a way that the resulting angle
\begin{equation}\label{eq:epsilon-k}
\ANG \big(\R e^{i\epsilon_k  \theta_k}\G'(\sigma_k),\R\G'(\tau_k)\big)=\frac{\pi}2\quad
\Fo k=1,\ldots,n.
\end{equation}
Now choose a cut-off function $\phi_{\eta_0}\in C^\infty_0((-\eta_0,\eta_0))$
(extended to all of $\R$ by the value $0$) such that $0\le\phi_{\eta_0}(x)\le 1$
for all $x\in\R$, $\phi_{\eta_0}(x)=1$ for all $x\in [-\frac{\eta_0}2,
\frac{\eta_0}2]$. 
Define the deformed  curves for $\lambda\in [0,1]$ (see Figure
\ref{fig:construct-comparison-curve-step}(b) for $\lambda=1$)
$$
\tilde{\G}_\lambda(t):=\begin{cases}
e^{i\lambda\phi_{\eta_0}(t-\sigma_k)\epsilon_k\theta_k}(\G(t)-\G(\sigma_k))+\G(\sigma_k) &
\Fo t\in B_{\eta_0}(\sigma_k)\\
\G(t) & \Fo t\not\in\bigcup_{k=1}^n B_{\eta_0}(\sigma_k),
\end{cases}
$$
so that $\tilde{\G}_\lambda\in C^{1,1}(\R/\Z,\R^2)$ 
has for each $\lambda\in [0,1]$
exactly the same double points as $\G$, i.e.,
$$
\tilde{\G}_\lambda(\sigma_k)=\G(\sigma_k)=\G(\tau_k)=\tilde{\G}_\lambda(\tau_k)
\quad\Fo k=1,\ldots,n,
$$
and no additional self-intersections by means of \eqref{eq:eta-ball-intersection}.
Since $\G(t)-\G(\sigma_k)=\pm \G'(t)$ for $t\in B_{\eta_0}(\sigma_k)$ we
find for the  speed of $\tilde{\G}_\lambda $
\begin{align*}
|\tilde{\G}_\lambda'(t)| & = |e^{i\lambda\phi_{\eta_0}(t-\sigma_k)\epsilon_k\theta_k}
(1\pm i\lambda\epsilon_k
\theta_k\phi_{\eta_0}'(t-\sigma_k))\G'(t)|\ge |\G'(t)|=1
\Foa t\in B_{\eta_0}(\sigma_k),
\end{align*}
and $|\tilde{\G}_\lambda'(t)|=|\G'(t)|=1$
for $t\not\in\bigcup_{k=1}^nB_{\eta_0}(\sigma_k),$ so that the map
$\lambda\mapsto \tilde{H}(\cdot,\lambda):=\tilde{\G}_\lambda(\cdot)$
is a regular homotopy between $\tilde{H}(\cdot,0)=\G(\cdot)$ and
$\tilde{H}(\cdot,1)=\tilde{\G}_1$.
With $\tilde{\G}_1'(\sigma_k)=e^{i\epsilon_k\theta_k}\G'(\sigma_k)$ and
$\tilde{\G}_1'(\tau_k)=\G'(\tau_k)$ we obtain for the intersection angles by 
means of \eqref{eq:epsilon-k}
$$
\ANG\big(\R\tilde{\G}_1'(\sigma_k),\R\tilde{\G}_1'(\tau_k)\big)
=
\ANG\big(e^{i\epsilon_k\theta_k}\R\G'(\sigma_k),\R\G'(\tau_k)\big)=\frac{\pi}2\Foa k=
1,\ldots,n.
$$
Rescaling $\tilde{\G}_1$ to unit length and then reparametrizing to arclength yields
the desired curve of class $C^{1,1}(\R/\Z,\R^2)$ in
$\HC\cap\HF(\eta_1,j_\pm,s,\omega)$ for
some $\eta_1\in (0,\eta_0)$ such that the straight line segments at each
intersection point are perpendicular to each other. 
This shows again
by \eqref{eq:nested-F} that for all $\eta\in (0,\eta_1]$
the admissibility class $\HC\cap\HF(\eta,j_\pm,s,\omega)$ contains
curves with right  intersection angles.
\qed

\section{\texorpdfstring{Energy space of the $\delta$-renormalized tangent-point energy}{Energy space of the delta-renormalized tangent-point energy}}
\label{def_delta_energy}

In this section, we characterize the  energy space of the energy $\TP_{q,\delta}$
following the ideas of Blatt and Reiter in \cite[Section 2]{blatt_reiter_regularity}. Crucial for our
proofs are bilipschitz estimates for planar immersions away from their self-intersections.
Notice that  for any arclength parametrized
curve $\g\in C^1(\R/\Z,\R^d)$ we have
\begin{align}
&|\g(s +w) -\g(s)| 
 =\textstyle \Big|\int_s^{s+w}\g'(t^*)\,dt+\int_s^{s+w}\big(\g'(t)-\g'(t^*)
 \big)\,dt\Big|\label{eq:omega}\\
& \textstyle\ge |w|\big(1-\max_{t\in [s,s+w]}|\g'(t)-\g'(t^*)|\,\big)\,\,\Foa
(s,w)\in\R/\Z\times [-\frac12,\frac12],\,t^*\in [s,s+w].\notag
\end{align}
This inequality  will be used repeatedly throughout the paper to
obtain bilipschitz estimates, like in the following \emph{individual
global bilipschitz estimate} for planar 
immersed curves away from their self-intersections.

\begin{lemma}[Global individual bilipschitz estimate]
  \label{bilipschitz_1}
  Let $\gamma \in C^1\left(\R/\Z,\R^2\right)$ be an immersion with finitely 
  many self-intersections. Then for every $0<\delta<\frac{\Lambda(\gamma)}{2}$ there exists a constant $c=c(\gamma,\delta)>0$ such that
  \begin{align}
 |x-y|_{\Sc}\leq c(\g,\delta)|\gamma(x)-\gamma(y)| \quad \text{for all }(x,y) \in \Ydc. 
 \label{eq:global-individual-bilip}
  \end{align}
  \begin{proof}
 Since $\g'$ is uniformly continuous on $\R/\Z$ we find by means of 
 \eqref{eq:omega} a constant $\tau_\g>0$ such that
 \begin{equation}\label{eq:local-individual-bilip} 
  |\g(s+w)-\g(s)|\ge\frac12 |w|\quad\Foa (s,w)\in\R/\Z\times [-\tau_\g,\tau_\g].
  \end{equation}
  The function $(s,w)\mapsto |\g(s+w)-\g(s)|$ is positive and uniformly continuous
  on the compact set
  $
  K(\delta,\tau_\g):=\{(s,w)\in (\R/\Z\times [-\frac12,\frac12])\setminus \hat{Y}_\delta
  (\g):
  |w|\ge \tau_\g\}
  $
  so that there exists a constant $c_0(\g,\delta)>0$ such that
  \begin{equation}\label{eq:off-diagonal}
  |\g(s+w)-\g(s)|\ge c_0(\g,\delta)\quad\Foa (s,w)\in K(\delta,\tau_\g).\notag
  \end{equation}
  Combining this with \eqref{eq:local-individual-bilip} and the fact that
  $|w|\le\frac12$ yields  the desired global bilipschitz estimate 
  \eqref{eq:global-individual-bilip} with 
  $c(\g,\delta):=\max\{2,(2c_0(\g,\delta)^{-1}\},$  replacing $s$ by $x$ and
  $s+w$ by $y$.
  \end{proof}
  \end{lemma}

We prove most statements in this section in a larger class than the set 
$\HF(\eta,j_\pm,s,\omega)$ introduced in \cref{admissibile_curves}.

\begin{definition}\label{def:Fn}
Let $n \in \N$. Define $\HF_n \subset C^1\left(\Sc,\R^2\right)$ to be  the subset of
all arclength parametrized immersions $\gamma \in C^1\left(\Sc,\R^2\right)$ that
have 
exactly
$n$ transverse self-intersections of multiplicity two and no other self-intersections.
\end{definition}

Note that  $\HF(\eta,j_\pm,s,\omega) \subset \HF_{j_+-j_-}$ by virtue
of \eqref{connection_AI_number_self_intersections}.
Restricting to arclength para\-metrized curves fixes the length which takes care
of the missing scale-invariance of the energy. Indeed, one has
\begin{equation}\label{eq:scaling}
\TP_{q,\delta}(R\g)=R^{2-q}\TP_{q,\delta}(\g)\quad\Foa R>0.
\end{equation}

Observe that   the truncated tangent-point energy can be rewritten as
\begin{align}
\textstyle
\TP_{q,\delta}(\gamma)=\iint\limits_{\Ydch} \left(\frac{2\; |P^\perp_{\g'(x)}
(\g(x+w)-\g(x))|}{|\g(x+w)-\g(x)|^2}\right)^q |\gamma'(x+w)||\gamma'(x)| dw dx,
\label{eq:new-TP-formula}
\end{align}
where  for any unit vector $\nu\in\S^1\subset\R^2$ the expression
$P^\perp_\nu$ denotes the orthogonal
projection onto its orthogonal complement $\nu^\perp$. The expression \eqref{eq:new-TP-formula} is quite useful to relate the energy to a fractional
Sobolev seminorm and will be used quite frequently from now on.

\begin{lemma}
\label{Sobolev_endliche_Energie}
Let $q> 2$, $n \in \N$, $\gamma \in 
W^{2-\frac{1}{q},q}\left(\Sc,\R^2\right)\cap \HF_n$ and 
$0<\delta<\frac{\Lambda(\gamma)}{2}$. 
Then $\TP_{q,\delta}(\gamma)<\infty$.
\begin{proof}
With the bilipschitz estimate \eqref{eq:global-individual-bilip}
from \cref{bilipschitz_1} we can bound the integrand in \eqref{eq:new-TP-formula}
from above
by
$$
\textstyle 2^qc^{2q}|w|^{-2q}\left|P_{\gamma'(x)}^\perp \left(\int_0^1 \frac{d}{dt}\gamma(x+tw)dt \right) 
\right|^q = 2^qc^{2q}|w|^{-q}
\left|\int_0^1 P_{\gamma'(x)}^\perp (\gamma'(x+tw)-\gamma'(x))  dt\right|^q,
$$
which is now integrable over all of $\R/\Z\times [-\frac12,\frac12].$ 
Indeed, by Jensen's inequality,  $\|P^\perp\|\le 1$, and Fubini, we estimate
\begin{align*}
\textstyle\frac{1}{2^q}\TP_{q,\delta}(\gamma)
&\leq \textstyle
c^{2q} \int_0^1 \int_{\R/\Z} \int_{-\frac{1}{2}}^{\frac{1}{2}} \frac{|\gamma'(x+tw)-\gamma'(x)|^q}{|w|^q}dw dx dt.
\end{align*}
Now we change variables according to $\sigma(w):=tw$ to arrive at
\begin{align*}
\textstyle 
 \int_0^1 \int_{\R/\Z} \int_{-\frac{t}{2}}^{\frac{t}{2}}\frac{|\gamma'(x+\sigma)-\gamma'(x)|^q}{|\sigma|^q}t^{q-1} d\sigma dx dt
\leq \frac{1}{q} \int_{\R/\Z} \int_{-\frac{1}{2}}^{\frac{1}{2}}\frac{|\gamma'(x+\sigma)-\gamma'(x)|^q}{|\sigma|^q}d\sigma dx
\textstyle =\frac{1}{q} {\left[\gamma'\right]^q}_{1-\frac{1}{q},q}. 
\end{align*}
as an upper bound for $(2c^2)^{-q}\TP_{q,\delta}$. Notice that the constant $c$ depends 
on $\g$ and $\delta$ so that this  bound does not yield an a priori
estimate on the truncated energy.
\end{proof}
\end{lemma}

If we want to investigate, on the other hand,
whether the truncated tangent point energy $\TP_{q,\delta}$
also regularizes the curve to belong to the 
fractional Sobolev space $W^{2-\frac{1}{q},q}(\Sc,\R^2)$, we need to 
first relate
the minimal parameter distance of intersection pairs to the truncation
parameter $\delta$. Namely, for $n\in\N$ and $\g\in\HF_n$ set
\begin{equation}\label{eq:lambda}
\lambda:=\min_{(u,v) \in S(\gamma)}\left|u-v \right|_{\Sc}, \quad\textnormal{with $S(\gamma)$ as in \eqref{def_S_gamma}.}\notag
\end{equation}
It turns out that 
the truncated parts of the domain stay away from the diagonal, if $\delta$ is chosen
sufficiently small.
\begin{lemma}\label{lem:extra-parameter-lemma}
For $\g\in\HF_n$, $\delta\in (0,\frac{\Lambda(\g)}4]$, and any $\tau\in [0,\frac12]$ such that $|\g'(s+w)-\g'(s)|\le\frac12$
for all $s\in\R/\Z$ and $|w|\le\tau$, one has the inequality
$\lambda-2\delta >\tau$. In particular, $\overline{\hat{Y}_\delta(\g)}\subset \R/\Z\times \{\tau<|w|\le\frac12\}.$
\begin{proof}
It suffices to show that $\tau <\frac{\lambda}2$ since then we find by means
of $\delta\le\frac{\Lambda(\g)}4\le\frac{\lambda}4$ the desired inequality
$\lambda-2\delta\ge\frac{\lambda}2>\tau.$
Take $(u,v) \in S(\g)$ with $0\le u<v<1$ and $|u-v |_{\Sc}=\lambda$, 
and assume to the contrary that
$\tau\ge \frac{\lambda}2$. This implies $|t-\frac12(u+v)|_{\Sc}\le\tau$ for all
$t\in [u,v]$. Hence we infer from \eqref{eq:omega}
for $s:=u,$ $w:=v-u$ and $t^*:=\frac12 (u+v)$
the contradiction 
$$
|\g(v)-\g(u)|\ge\frac12 |v-u|_{\Sc}>0,
$$
which proves the inequality.
For the last statement take a point $(x,w)$ in the closure of 
$\hat{Y}_\delta(\g)$ with $|x-u_i|_{\Sc}\le\delta$ and $|x+w-v_i|_{\Sc}\le\delta$ for some $(u_i,v_i)\in S(\g)$. Then estimate
$\frac12\ge|w|\ge |u_i-v_i|_{\Sc}-|v_i-(x+w)|_{\Sc}-|x-u_i|_{\Sc}\ge\lambda-2\delta >\tau.$
\end{proof}
\end{lemma}

The following a priori estimates for the seminorm and the oscillation of
 the tangent in terms of the truncated  tangent point energy $\TP_{q,\delta}$
 will be important in the proof of 
$C^1$-compactness of sublevel sets in $\HF(\eta,j_\pm,s,\omega)\subset\HF_n$; 
see \cref{compactness}. 

\begin{lemma}[Energy controls seminorm and oscillation of tangents]
\label{regularity}
Let $n \in \N$, $q>2$,  $\gamma \in \HF_n$, and  $\delta \in 
(0,\frac{\Lambda(\gamma)}{4}]$. If $\TP_{q,\delta}(\gamma)<\infty$, 
then there exist constants $c_1(q),c_2(q)>0$ only depending on $q$ such that
\begin{align}\label{eq:seminorm-estimate-TP}
[\gamma']_{1-\frac{1}{q},q}^q\leq c_1(q)\Big[\TP_{q,\delta}(\gamma)+\left(\TP_{q,\delta}(\gamma)\right)^{\frac{q-1}{q-2}} \Big].
\end{align}
Furthermore, 
\begin{align}\label{eq:tangent-oscillation-TP}
|\gamma'(x)-\gamma'(y)|\leq c_2(q)\TP^\frac1{q}_{q,\delta}(\gamma) |x-y|_{\Sc}^{1-\frac{2}{q}} \quad \text{for all } x,y \in \Sc. 
\end{align}
\begin{proof}
Choose $\tau\in [0,\frac12]$ such that, on the one hand,
\begin{equation}\label{eq:tan-osc-grob}
|\g'(s+w)-\g'(s)|\le\frac12\quad\Foa s\in\R/\Z,\,|w|\le\tau,
\end{equation}
and such that there exist\footnote{Such parameters do exist since $\g$ is a closed curve
so that the unit tangent vectors cannot be contained in a cone with opening angle
$\alpha\le\frac{\pi}2\sin\frac{\alpha}2=\frac{\pi}4|\g'(s+w)-\g'(s)|\le\frac{\pi}8.$
}, on the other hand, parameters $s_0\in\R/\Z$ and $|w_0|\le\tau$
with
\begin{equation}\label{eq:tan-equality}
|\g'(s_0+w_0)-\g'(s_0)|=\frac12.
\end{equation}
Now we proceed in two steps. We first establish a preliminary upper bound for
the seminorm of $\g'$ in terms of the truncated energy and the parameter $\tau$
that still depends on $\g$. To bound $\tau$ uniformly from below again in terms
of $\TP_{q,\delta}$ we insert \eqref{eq:tan-equality} into
a Morrey-type inequality which we finally prove in the second step. This inequality
also yields
\eqref{eq:tangent-oscillation-TP}.

{\it Step 1.}\,
Lemma \ref{lem:extra-parameter-lemma} implies for the seminorm on the truncated region
$\hat{Y}_\delta(\g)$
\begin{align}\label{eq:seminorm-truncated-region}
\textstyle\iint\limits_{\hat{Y}_\delta(\gamma)} 
\frac{|\gamma'(x+w)-\gamma'(x)|^q}{|w|^{q}}dw dx
\leq 2^q\iint\limits_{\substack{\hat{Y}_\delta(\gamma)\cap \{|w|> \tau \}}}
\frac{1}{|w|^q}dw dx\leq \frac{2^{q+1}}{q-1}\tau^{1-q}.
\end{align}
Similarly, we obtain on the domain of $\TP_{q,\delta}$ away from the diagonal
\begin{equation}\label{eq:seminorm-off-diagonal}
\textstyle\iint\limits_{\substack{\Ydch \\ \cap\{|w|> \tau \}}} \hspace{-2em} \frac{|\g'(x+w)-\g'(x)|^q}{|w|^{q}}dx dw
\leq \frac{2^{q+1}}{q-1}\tau^{1-q},
\end{equation}
so that it remains to investigate the domain of $\TP_{q,\delta}$ near the diagonal. 
For that we change variables to split the truncated energy \eqref{eq:new-TP-formula}
into two equal summands
as
\begin{equation}\label{eq:split-energy}
\textstyle\TP_{q,\delta}(\g)=\frac12 \iint\limits_{\Ydch} \hspace{-2em}
\frac{|2P_{\gamma'(x+w)}^{\perp}(\gamma(x+w)-\gamma(x))|^q +
|2P_{\gamma'(x)}^{\perp}(\gamma(x+w)-\gamma(x))|^q}{|\gamma(x+w)-\gamma(x)|^{2q}} dwdx,
\end{equation}
and bound the resulting numerator (including the prefactor $\frac12$) from below by  the expression
\begin{align}
\label{proj_q}
N^\frac{q}2:=\textstyle\big|P_{\gamma'(x+w)}^{\perp}\big(\gamma(x+w)-\gamma(x)\big)
-P_{\gamma'(x)}^{\perp}\big(\gamma(x+w)-\gamma(x)\big) \big|^q.
\end{align}
Abbreviating the difference $\triangle f\equiv\triangle_w f(x):=f(x+w)-f(x)$ for any
function $f:\R/\Z\times [-\frac12,\frac12]$ and using the explicit form
$P^\perp_\nu(\xi)=\xi-\langle\xi,\nu\rangle\nu$ of the orthogonal projection
onto the orthogonal complement of $\nu\in\S^1$ we can estimate $N$ from below as
\begin{align*}
N &= \langle\triangle\g,\g'(x+w)-\g'(x)\rangle^2+2\langle\triangle\g,\g'(x+w)\rangle
\langle\triangle\g,\g'(x)\rangle \underbrace{(1-\langle\g'(x+w),\g'(x)\rangle)}_{=\frac12 |\triangle\g'|^2}\\
&\ge \langle\triangle\g,\g'(x+w)\rangle
\langle\triangle\g,\g'(x)\rangle |\triangle\g'|^2.
\end{align*}
With $\triangle\g=w\int_0^1\g'(x+\sigma w)\,d\sigma$ one obtains
$$
\textstyle N\ge w^2\int_0^1\langle\g'(x+\sigma w),\g'(x+w)\rangle\,d\sigma\int_0^1
\langle\g'(x+\theta w),\g'(x+w)\rangle\,d\theta\cdot|\triangle\g'|^2,
$$
and again by means of $\langle a,b\rangle=1-\frac12 |a-b|^2$ for $a,b\in\S^1$ applied
to both integrands we finally obtain for 
 $x \in \Sc$ and $w \in [-\tau,\tau]$ because of \eqref{eq:tan-osc-grob}
\begin{align}
\label{proj_qq}
N^\frac{q}2\ge 
\textstyle\big(\frac{7}{8} \big)^q|w|^q|\gamma'(x+w)-\gamma'(x)|^q.
\end{align}
Combining \eqref{eq:split-energy},\eqref{proj_q}, and  \eqref{proj_qq} with
the Lipschitz estimate $|\g(x+w)-\g(x)|\le|w|$ in the denominator we arrive at
\begin{equation}\label{eq:near-diagonal}
\textstyle\iint\limits_{\substack{\Ydch \\ \cap\{|w|\leq \tau \}}}\hspace{-2em} 
\frac{|\gamma'(x+w)-\gamma'(x)|^q}{|w|^{q}}dwdx\le \left(\frac{8}{7}\right)^q
\TP_{q,\delta}(\gamma).
\end{equation}
Adding up \eqref{eq:seminorm-truncated-region}, 
\eqref{eq:seminorm-off-diagonal}, and \eqref{eq:near-diagonal} yields
\begin{equation}
\label{zwischen}
[\g']_{1,1-\frac1{q}}^q=\textstyle\iint\limits_{\R/\Z\times [-\frac12,\frac12]}
\frac{|\gamma'(x+w)-\gamma'(x)|^q}{|w|^{q}}dw dx
\leq \left(\frac{8}{7}\right)^q\TP_{q,\delta}(\gamma)+\frac{2^{q+2}}{q-1}\tau^{1-q}.
\end{equation}
Note that $\tau$ still depends on the curve $\gamma$, so this is not yet the 
desired a priori  estimate for the seminorm. Furthermore, the second term blows up if $\tau$ tends to zero.

We need to establish a uniform lower bound on $\tau$ 
in order to prove the a priori estimate \eqref{eq:seminorm-estimate-TP}. 
In step 2 we will show that there is a constant  $c_M(q)>0$ independent of $\g$
such that
\begin{align}
\label{Hoelder1}
\sup_{x \in \Sc}|\gamma'(x+w)-\gamma'(x)|\leq c_M(q)\TP^\frac1{q}_{q,\delta}(\gamma)|w|^{1-\frac{2}{q}} 
\quad \text{for all  $ |w|\leq \frac{\tau}{4}.$}
\end{align}
Then the explicit choice of $s_0$ and $w_0$ satisfying \eqref{eq:tan-equality}
 and the triangle inequality yield
\begin{align}
\label{delta0}
\textstyle\frac12=|\gamma'(s_0+w_0)-\gamma'(s_0)| 
\bleq[\eqref{Hoelder1}] 4c_M(q)\TP^\frac1{q}_{q,\delta}(\gamma)\big|\frac{w_0}{4} \big|^{1-\frac{2}{q}}
\leq 4^{\frac{2}{q}}c_M(q)\TP^\frac1{q}_{q,\delta}(\gamma)\tau^{1-\frac{2}{q}}.
\end{align}
Hence, setting $\tilde{c}(q):=(2c_M(q))^{\frac{q(q-1)}{q-2}}16^{\frac{q-1}{q-2}}$ leads to 
$
\tau^{1-q}\leq\tilde{c}(q)\left(\TP_{q,\delta}(\gamma)\right)^{\frac{q-1}{q-2}}.
$
Finally, setting $c_1(q):=\max\left\{\left(\frac{8}{7}\right)^q,\frac{2^{q+2}}{q-1}\tilde{c}(q) \right\}$ implies by virtue of \eqref{zwischen} inequality 
\eqref{eq:seminorm-estimate-TP}. 
Moreover, \eqref{Hoelder1} is also the essential local
estimate which leads to the global oscillation bound \eqref{eq:tangent-oscillation-TP}
for the tangent. Indeed, for $|x-y|_{\Sc}>\frac{\tau}4$ one simply estimates
$|\g'(x)-\g'(y)|\le 2<2(\frac4{\tau}|x-y|_{\Sc})^{1-\frac2{q}}$, which 
together 
with \eqref{delta0} 
leads to
$$
|\g'(x)-\g'(y)|<16 c_M(q)\TP_{q,\delta}^\frac1{q}(\g)|x-y|^{1-\frac2{q}}_{\Sc}\quad\Foa |x-y|_{\Sc}>\frac{\tau}4.
$$
Combining this with \eqref{Hoelder1} yields \eqref{eq:tangent-oscillation-TP} 
if we set
$c_2(q):=16 c_M(q).$

{\it Step 2.}\,
To complete the proof, we need to show inequality \eqref{Hoelder1}. 
For $r\in (0,\frac{\tau}2]$ define the open neighborhood
$M_r:=\bigcup_{x \in \Sc} \left(B_r(x)\times B_r(x) \right)$ of the diagonal in $\R/\Z\times
\R/\Z$, where $B_r(x):=\{y\in\R/\Z: |y-x|_{\Sc}<r\}$.
Lemma \ref{lem:extra-parameter-lemma} implies that $M_r$ is contained in $\Ydc$.
Set $\gamma'_{B_r(x)}:=\mvint_{B_r(x)}\g'(z)\,dz\equiv
\frac{1}{2r}\int_{B_r(x)}\g'(z)\,dz$ and 
apply H\"older's inequality to estimate for $x \in \Sc$
\begin{align}
\label{l_p_diff}
&\textstyle \frac{1}{2r}\int_{B_r(x)}|\gamma'(y)-\gamma'_{B_r(x)} |dy 
\textstyle\bleq \frac{1}{4r^2}
\iint\limits_{B_r(x)\times B_r(x)}|\gamma'(y)-\gamma'(z)|dz dy\nonumber\\
&\textstyle\leq \Big(\frac{1}{4r^2}\iint\limits_{B_r(x)\times B_r(x)}|\gamma'(y)-\gamma'(z)|^q dz dy\Big)^{\frac{1}{q}}
\leq \textstyle \Big((2r)^{q-2}\iint\limits_{B_r(x)\times B_r(x)}\frac{|\gamma'(y)-\gamma'(z)|^q}{|y-z|^q}dz dy\Big)^{\frac{1}{q}}\nonumber\\
&\leq \textstyle
2^{1-\frac{2}{q}}r^{1-\frac{2}{q}}\Big( \; \; 
\iint\limits_{\substack{\Ydc \\ \cap\{|y-z|< 2r \}}}
\frac{|\gamma'(y)-\gamma'(z)|^q}{|y-z|^q}dz dy \Big)^{\frac{1}{q}}
\textstyle\bleq[\eqref{eq:near-diagonal}]2^{1-\frac{2}{q}}r^{1-\frac{2}{q}}\frac{8}{7}\TP^\frac1{q}_{q,\delta}(\gamma).
\end{align}
For two  points $x,y \in \Sc$ 
with $r:=|x-y|\in (0,\frac{\tau}{4}]$ one finds
\begin{align}
\label{reihe}
&|\gamma'(x)-\gamma'(y)| \le \limsup_{k\to\infty}|\g'_{B_{2^{-k}r}(x)}-
\g'_{B_{2^{-k}r}(y)}|\\
&\textstyle\le \sum_{l=0}^\infty\big|\gamma'_{B_{2^{-l}r}(x)}-\gamma'_{B_{2^{-l+1}r}(x)}
\big|
+\big|\gamma'_{B_{2r}(x)}-\gamma'_{B_{2r}(y)}\big|
+\sum_{l=0}^\infty\big|\gamma'_{B_{2^{-l}r}(y)}-\gamma'_{B_{2^{-l+1}r}(y)}\big|.\notag
\end{align}
Averaging the inequality $|\g'_{B_{2r}(x)}-\g'_{B_{2r}(y)}|\le
|\g'_{B_{2r}(x)}-\g'(z)|+|\g'(z)-\g'_{B_{2r}(y)}|$ over all $z\in B_{2r}(x)\cap
B_{2r}(y)$ in combination with the inclusions $B_r(x),B_r(y)\subset
B_{2r}(x)\cap B_{2r}(y)\subset B_{2r}(x),B_{2r}(y)$ bounds the middle summand in
\eqref{reihe} as
\begin{align}
\label{middle}
&\textstyle\big|\gamma'_{B_{2r}(x)}\!-\!\gamma'_{B_{2r}(y)}\big| \le \textstyle
\frac{1}{|B_{2r}(x)\cap B_{2r}(y)|}\Big[
\int\limits_{B_{2r}(x)}|\gamma'_{B_{2r}(x)}\!-\!\gamma'(z)|dz
+\int\limits_{B_{2r}(y)}|\gamma'_{B_{2r}(y)}\!-\!\gamma'(z)|dz\Big] \nonumber\\
\leq&\textstyle 
2\Big[\mvint_{B_{2r}(x)}|\gamma'_{B_{2r}(x)}\!-\!\gamma'(z)|dz
+\mvint_{B_{2r}(y)}|\gamma'_{B_{2r}(y)}\!-\!\gamma'(z)|dz \Big]\notag\\
&
\overset{\eqref{l_p_diff}}{\le}\textstyle\frac{16}{7}(4r)^{1-\frac2{q}} 
\TP^\frac1{q}_{q,\delta}(\gamma)
=:\hat{c}_1(q)r^{1-\frac2{q}}\TP^\frac1{q}_{q,\delta}(\gamma). 
\end{align}
For the remaining terms in the series, we use the equality 
$\left|B_{2R}(x)\right|=2\left|B_R(x)\right|$ for any $R \in ( 0,\frac{\tau}{4}]$ to estimate again by means of \eqref{l_p_diff}
\begin{align}
|\gamma'_{B_{2R}(x)}& -\gamma'_{B_R(x)}|\textstyle
\leq \mvint_{B_R(x)}|\gamma'_{B_{2R}(x)}-\gamma'(z)|dz+\mvint_{B_R(x)}
|\gamma'_{B_{R}(x)}-\gamma'(z)|dz\notag\\
&\leq\textstyle 2\mvint_{B_{2R}(x)}|\gamma'_{B_{2R}(x)}-\gamma'(z)|dz
+\mvint_{B_R(x)}|\gamma'_{B_{R}(x)}-\gamma'(z)|dz\notag\\
&\textstyle\overset{\eqref{l_p_diff}}{\le} \frac{8}{7}R^{1-\frac{2}{q}}
\TP^\frac1{q}_{q,\delta}(\gamma)
\big(2 \cdot 4^{1-\frac{2}{q}}+ 2^{1-\frac{2}{q}}\big)
=:\hat{c}_2(q)R^{1-\frac{2}{q}}\TP^\frac1{q}_{q,\delta}(\gamma).
\label{eq:summands-series-estimate}
\end{align}
Inserting \eqref{middle} and\eqref{eq:summands-series-estimate} 
for $R:=2^{-l}r$, $l=0,1,2,\ldots$ into \eqref{reihe}
yields
\begin{align*}
|\gamma'(x)-\gamma'(y)|&\textstyle\leq \max\{\hat{c}_1(q),\hat{c}_2(q)\}
\TP^\frac1{q}_{q,\delta}(\gamma)r^{1-\frac{2}{q}}
\Big(\sum_{l=0}^\infty 2^{-l(1-\frac{2}{q})}+ \sum_{k=0}^\infty 2^{-k(1-\frac{2}{q})}
\!+\!1\Big)\\
&=:\textstyle c_M(q)\TP^\frac1{q}_{q,\delta}(\gamma)|x-y|^{1-\frac{2}{q}}\quad
\Foa
|x-y|_{\Sc}\le\frac{\tau}4.\qedhere
\end{align*}
\end{proof}
\end{lemma}

Our next goal is to prove a \emph{global bilipschitz estimate} 
away from self-intersection points for curves of finite
truncated tangent-point energy where the bilipschitz constant depends \emph{only} on the
energy bound and the intrinsic distance between any two self-intersections; see Proposition \ref{bilipschitz} below. 
 It is crucial for the compactness
of sublevel sets of $\TP_{q,\delta}$ in $\HF_n(\eta,j_\pm,s,\omega)$ proven in Theorem \ref{compactness} in Section
\ref{chapter_existence_minimizers} that the bilipschitz constant in
Proposition \ref{bilipschitz} does \emph{not}
depend on the  shape of the individual curve as it does in Lemma \ref{bilipschitz_1}.
To prepare  this we bound 
  $\TP_{q,\delta}$ uniformly from below for normalized curves that possess a 
  secant perpendicular to the curve in one of its two endpoints.
To make this precise, define for $d\ge 2$,
$\alpha\in (0,1]$, $L>2$, and $H\in (0,\infty)$ the set\footnote{Note that without the  condition $L>2$ the set
$\mathscr{L}(\alpha,H,L,d)$ might be empty because of condition \ref{eq:L1}.} 
$\mathscr{L}(\alpha,H,L,d)$ to be 
the collection of all arclength
parametrized loops $\G\in C^1(\R/L\Z,\R^d)$
such that
\begin{enumerate}
\Item[\rm ($\mathscr{L} 1$)]\label{eq:L1}\quad
there are $s,t\in\R/L\Z$ with $|\G(s)-\G(t)|=1$;
\Item[\rm ($\mathscr{L} 2$)]
\quad \label{eq:L2}
$\G'(s)\perp \big(\G(s)-\G(t)\big)$;
\Item[\rm ($\mathscr{L} 3$)]\quad \label{eq:L3}
$|\G'(\sigma)-\G'(\tau)|\le H |\sigma-\tau|_{\Sc}^\alpha$ for all $\sigma,\tau
\in\R/L\Z$.
\end{enumerate}
\begin{lemma}[Least energy contribution]\label{lem:least-energy-contribution}
Let $q\in (0,\infty)$  and $\G\in\mathscr{L}
(\alpha,H,L,d)$.
Then
\begin{equation}\label{eq:least-energy}
\int_{I_\epsilon(s)}\int_{I_{\rho}(t)}\frac{d\tau d\sigma}{
(\rtp[\G](\G(\sigma),\G(\tau)))^q}>\frac1{6^q}\epsilon\rho
\end{equation}
for all $\rho\le\frac1{8},$ $\epsilon\le\epsilon_0(\alpha,H):=
\min\{\frac18,(\frac1{8H})^\frac1{\alpha}\},$ where
$I_\epsilon(s)\in\{(s-\epsilon,s),(s,s+\epsilon)\}$ and $
I_{\rho}(t)\in\{(t-\rho,t),(t,t+\rho)\}$.
\begin{proof}
Fix an arbitrary loop $\G\in\mathscr{L}(\alpha,H,L,d)$ with two parameters
$s,t\in\R/L\Z$ satisfying conditions \ref{eq:L1} and \ref{eq:L2}.
In particular, \ref{eq:L2} implies
\begin{equation}\label{eq:000}
\textstyle\rtp[\G](\G(s),\G(t))=
\frac12 |\G(s)-\G(t)|=\frac12.
\end{equation}
Keep $t$ fixed for a moment and abbreviate the numerator and the denominator of 
the tangent-point radius for any  
$\sigma\in\R/L\Z$  as
$$
\textstyle
\rtp[\G](\G(\sigma),\G(t))=\frac{|\G(t)-\G(\sigma)|^2}{2\dist\big(\G(t),\G(\sigma)+\R\G'(\sigma)\big)}\coloneqq \frac{N(\sigma)}{D(\sigma)},
$$
and use the triangle inequality and \eqref{eq:000} to find
$$
\textstyle\frac{N(\sigma)}{D(\sigma)}\le\frac12 + \left|\frac{N(\sigma)}{D(\sigma)}-
\frac{N(s)}{D(s)}\right|\le\frac12+\frac1{2D(\sigma)} \big|2N(\sigma)-D(\sigma)\big|,
$$
where we also used condition \ref{eq:L1}   which
implies $2N(s)=D(s)=2$. 
If $N(\sigma)>N(s)=1$ we estimate by means of the binomial identity and \ref{eq:L1}  
\begin{align}
|N(\sigma)-N(s)| &= \big(|\G(t)-\G(\sigma)|+1\big)\cdot\big(|\G(t)-\G(\sigma)|-1\big)
\notag\\
& \le \big(2+|\G(\sigma)-\G(s)|\big)|\G(s)-\G(\sigma)|<3|s-\sigma|
\label{eq:002}
\end{align}
for all $|s-\sigma|<\epsilon_1\coloneqq 1<\frac{L}2$ by means of  $L>2
$.
In any case, $N(\sigma)<2$ for $|s-\sigma|<\epsilon_2:=\frac13<\epsilon_1$.
With conditions \ref{eq:L2} and \ref{eq:L3}   
one  finds
\begin{align*}
\frac12 |D(\sigma)-D(s)| & \le |\G(t)-\G(\sigma)-\langle\G(t)-\G(\sigma),
\G'(\sigma)\rangle \G'(\sigma)-(\G(t)-\G(s))|
\notag\\
&\hspace{-2cm} = \big|\G(s)-\G(\sigma)-\big[\langle\G(t)-\G(\sigma),
\G'(\sigma)-\G'(s)\rangle+\langle\G(s)-\G(\sigma),
\G'(s)\rangle\big]\G'(\sigma)\big|\notag\\
& \hspace{-2cm} \le 2|s-\sigma|+2H|s-\sigma|^\alpha, 
\end{align*}
where we also used again
that $\G$ has Lipschitz constant $1$, and therefore also
$|\G(t)-\G(\sigma)|\le |\G(t)-\G(s)|+|\G(s)-\G(\sigma)|\le 1+|s-\sigma|<2$ for 
all $|s-\sigma|<\epsilon_1$. 
This implies in particular that $D(\sigma)\in (1,3)$ for all $|s-\sigma|<
\epsilon_0:=\min\{\frac18,\left(\frac{1}{8H}\right)^{\frac1{\alpha}}\}
<\epsilon_2$  so that we infer
$
\rtp[\G](\G(\sigma),\G(t))=\frac{N(\sigma)}{D(\sigma)} < 2$ for all $ |s-\sigma|<\epsilon_0.
$
Now abbreviate $
\frac{Z(\sigma,\tau)}{W(\sigma,\tau)}\coloneqq \rtp[\G](\G(\sigma),\G(\tau))$.
Similarly as 
in \eqref{eq:002} we estimate  in case $Z(\sigma,\tau)>N(\sigma)$
\begin{align*}
|Z(\sigma,\tau)-N(\sigma)| & = \big( |\G(\tau)-\G(\sigma)|+\sqrt{N(\sigma)}\big)\cdot
\big(|\G(\tau)-\G(\sigma)|-\sqrt{N(\sigma)}\big)\notag\\
& \le \big(2\sqrt{2}+|t-\tau|\big)|t-\tau|<3|t-\tau|
\end{align*}
for all $|s-\sigma|<\epsilon_0$ and $|t-\tau|<\frac1{8}\eqqcolon \rho_1$. In particular,
$Z(\sigma,\tau)<2+3|t-\tau|<3$ for all $|t-\tau|<\rho_1$. 
Moreover,
\begin{align*}
\textstyle\frac12 |W(\sigma,\tau)-D(\sigma)| & \le |\G(\tau)-\G(\sigma)
-\langle\G(\tau)-\G(\sigma),
\G'(\sigma)\rangle\G'(\sigma)-D(\sigma)|\notag\\
&\hspace{-1cm} =\textstyle
|\G(\tau)-\G(t)-\langle\G(\tau)-\G(t),\G'(\sigma)\rangle\G'(\sigma)|
\le 2|\tau-t|<\frac1{4},
\end{align*}
in particular, $W(\sigma,\tau)>\frac12$ for $|s-\sigma|<\epsilon_0$ and $|t-\tau|<\rho_1$. 
Thus, we easily obtain
the rough estimate $\rtp[\G](\G(\sigma),\G(\tau))<6$ for all $|s-\sigma|<\epsilon_0$
and $|t-\tau|<\rho_1$, from which we infer the claim by integrating over 
$I_\epsilon(s)\times I_{\rho}(t)$ for $\epsilon\in (0,\epsilon_0]$ and 
$\rho\le\rho_1=\frac1{8}$. 
\end{proof}
\end{lemma}

\begin{proposition}[Global uniform bilipschitz estimate]
\label{bilipschitz}
Let $q>2$,  $\eta\in (0,1)$, and $E>0$.  
Then there exists a constant $c=c(q,E,\eta)>0$ such that for all $\gamma \in 
\HF(\eta,j_\pm,s,\omega)$ with $\TP_{q,\delta}(\gamma)\leq E$ 
\begin{align*}
|x-y|_{\Sc}\leq c(q,E,\eta)|\gamma(x)-\gamma(y)|
\end{align*}
for all $(x,y) \in \Ydc$ and all $\delta\in\left(0,\frac{\eta}{2} \right]$.
\begin{proof}
Combining the basic inequality \eqref{eq:omega} with the uniform
oscillation estimate \eqref{eq:tangent-oscillation-TP} for the tangents $\g'$
of any curve $\g\in\HF(\eta,j_\pm,s,\omega)\subset\HF_n$ for $n=j_+-j_-$
we find a constant
$\tau_1=\tau_1(q,E)>0$ such that
\begin{equation}
\label{eq:unif-local-bilip}
\textstyle|\g(y)-\g(x)|\ge\frac12 |y-x|_{\Sc}\quad\Foa |y-x|_{\Sc}\le\tau_1.
\end{equation}
As in the proof\footnote{With the difference that $\tau_1$ does
\emph{not} depend on the individual shape of the curve in contrast to
$\tau_\g$ in \eqref{eq:local-individual-bilip} in the proof of Lemma \ref{bilipschitz_1}.} of Lemma 
\ref{bilipschitz_1}
we find that
\begin{align*}
C:=\inf\big\{|\gamma(x)-\gamma(y)| \; : \; (x,y) \in \Ydc \text{ and } 
|x-y|_{\Sc}\geq \tau_1 \big\}
\end{align*}
is strictly positive and attained by at least one parameter pair $(x_0,y_0)$
because of uniform continuity of the mapping $(x,y)\mapsto
|\gamma(x)-\gamma(y)|^2$ on compact subsets of $(\R/\Z)^2$.
We can assume $C\in (0,\frac{1}{2})$ since otherwise we 
could conclude the
proof with
\begin{align*}
  |\gamma(x)-\gamma(y)|\ge \frac{1}{2}\geq |x-y|_{\Sc} \quad \text{for all } (x,y) \in \Ydc \text{ with } |x-y|_{\Sc}\geq \tau_1.
\end{align*}
If $|x_0-y_0|_{\Sc}=\tau_1$, we estimate with \eqref{eq:unif-local-bilip}
$2C=2|\gamma(x_0)-\gamma(y_0)|\geq \tau_1$ and thus
\begin{align*}
\textstyle |x-y|_{\Sc}\leq \frac{1}{2}\leq \frac{C}{\tau_1}\leq 
\frac{|\gamma(x)-\gamma(y)|}{\tau_1} 
\quad \text{for all } (x,y) \in \Ydc,  |x-y|_{\Sc}\geq \tau_1,
\end{align*}
which finishes the proof in that case.
Now, assume $|x_0-y_0|_{\Sc}>\tau_1$. The goal is to derive a uniform lower 
bound on $C$ in order to finish the proof.
Define the intervals $I^+_{r}(x):=(x,x+r)$ and $I^-_r(x):=(x-r,x)$ for $x\in\R/\Z$ and
$r>0$. Then it  is easy to check by means of $0<\delta\le\frac{\eta}2 
<\frac{\Lambda}4$ (cf. Definition \ref{admissibile_curves})
%
that 
at least one of the following two inclusions is true
\begin{align}
\textstyle I^+_{\eta}(x_0)\times
I^+_{\eta}(y_0)
\subset \Ydc ,\label{eq:inclusion1}\\
\textstyle I^-_{\eta}(x_0)\times
I^-_{\eta}(y_0)\label{eq:inclusion2}
\subset \Ydc.
\end{align}
If $(x_0,y_0)\in (\R/\Z)^2\setminus\overline{Y_\delta(\g)}$ then it realizes the
infimum $C$ as an interior point so that the gradient of the mapping
$(x,y)\mapsto |\g(x)-\g(y)|^2$ vanishes at $(x_0,y_0)$ which leads to 
 $\g'(x_0)\perp (\g(x_0)-\g(y_0))\perp\g'(y_0)$.
But also if $(x_0,y_0)$ lies in the boundary of $Y_\delta(\g)$, we can infer 
the orthogonality of the chord $\g(x_0)-\g(y_0)$ to  one of the tangent
vectors $\g'(x_0)$ or $\g'(y_0)$ since with $\delta<\eta$ both points are contained
in straight segments of the curve $\g$ near a self-intersection pair $(u,v)\in S(\g)$ 
(recall
\eqref{def_S_gamma} and Definition \ref{admissibile_curves}). Indeed, at least one of
the points, say $\g(x_0)$ lies on the intersection of one straight segment of $\g$ 
with the
circle of radius $\delta$ centered at the
double point $\g(u)=\g(v)$. On the other hand, $\g(y_0)$ sits on the other straight
segment of $\g$ and may be closer to the double point. But since $|\g(x_0)-\g(y_0)|=C$
is minimal the chord $\g(x_0)-\g(y_0)$ must be perpendicular to $\g$ at $\g(y_0)$.
Summarizing we can say that in any case one has $\g'(y_0)\perp (\g(x_0)-\g(y_0)).$
To bring Lemma \ref{lem:least-energy-contribution} into play we first look at the reparametrization $\eta(x):=
\g(Cx)$ for $x\in\R/\frac1{C}\Z$ and change variables according to $z(x,y)=(Cx,Cy)$ 
to compute 
\begin{equation}\label{eq:trafo}
\textstyle\TP_{q,\frac{\delta}C}(\eta)=\iint\limits_{(\R/\frac1{C}\Z)^2\setminus Y_\frac{\delta}C(\eta)}
\big(\frac{2|P^\perp_{\eta'(x)}(\eta(y)-\eta(x))|}{|\eta(y)-\eta(x)|^2}\big)^q
|\eta'(x)||\eta'(y)|dxdy=\TP_{q,\delta}(\g).\notag
\end{equation}
Together with the scaling behavior \eqref{eq:scaling} of the tangent-point 
energy we obtain for the rescaled curve $\Upsilon:=\frac1{C}\eta$
\begin{equation}\label{eq:energy-identity}
\textstyle\TP_{q,\frac{\delta}C}(\Upsilon)=\big(\frac1{C}\big)^{2-q}
\TP_{q,\frac{\delta}C}(\eta)=\big(\frac1{C}\big)^{2-q}\TP_{q,\delta}(\g)
\le \big(\frac1{C}\big)^{2-q}E.
\end{equation}
On the other hand, $\Upsilon$ is contained in the set $\mathscr{L}(\alpha,H,L,d)$
referred to in Lemma \ref{lem:least-energy-contribution}
for $d:=2$, $L:=\frac1{C}>2$. 
$\alpha:=1-\frac2{q},$ $H:=(\frac12)^{1-\frac2{q}}c_2(q)E^{\frac1{q}}$, where $c_2(q)$ is
the constant in the oscillation inequality \eqref{eq:tangent-oscillation-TP} of Lemma 
\ref{regularity}. Properties \ref{eq:L1} and \ref{eq:L2}
 hold for
$s:=\frac{y_0}C$ and $t:=\frac{x_0}C$.
From \eqref{eq:inclusion1} and \eqref{eq:inclusion2} we infer that at least one of
the following two inclusions holds true as well:
\begin{align}
\textstyle I^+_{\frac{\eta}{C}}(\frac{x_0}C)\times
I^+_{\frac{\eta}{C}}(\frac{y_0}C)
\subset (\R/\frac1{C}\Z)^2\setminus Y_{\frac{\delta}C}(\Upsilon),\label{eq:inclusion11}\\
\textstyle I^-_{\frac{\eta}{C}}(\frac{x_0}C)\times
I^-_{\frac{\eta}{C}}(\frac{y_0}C)
\subset  (\R/\frac1{C}\Z)^2\setminus Y_{\frac{\delta}C}(\Upsilon).\label{eq:inclusion22}
\end{align}
Assuming without loss of generality that inclusion \eqref{eq:inclusion11}
holds,  we can combine the lower bound \eqref{eq:least-energy} of Lemma \ref{lem:least-energy-contribution} for
$\G:=\Upsilon$ with the energy estimate \eqref{eq:energy-identity} to obtain
\begin{align}
\textstyle
C^{q-2}E &
\ge\TP_{q,\frac{\delta}C}(\Upsilon)\ge\textstyle
\int_{I^+_{\frac{\eta}{C}}(\frac{y_0}C)} \int_{I^+_{\frac{\eta}{C}}(
\frac{x_0}C)}
\frac{d\tau d\sigma}{(\rtp[\Upsilon](\Upsilon(\sigma),\Upsilon(\tau)))^q}\notag\\
&\textstyle >\frac1{6^q}\min\big\{\frac{\eta}{C},\epsilon_0\big\}
\min\big\{\frac{\eta}{C},\frac1{8}\big\} >\min\big\{2\eta,\epsilon_0\big\}
\min\big\{2\eta,\frac1{8}\big\},
  \label{eq:crucial-bound}
\end{align}
since $C\in (0,\frac12)$, 
where  $\epsilon_0$ is the constant from Lemma \ref{lem:least-energy-contribution}
depending on $\alpha$ and $H$, which themselves were identified above and depend 
only on $q$ and on $E$. Since $q>2$ the right-hand side of \eqref{eq:crucial-bound}
 gives us the desired lower bound for $C$ solely 
in terms of $\eta,$ $q,$ and $E$. 
\end{proof}
\end{proposition}

\section{\texorpdfstring{Minimizing the $\delta$-truncated energy}{Minimizing the delta-truncated energy}}
\label{chapter_existence_minimizers}
So far, we have proven  an a priori estimate for the fractional Sobolev seminorm for
curves in $\HF_n$ (Lemma \ref{regularity}), and a uniform bilipschitz estimate 
for curves in the subclass $\HF(\eta,j_\pm,s,\omega)\subset\HF_{j_+-j_-}$ (Proposition \ref{bilipschitz}), both under the assumption of uniformly bounded truncated energy $\TP_{q,\delta}$. A key ingredient in the direct method in
the calculus of variations to prove existence of minimizers 
is the lower semi-continuity of the energy. The truncated tangent-point energy is 
in fact lower semi-continuous if the lengths of the subarcs between self-intersections 
of the converging curves are bounded from below.

\begin{lemma}
\label{lower_sc}
Let $q,\Lambda_0>0$, 
$n \in \N$ and $\delta\in (0,\frac{\Lambda_0}2)$.
Furthermore, let $(\gamma_k )_{k \in \N} \subset \HF_n$ be a sequence with $\Lambda(\gamma_k)\geq \Lambda_0$ and $\gamma_k \to \gamma \in 
\HF_n$ in $C^1$ as $k \to \infty$. Then
\begin{align}
\TP_{q,\delta}(\gamma)\leq \liminf_{k \to \infty} \TP_{q,\delta}(\gamma_k).\label{eq:UHS}
\end{align}
\begin{proof}
We may assume that $\liminf_{k \to \infty} \TP_{q,\delta}(\gamma_k)<\infty$.
The parameter pairs of  self-intersection points are collected in the sets 
$S(\gamma_k)=\{ (u_1^k,v_1^k),\ldots,(u_n^k,v_n^k) \}$ and $S(\gamma)= 
 \{ (u_1,v_1),\ldots,(u_n,v_n) \}$ as defined in \eqref{def_S_gamma}.
  As the sequences $(u_i^k)_{k \in \N}$ and 
  $(v_i^k)_{k \in \N}$ are bounded, we can assume that 
  $u_i^k \to \tilde{u}_i$ and $v_i^k \to \tilde{v}_i$ 
  for some $\tilde{u}_i,\tilde{v}_i\in\R/\Z$, $i=1,\ldots,n$, as $k \to \infty$. 
  Since $\Lambda(\gamma_k)\geq \Lambda_0$ for all $k \in \N$ 
  it is ensured that the set 
  $\{\tilde{u}_i,\tilde{v}_i \; | \; i=1,\dots,n \}$ has exactly $2n$ 
  elements. The $C^1$-convergence then implies 
  $\gamma(\tilde{u}_i)=\gamma(\tilde{v}_i)$ for all $i=1,\dots,n$. 
  Since $\gamma$ has exactly $n$ self-intersections, this implies
  that for each $i\in\{1,\ldots,n\}$ there is exactly one 
  $j=j(i)\in
  \{1,\ldots,n\}$ such that  $(\tilde{u}_i,\tilde{v}_i )=(u_j,v_j)$.
  Assume for a moment that the integrands converge pointwise, i.e.,
\begin{align}
  \label{pw_convergence_lower_sc}
  \textstyle
\left(\frac{2
\dist(\ell_k(x),\gamma_k(y))}{|\gamma_k(x)-\gamma_k(y)|^2} \right)^q
\rchi_{(\R/\Z)^2\setminus Y_\delta(\gamma_k)}(x,y) 
\to \left(\frac{2\dist(\ell (x),\gamma(y))}{|\gamma(x)-\gamma(y)|^2} \right)^q \rchi_{(\R/\Z)^2\setminus Y_\delta(\gamma)}(x,y)
\end{align}
for almost all $(x,y ) \in (\R/\Z)^2$
as $k \to \infty$. Fatou's lemma then implies \eqref{eq:UHS}.

 Hence, it remains to prove the pointwise convergence in 
 \eqref{pw_convergence_lower_sc}. If $(x,y) \in (\R/\Z)^2\setminus 
 \overline{Y_\delta(\gamma)}$ with $x \neq y$, then due to the 
 convergence of the sequences $(u_i^k)_{k \in \N}$ and 
 $(v_i^k)_{k \in \N}$ there exists an index $k_1 \in \N$ such that 
 $(x,y) \in (\R/\Z)^2\setminus 
 \overline{Y_\delta(\gamma_k)}$ for every $k \geq k_1$. 
 The $C^1$-convergence $\g_k\to\g$ as $k\to\infty$
 then implies pointwise convergence in \eqref{pw_convergence_lower_sc} 
 for those tuples.  In the same manner one can argue that if 
 $(x,y) \in Y_\delta(\gamma)$, there exists a 
 $k_2 \in \N$ such that $(x,y) \in Y_\delta(\gamma_k)$ for all 
 $k \geq k_2$. 
 Hence, the pointwise convergence in \eqref{pw_convergence_lower_sc} is 
 proven 
 for all $(x,y)\in (\R/\Z)^2\setminus\partial
 Y_\delta(\g)$, and the exceptional set $\partial Y_\delta(\g)$ has measure zero.
\end{proof}
\end{lemma}

The following compactness theorem will be essential to show the existence of minimizers.  

\begin{theorem}
\label{compactness}
Let $q>2$, $j_+,j_-,s,l \in \Z$, and $\eta \in (0,1)$  
such that 
$\HF(\eta,j_\pm,s,\omega)\neq \emptyset$. 
If  $\delta \in (0,\frac{\eta}{2}]$ 
then for every $E<\infty$ the set
\begin{align*}
\mathscr{A}(\delta,E,\eta,j_\pm,s,\omega)
:=\left\{\gamma \in \HF(\eta,j_\pm,s,\omega) :  
\TP_{q,\delta}(\gamma)\leq 
E \textnormal{\,\,and\,\,} 0\in\gamma(\R/\Z) \,\right\}
\end{align*}
is sequentially compact with respect to the $C^1$-norm.
\begin{proof}
If $\mathscr{A}(\delta,E,\eta,j_\pm,s,\omega)=\emptyset$ there is
nothing to prove, so let
$(\gamma_k )_{k \in \N}$ be a sequence contained in 
$\mathscr{A}(\delta,E,\eta,j_\pm,s,\omega)$. 
Let $n:=j_+-j_-$ be the number of the exclusively
transverse self-intersections of each curve $\g_k$ and
denote by $(u_i^k,v_i^k)$ the distinct 
parameters with $\gamma_k(u_i^k)=
\gamma_k(v_i^k)$
for every $k \in \N$ and $i=1,\dots,n$. Since the origin is contained 
in the 
image of every curve $\gamma_k$ and the curves 
are $1$-lipschitz due to arclength
parametrization, we have 
\begin{align*}
|\gamma_k(x)|\leq \textstyle\frac12\AND |\gamma_k(x)-\gamma_k(y)|
&\leq |x-y|_{\Sc} \quad \text{for all } x,y \in \Sc \text{ and } k \in \N. 
\end{align*}
By \cref{regularity}, we obtain a uniform H\"older estimate 
\begin{align*}
|\gamma_k'(x)-\gamma_k'(y)|&\leq c(q,E)|x-y|_{\Sc}^{1-\frac{2}{q}} \quad \text{for all } x,y \in \Sc \text{ and } k \in \N. 
\end{align*}
By means of the Arzel\`a-Ascoli theorem we obtain a subsequence, 
again denoted by $(\gamma_k )_{k \in \N}$, and a curve $\gamma \in C^1(\Sc,\R^2)$ 
such that $\gamma_k \to \gamma$ in $C^1$ as $k \to \infty$. It remains to check that 
$\gamma \in \mathscr{A}(\delta,E,\eta,j_\pm,s,\omega)$.
The lower-semicontinuity 
of $\TP_{q,\delta}$ established in Lemma \ref{lower_sc} implies
$\TP_{q,\delta}(\gamma)\le E$. 
  As a $C^1$-limit of the $\g_k$
 the curve $\g$ satisfies $|\g'|=1$ on $\R/\Z$ and $0\in\g(\R/\Z)$.
 Since the sequences $(u_i^k)_{k \in \N}$ and  $(v_i^k)_{k \in \N}$
 with $\gamma_k(u_i^k)=\gamma_k(u_i^k)$ 
 are bounded in $\Sc$, there exist subsequences, 
 again denoted with index $k$, such that 
 $u_i^k \to u_i$ and $v_i^k \to v_i$ for some $u_i,v_i\in\R/\Z$ as $k \to \infty$.
 The $C^1$-convergence then implies $\gamma(u_i)=\gamma(v_i)$ for all $i=1,\ldots,n$. Recall from \eqref{def_T_gamma} that $T(\gamma_k)=
 \{u_1^k,v_1^k,\dots,u_n^k,v_n^k \}$. Define the set $\tilde{T}:=
 \{u_1,v_1,\dots,u_n,v_n \}$ and observe
\begin{align}
2\eta\le \Lambda(\gamma_k)=&\min\{|a-b|_{\Sc} \; : \; a,b 
\in T(\gamma_k), a \neq b \} \notag\\
\to &\min\{|a-b|_{\Sc} \; : \; a,b \in \tilde{T}, 
a \neq b \}\quad\textnormal{
as $k\to\infty$.}\label{eq:tildeT-parameterabstand}
\end{align}
Hence, the set $\tilde{T}$ contains $2n$ distinct elements and $\gamma$ has at least $n$ self-intersections. 
The affine linearity of $\g_k|_{B_\eta(u_i^k)}$ and of
$\g_k|_{B_\eta(v_i^k)}$  transfers via the
$C^1$-convergence as $k\to\infty$
to affine linearity of $\g|_{B_\eta(u_i)}$ and of
$\g|_{B_\eta(v_i)}$ for $i=1,\ldots,n$.
It still needs to be shown that $\gamma$ has \emph{exactly} $n$ 
self-intersections 
and that every self-intersection is transverse. From 
\cref{bilipschitz}, we obtain the uniform bilipschitz estimate
$
|x-y|_{\Sc}\leq c(q,E,\eta)|\gamma_k(x)-\gamma_k(y)|
$
for all $(x,y) \in (\Sc)^2 \setminus  Y_\delta(\gamma_k)$ and $k \in \N$. 
Passing to the limit $k \to \infty$ yields 
\begin{align}
  \label{bilip_gamma_compactness}
  |x-y|_{\Sc}\leq c(q,E,\eta)|\gamma(x)-\gamma(y) | 
\end{align}
for all $(x,y) \notin 
\tilde{M}:=\bigcup_{i=1}^n \big(B_\delta(u_i)\times
B_\delta(v_i)\big)\cup \big(B_\delta(v_i)\times B_\delta(u_i)\big)$, 
where
as before we denote $B_r(x):=\{y\in\R/\Z:|y-x|_{\Sc}<r\}$.
Note that for $i=1,..,n$ we have
$
|u_i-v_i |_{\Sc}=\lim_{k \to \infty} |u_i^k-v_i^k |_{\Sc}
\geq \lim_{k \to \infty} \Lambda(\gamma_k) \ge 2\eta 
$
by definition of the set $\HF(\eta,j_\pm,s,\omega)$. 
With the bilipschitz estimate  \eqref{bilip_gamma_compactness} 
we infer 
\begin{align*}
c(q,E,\eta)|\gamma(u_i+\delta)-\gamma(v_i+\delta)|&\geq 
|u_i-v_i|_{\Sc}
\geq 2\eta\Fo i=1,\ldots,n.
\end{align*}
Thus, since $\delta\le\frac{\eta}2$ we can use the fact that 
$\g$ is linear
near each self-intersection point to 
 estimate the intersection angle by
\begin{align}
\textstyle\ANG
\big(\gamma'(u_i),\gamma'(v_i)\big)&
\textstyle =2\arcsin\Big( \frac{|\gamma(u_i+\delta)-\gamma(v_i+\delta)|}{2\delta}\Big) 
\textstyle\geq 2\arcsin\Big(\frac{\eta}{\delta c(q,E,\eta)} \Big)
\notag\\
&
\textstyle\geq 2\arcsin\Big(\frac{2}{c(q,E,\eta)} \Big)>0. \label{estimate_angle}
\end{align}
By means of \eqref{bilip_gamma_compactness}
the only possible self-intersections of $\gamma$ other than at the parameter pairs
$(u_i,v_i)$, $i=1,\ldots,n$, 
would be at  $(x,y) \in \tilde{M}$. Without loss of generality we assume 
$(x,y)\in B_\delta(u_i)\times B_\delta(v_i)\setminus\{(u_i,v_i)\}$ 
for some $i=1,\dots,n$, 
otherwise we exchange the roles of $x$ and $y$. 
Due to linearity of $\g$ near the self-intersection at $(u_i,v_i)$, 
this is only  possible if the two lines lie on one straight line. 
However, this is prevented by \eqref{estimate_angle}. 
Hence, $\gamma$ has not more than and therefore exactly 
$n$ self-intersections at $(u_i,v_i)$ for $i=1,\ldots,n$, 
with minimal parameter
distance $\Lambda(\g)\ge 2\eta$ according to
\eqref{eq:tildeT-parameterabstand},
and all these self-intersections are transverse.
In particular,  $\gamma$ is a generic curve and therefore its Arnold invariants 
are well-defined. Since the Arnold invariants are locally constant under 
$C^1$-convergence, we have $J^+(\gamma)=j_+$, $J^-(\gamma)=j_-$ and $\st(\gamma)=s$, 
as these values were prescribed on the curves $\gamma_k$ for each $k\in\N$. 
By continuity of the mapping degree we obtain for
the winding number $W(\gamma)=\lim_{k \to \infty}W(\gamma_k)=
\omega$. 
Thus, we have proven $\gamma \in \HF(\eta,j_\pm,s,\omega)$. 
\end{proof}
\end{theorem}

We have now all tools together to prove existence of minimizers of 
the truncated tangent-point energy $\TP_{q,\delta}$ in the set of admissible curves. 

\textit{Proof of \cref{ex_min}: }
The infimum $\mathcal{I}
:=\inf_{\HF(\eta,j_\pm,s,\omega)\cap\HC}
\TP_{q,\delta}$ is non-negative, since the energy is.
Moreover, $\mathcal{I}$ is finite, because 
the admissibility class $\HF(\eta,j_\pm,s,\omega)\cap\HC$ contains
a $C^{1,1}$-curve by Theorem \ref{thm:admissible-curves}, for which
the energy is finite according to \cref{Sobolev_endliche_Energie}
since $C^{1,1}(\R/\Z,\R^2)\subset W^{2-\frac1{q},q}(\R/\Z,\R^2)$. 
Let $(\gamma_k)_{k \in \N} \subset \HF(\eta,j_\pm,s,\omega)\cap\HC$ 
be a minimizing sequence such that $\lim_{k \to \infty} 
\TP_{q,\delta}(\gamma_k)=\mathcal{I}$. Hence, there exists a constant $E>0$ such that $
\TP_{q,\delta}(\gamma_k)\leq E$ for all $k$ sufficiently large.
Since $\TP_{q,\delta}$ and the admissibility class are
invariant under translations we may assume that
$0\in\g_k(\R/\Z)$ for all $k\in\N.$
By the compactness result of \cref{compactness} there exists a 
subsequence, 
again denoted by  $(\gamma_k)_{k \in \N}$, and a limit curve 
$\gamma_\delta^\eta \in \HF(\eta,j_\pm,s,\omega)$ with 
$\gamma_k \to \gamma_\delta^\eta$ in $C^1$ as $k \to \infty$. 
In addition, $\g_\delta^\eta$ is a generic immersion
contained in some 
$(j_\pm,s,\omega)$-compartment $\tilde{\HC}$,  
which is open with respect
to the $C^1$-topology, so that $\g_k\in\tilde{\HC}$ for
$k\gg 1$; hence $\tilde{\HC}=\HC$.
Applying \cref{lower_sc} for $\Lambda_0:=2\eta\le\Lambda(\g_k)$,
we deduce
$
\mathcal{I}\le\TP_{q,\delta}(\gamma_\delta^\eta )\leq 
\liminf_{k \to \infty} \TP_{q,\delta}(\gamma_k)=\mathcal{I}.
$
The $W^{2-\frac1{q},q}$-regularity of $\g_\delta^\eta$ now
follows from Lemma \ref{regularity}.
\qed

\section{The limit $\delta \to 0$ and an optimally immersed curve}
\label{sending_delta_to_zero}
\subsection{The renormalized energy and Gamma convergence.}
\label{sec:5.1}
In general, the truncated energy
$\TP_{q,\delta}$ does not allow for the limiting 
process $\delta$ to zero since the limit energy
would be infinite due to the self-intersections. To take this blow-up 
into account, we scale the energy by a factor following an idea
of Dennis Kube \cite{dennis}, who 
worked with  a logarithmic rescaling of  a suitably truncated
Möbius energy 
on the space of figure-eight shaped curves with only one
self-intersection and vanishing winding number.
To deduce the correct scaling  
factor for the tangent-point energy we 
focus on the dominating interactions of different curve strands near self-intersection
points to study the energy's blow-up behavior. For that
 we restrict the tangent-point energy to annular regions of
   arclength
   near self-intersections of curves in the class $\HF_n$ defined in
   Definition \ref{def:Fn}. 
\begin{definition}\label{def:annular}
Let $q>0$, $n \in \N$, $0<\delta<\theta< \eta$ and 
$\gamma \in \HF_n$. Define the 
\emph{$(\theta,\delta)$-annular tangent-point energy}
$A_{q,\theta,\delta}$ of $\g$ as
\begin{align}
A_{q,\theta,\delta}(\gamma):=\iint \limits_{Y_\theta(\gamma) 
\setminus Y_\delta(\gamma)} 
\left(\frac{2\dist(\ell(x),\gamma(y))}{|\gamma(x)-\gamma(y)|^2} 
\right)^q dx dy=\TP_{q,\delta}(\g)
-\TP_{q,\theta}(\g).\label{eq:def-annular}
\end{align}
\end{definition}
The following representation of the annular tangent-point
energy on curves in the subclass
 \begin{equation}\label{eq:Fneta}
\HF_{n,\eta}:=\{\g\in\HF_n: \textnormal{$\gamma$ is linear within
 arclength $\eta$ of each self-intersection}\}
 \end{equation}
is of central significance. 
\begin{proposition}
\label{gamma_c_min}
Let $q \geq 2$, $n \in \N$, $0<\delta<\theta< \eta $, and
 $\gamma \in \HF_{n, \eta}$. Then 
\begin{align}
  A_{q,\theta,\delta}(\gamma) & =\frac{2^{q+2}}{q-2}\left( \delta^{2-q}-\theta^{2-q} \right)\sum_{i=1}^{n}F(\alpha_i(\gamma)),
    \label{eq:annular-blow-up}
\end{align}
where  $\alpha_i(\g)$ is the angle between the two
tangent lines of
$\g$ at the $i$-th self-intersection for $i\in\{1,\ldots,n\}$, and
\begin{align}
  \begin{split}
    \label{def_F_angle}
    F(\alpha):={}&\textstyle
    \frac{1}{(\sin\alpha)^{q-1}}\Big[
    \int_{-\cot\alpha}^{\tan\frac{\alpha}{2}} 
    \Big(\big(\frac{1}{z^2+1}\big)^q 
    +\big(\frac{z\sin\alpha+\cos\alpha}{z^2+1} \big)^q\Big)dz \\
&\textstyle
+ \int_{\cot\alpha}^{\cot\frac{\alpha}{2}} 
\Big(\big(\frac{1}{z^2+1}\big)^q 
+\big(\frac{z\sin\alpha-\cos\alpha}{z^2+1} \big)^q\Big)dz
\Big]\quad\Fo \alpha \in \left(0,\frac{\pi}{2} \right].
  \end{split}
\end{align}
\end{proposition}
Before proving this proposition let us draw some immediate conclusions.
We show in Lemma \ref{angle_energy_monotone} in the appendix
that the map $F$ defined in \eqref{def_F_angle}
is strictly
monotone decreasing and minimal for $\alpha=\frac{\pi}2$.
\begin{corollary}\label{cor:minimal-right-angle}
A curve $\g\in\HF_{n,\eta}$ minimizes the annular tangent-point energy
$A_{q,\theta,\delta}$
if and only if $\alpha_i(\g)=\frac{\pi}2$ for all $i=1,\ldots,n$.
\end{corollary}
Observe as a direct consequence of the representation
\eqref{eq:annular-blow-up} the following
expression for the limit of the rescaled annular tangent-point
energy as $\delta\to 0$:
\begin{align}
  \label{eq:lim-annular}
\lim_{\delta\to 0}\delta^{q-2}A_{q,\theta,\delta}(\g) & =\frac{2^{q+2}}{q-2}\sum_{i=1}^nF(\alpha_i(\g))\quad\Foa \g\in\HF_{n,\eta}.
\end{align}
This motivates the renormalization
of the truncated tangent-point energy as defined in the introduction.
\begin{definition}
  \label{def_angle_energy}
 For  $q>2$, $n \in \N$, $\eta\in (0,1)$, 
and $\gamma \in \HF_{n,\eta}$ define the 
 \emph{renormalized 
 tangent-point energy}
  \begin{align}\label{eq:Rq}
  R_q(\gamma):=\lim_{\delta \to 0} \delta^{q-2} \TP_{q,\delta}(\gamma) 
  \in [0,\infty].
  \end{align}
  \end{definition}
Note that the definition of the annular energy \eqref{eq:def-annular} 
implies the relation
\begin{align}
\label{relation_energies}
\TP_{q,\theta}(\gamma)=\TP_{q,\delta}(\gamma)
-A_{q,\theta,\delta}(\gamma)\Fo 0<\delta<\theta
<\eta, \,\g\in\HF_{n,\eta}.
\end{align}
This together with
 the additional assumption that $\TP_{q,\frac{\eta}2}(\gamma)<\infty$ 
 leads to
  \begin{align}
  \label{alternativ_limit}
 R_q(\gamma)&=\lim_{\delta \to 0}\delta^{q-2} \TP_{q,\delta}(\gamma)
 \overset{\eqref{relation_energies}}{=}
 \lim_{\delta \to 0}\delta^{q-2} (\TP_{q, \frac{\eta}2}(\gamma) 
 + A_{q,\frac{\eta}{2},\delta}(\gamma) )\notag\\
  &=\lim_{\delta \to 0} \delta^{q-2} A_{q,\frac{\eta}2,\delta}(\gamma)
  \overset{\eqref{eq:lim-annular}}{=}\frac{2^{q+2}}{q-2}\sum_{i=1}^n F(\alpha_i(\gamma)).
  \end{align}
  Thus, we have established an explicit formula for the
  renormalized tangent-point energy $R_q$ 
  in terms of the function $F$, if the curve $\g$ has sufficient 
  regularity as stated in the
  following lemma.
\begin{lemma}
\label{limit_energy}
For $q>2$, $n \in \N$, $\eta\in (0,1)$, and $\gamma \in
W^{2-\frac{1}{q},q}(\R/\Z,\R^2 ) \cap  \HF_{n,\eta}$ one has 
\begin{align*}
R_q(\gamma)=\frac{2^{q+2}}{q-2}\sum_{i=1}^n F(\alpha_i(\gamma)).
\end{align*}
Moreover,  $\gamma \in
W^{2-\frac{1}{q},q}(\R/\Z,\R^2 ) \cap  \HF_{n,\eta}$ is a global
minimizer of $R_q$ if and only if  
$\alpha_i(\gamma)=\frac{\pi}{2}$ for all $i=1,\dots,n$.
\begin{proof}
By \cref{Sobolev_endliche_Energie}, the 
assumed fractional Sobolev regularity ensures that 
$\TP_{q,\frac{\eta}2}(\gamma)$ is finite, so we can apply 
\eqref{alternativ_limit}.
The last statement follows from \cref{angle_energy_monotone} in the
appendix.
\end{proof}
\end{lemma}

In \cref{def_angle_energy}, the renormalized tangent-point energy 
$R_q$ is defined as a pointwise limit. However, one can prove more: 
The energy is in fact the $\Gamma$-limit of the scaled truncated tangent-point
energies $(\delta^{q-2}  \TP_{q,\delta}  )_{\delta>0}$. In order to show
this stronger convergence result stated in \cref{gamma_convergence}, 
we need to verify the $\liminf$ and $\limsup$ inequalities, see 
\cite[Definition 1.5]{gamma_beginner}.

\textit{Proof of \cref{gamma_convergence}:}
We start with 
the $\liminf$-inequality. Let $(\xi_\delta)_\delta$ be a sequence of 
curves contained in the
space $
 \mathcal{W}:= \HF(\eta,j_\pm,s,\omega)\cap W^{2-\frac{1}{q},q}
 (\R/\Z,\R^2)\cap\HC$ such that $\xi_\delta \to \xi$ in $C^1$ as 
 $\delta \to 0$, and with $\xi\in\mathcal{W}$ as well. The latter 
 implies by means of
 \cref{Sobolev_endliche_Energie} that
   $\TP_{q,\frac{\eta}2}(\xi)<\infty$. 
   By virtue of \eqref{relation_energies} for $\theta:=\frac{\eta}2$
   one has
\begin{align*}
\TP_{q,\delta}(\xi_\delta)=\TP_{q,\frac{\eta}2}(\xi_\delta)+
A_{q,\frac{\eta}2,\delta}(\xi_\delta)\quad\Fo 0<\delta <
\frac{\eta}2.
\end{align*}
Then we obtain by \eqref{eq:annular-blow-up} in \cref{gamma_c_min}
\begin{align*}
\delta^{q-2}\TP_{q,\delta}(\xi_\delta)&
\geq \delta^{q-2} A_{q,\frac{\eta}2,\delta}(\xi_\delta)
\overset{\eqref{eq:annular-blow-up}}{=}\textstyle
\frac{2^{q+2}}{q-2}\big(1-\big(\frac{2\delta}{\eta}\big)^{q-2}\big)
\sum_{i=1}^n F\left(\alpha_i\left(\xi_\delta\right)\right),
\end{align*}
where $n:=j_+-j_-$ by 
\eqref{connection_AI_number_self_intersections}.
Due to $\xi\in\HW$ and the $C^1$-convergence $\xi_\delta\to\xi$
we have 
$\alpha_i(\xi_\delta) \to \alpha_i(\xi)>0$ 
as $\delta \to 0$ for every 
$i=1,\ldots,n$. By the continuity of $F$ on $(0,\frac{\pi}2]$
this leads to 
$F(\alpha_i(\xi_\delta)) \to F(\alpha_i(\xi))$ as $\delta \to 0$
for $i=1,\ldots,n$. Hence, applying \cref{limit_energy}
to $\xi\in\HW$,
\begin{align*}
\liminf_{\delta \to 0} \delta^{q-2}\TP_{q,\delta}(\xi_\delta) 
\geq\frac{2^{q+2}}{q-2}\sum_{i=1}^n F(\alpha_i(\xi)) =R_q(\xi).
\end{align*}
It remains to prove the $\limsup$-inequality. Let $\xi\in\mathcal{W}$.
 As a recovery sequence, we simply choose the constant sequence 
 $\xi_\delta =\xi$ for every $\delta>0$. Then by definition of 
 $R_q$ as a pointwise limit we infer  
 $
\lim_{\delta \to 0} \delta^{q-2}\TP_{q,\delta}(\xi_\delta) =
\lim_{\delta \to 0} \delta^{q-2}\TP_{q,\delta}(\xi) 
=R_q(\xi).$\hfill $\Box$

It remains to prove the central  representation formula 
\eqref{eq:annular-blow-up}.

\begin{figure}
  \centering
    \begin{tikzpicture}
    \draw[-] (-3,1.5) -- (-1.5,0.75);
    \draw[-] (-1.5,0.75) -- (0,0);
    \draw[-] (0,0) -- (1.5,-0.75);
    \draw[-] (1.5,-0.75) -- (3,-1.5);
    \draw[-] (-3,-1.5) -- (-1.5,-0.75);
    \draw[-] (-1.5,-0.75) -- (0,0);
    \draw[-] (0,0) -- (1.5,0.75);
    \draw[-] (1.5,0.75) -- (3,1.5);
    \draw (-2.25,1.125) node[red,above right,scale=0.75] {$1$};
    \draw (-0.75,0.375) node[red,above right,scale=0.75] {$2$};
    \draw (0.75,-0.375) node[red,below left,scale=0.75] {$3$};
    \draw (2.25,-1.125) node[red,below left,scale=0.75] {$4$};
    \draw (2.25,1.125) node[red,above left,scale=0.75] {$5$};
    \draw (0.75,0.375) node[red,above left,scale=0.75] {$6$};
    \draw (-0.75,-0.375) node[red,below right,scale=0.75] {$7$};
    \draw (-2.25,-1.125) node[red,below right,scale=0.75] {$8$};
    \draw[black] (0.89,0.45) arc (26.57:-26.57:1);
    \draw (0.3,0) node[black,right,scale=0.75] {$\alpha_i$};
    \draw (-3,1.5) node[black,above,scale=0.75] {$\gamma(u_i-\theta)$};
    \draw (3,-1.5) node[black,below,scale=0.75] {$\gamma(u_i+\theta)$};
    \draw (-3,-1.5) node[black,below,scale=0.75] {$\gamma(v_i+\theta)$};
    \draw (3,1.5) node[black,above,scale=0.75] {$\gamma(v_i-\theta)$};
      \draw (-1.5,0.75) node[black,above right,scale=0.75] {$\gamma(u_i-\delta)$};
    \draw (1.5,-0.75) node[black,below left,scale=0.75] {$\gamma(u_i+\delta)$};
    \draw (-1.5,-0.75) node[black,below right,scale=0.75] {$\gamma(v_i+\delta)$};
    \draw (1.5,0.75) node[black,above left,scale=0.75] {$\gamma(v_i-\delta)$};
    \fill[black] (-3,1.5) circle (2pt);
    \fill[black] (-3,-1.5) circle (2pt);
    \fill[black] (3,1.5) circle (2pt);
    \fill[black] (3,-1.5) circle (2pt);
    \fill[black] (0,0) circle (2pt);
    \fill[black] (-1.5,0.75) circle (2pt);
    \fill[black] (-1.5,-0.75) circle (2pt);
    \fill[black] (1.5,0.75) circle (2pt);
    \fill[black] (1.5,-0.75) circle (2pt);
    \fill[black] (0,0) circle (2pt);
    \draw (-0.2,0) node[black,left,scale=0.75] {$\gamma(u_i)=\gamma(v_i)$};
    \end{tikzpicture}  
    \caption{The local situation around the $i$-th 
    self-intersection of $\gamma$: 
    The numbers $k=1,\ldots,8$ (in red)} at the segments correspond 
    to the set $J_k^i$.
    \label{angle_energy_local_situation}
  \end{figure}
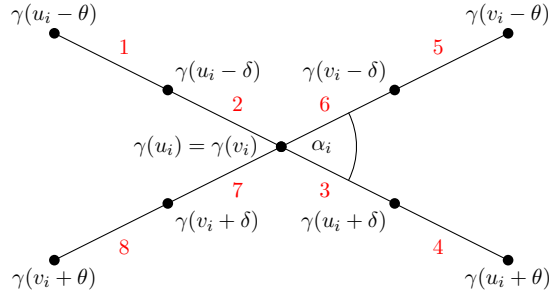

{\it Proof of Proposition \ref{gamma_c_min}.}\,
Let $\gamma \in \HF_{n,\eta}$ and let $u_i,v_i \in \Sc$ with 
$u_i \neq v_i$ be the distinct points with $\gamma(u_i)=\gamma(v_i)$ 
for every $i=1,\dots,n$. Then define for $i=1,\ldots,n$ the 
intervals 
\begin{alignat*}{4}
J_1^i&=[u_i-\theta,u_i-\delta], \qquad &J_5^i
&=[v_i-\theta,v_i-\delta],\\
J_2^i&=[u_i-\delta,u_i], \qquad &J_6^i&=[v_i-\delta,v_i],\\
J_3^i&=[u_i,u_i+\delta], \qquad &J_7^i&=[v_i,v_i+\delta],\\
J_4^i&=[u_i+\delta,u_i+\theta], \qquad &J_8^i&=[v_i+\delta,v_i+\theta].
\end{alignat*}
Denote by $\HE_{kl}^i$ the energy integral over  
$J_k^i \times J_l^i$ for $k,l \in \{1,\dots 8 \}$. Notice
that  there is no interaction between different parts of a single
straight segment, and that also the interactions $\HE^i_{26}$
and
$\HE^i_{27}$ and its symmetric counterparts between the 
central portions at the self-intersections vanish for each 
$i=1,\ldots,
n$ since the truncated annular energy does not see these portions
by definition. Therefore, the 
truncated annular energy can be written as 
\begin{align*}
A_{q,\theta,\delta}(\gamma)&=\textstyle
\sum_{i=1}^n 
\big(\HE_{15}^i+\HE_{16}^i+\HE_{17}^i+\HE_{18}^i+\HE_{25}^i+\HE_{28}^i+\HE_{35}^i+\HE_{38}^i+\HE_{45}^i+\HE_{46}^i+\HE_{47}^i+\HE_{48}^i\\
&\qquad \textstyle
+\HE_{51}^i+\HE_{52}^i+\HE_{53}^i+\HE_{54}^i+\HE_{61}^i+\HE_{64}^i+\HE_{71}^i+\HE_{74}^i+\HE_{81}^i+\HE_{82}^i+\HE_{83}^i+\HE_{84}^i \big)\\ 
&
\textstyle
=4\sum_{i=1}^n \big(\HE_{15}^i+\HE_{16}^i+\HE_{25}^i+\HE_{17}^i+\HE_{18}^i+\HE_{28}^i\big),
\end{align*}
where the last equality follows from the symmetric structure 
around each self-intersection. Since 
$A_{q,\theta,\delta}$ takes only the linear segments into account, we 
can rewrite it as a function that solely depends on the 
intersection angles $\alpha_i:=\alpha_i(\g)$, $i=1,\ldots,n$.
The local situation is shown in \cref{angle_energy_local_situation}. 
Consider first the energy
integrals $\HE_{15}^i+\HE_{16}^i$. 
For  $(x,y) \in [u_i-\theta,u_i-\delta]\times[v_i-\theta,v_i]$
set 
\begin{equation}\label{eq:s-t}
s(x):=|\g(x)-\g(u_i)|\quad\AND\quad
t(y):=|\g(y)-\g(v_i)|,
\end{equation}
and infer by basic trigonometry
$
\dist(\ell(x),\gamma(y))=t(y)\sin\alpha_i$
and 
\begin{align*}
|\gamma(x)-\gamma(y)|^2=&s^2(x)+t^2(y)-2s(x)t(y)\cos(\pi-\alpha_i)=
s^2(x)+t^2(y)-2s(x)t(y)\cos\alpha_i.
\end{align*}
Therefore,
\begin{align*}
\HE_{15}^i+\HE_{16}^i
&=\textstyle
2^q\int_{u_i-\theta}^{u_i-\delta} 
\int_{v_i-\theta}^{v_i} 
\big(\frac{t(y)\sin\alpha_i}{s^2(x)+t^2(y)-2s(x)t(y)\cos\alpha_i}
\big)^q\,dydx\\
&\textstyle
=2^q\int_\delta^\theta\int_0^\theta
\big(\frac{t\sin\alpha_i}{s^2+t^2+2st\cos\alpha_i}\big)^q\,dtds,
\end{align*}
where we changed variables in both integrals using
$$
t'(y)=\textstyle
\frac{\langle \gamma(y)-\gamma(v_i) , \gamma'(y)\rangle}{|\gamma(y)
-\gamma(v_i)|}=\big 
\langle \frac{\gamma(y)-\gamma(v_i)}{|\gamma(y)-\gamma(v_i)|}, \gamma'(y)
\big\rangle=-\langle \gamma'(y),\gamma'(y) \rangle
=-1,
$$
and analogously, $s'(x)=-1$,  
since the vectors
$\gamma(y)-\gamma(v_i)$ and
$\g(x)-\g(u_i)$ are contained in their respective 
linear segment and point in the opposite direction of the tangent-vector.
Furthermore, $t(v_i-\theta)=\theta$ and $t(v_i)=0$, as well as
$s(u_i-\theta)=\theta$ and $s(u_i-\delta)=\delta$,
due to the linear structure of the curve $\gamma$ and its arclength 
parametrization.
For $\HE_{25}^i$ we obtain in the same manner
\begin{align*}
\HE_{25}^i\textstyle
=2^q \int_0^\delta \int_\delta^\theta 
\big( \frac{t\sin\alpha_i}{t^2+s^2+2st\cos\alpha_i}\big)^q dt ds.
\end{align*}
We treat the remaining integrals similarly, 
where the angle $\alpha_i$ instead of $\pi- \alpha_i$ appears in the 
Euclidean distance $|\g(x)-\g(y)|$ for $y\in [v_i,v_i+\theta]$
and $x\in 
[u_i-\theta,u_i-\delta]$, or $x\in [u_i-\delta,u_i]$, respectively,
to obtain
\begin{align*}
\HE_{17}^i+\HE_{18}^i+\HE_{28}^i={}&\textstyle
2^q \int_\delta^\theta \int_0^\theta 
\big( \frac{t\sin\alpha_i}{t^2+s^2-2st\cos\alpha_i}\big)^q dt ds
+2^q \int_0^\delta \int_\delta^\theta 
\big( \frac{t\sin\alpha_i}{t^2+s^2-2st\cos\alpha_i}\big)^q dt ds.
\end{align*}
Abbreviating $
h(\alpha,s,t):=(t\sin(\alpha))^q \big((t^2+s^2+2st\cos\alpha)^{-q}
+(t^2+s^2-2st\cos\alpha)^{-q} \big)$
we can rewrite the truncated annular 
energy as
\begin{align*}
A_{q,\theta,\delta}(\gamma)={}&\textstyle
2^{q+2} \sum_{i=1}^n
\big( \int_0^\delta \int_\delta^\theta h(\alpha_i,s,t) dt ds 
+ \int_\delta^\theta \int_0^\theta h(\alpha_i,s,t) dt ds\big)\\
={}&\textstyle
2^{q+2} \sum_{i=1}^n \iint\limits_{Q_{\theta,\delta} }
h(\alpha_i,s,t) dt ds, 
\end{align*}
where 
$Q_{\theta,\delta}:=[0,\theta]^2 \setminus [0,\delta]^2$.  
Splitting the integration domain along the diagonal and
using the transformation formula for the transformation
$T:Q_{\theta,\delta}\cap\{t\ge s\}
\to Q_{\theta,\delta}\cap\{t\le s\}$ defined by $J(t,s):=(s,t)$
leads to
\begin{align*}
A_{q,\theta,\delta}(\gamma)
&\textstyle
=2^{q+2}\sum_{i=1}^n \iint\limits_{Q_{\theta,\delta} 
\cap\{t \leq s \}} [h(\alpha_i,s,t) + h(\alpha_i,t,s)] dt ds\\
&\textstyle
=2^{q+2} \sum_{i=1}^n \int_\delta^\theta \int_0^s 
[h(\alpha_i,s,t) + h(\alpha_i,t,s)] dt ds.
\end{align*}
To compute the inner integral consider the two functions
\begin{align*}
\textstyle
H_s(\alpha):=\int_0^s h(\alpha,t,s) dt,\quad
G_s(\alpha):=\int_0^s h(\alpha,s,t) dt, \quad \Fo s\in [\delta,\theta].
\end{align*}
Let us first rewrite  $H_s$ as
\begin{align*}
H_s(\alpha)
={}&\textstyle
\frac{1}{(s\sin\alpha)^q} 
\Big(\int_0^s \big[\left(\frac{t+s\cos\alpha}{s\sin\alpha}
\right)^2+1\big]^{-q} dt 
+ \int_0^s \big[\left(\frac{t-s\cos\alpha}{s\sin\alpha}
\right)^2+1\big]^{-q} dt \Big).
\end{align*}
Substituting $z(t)=\frac{t+s\cos\alpha}{s\sin\alpha}$ 
in the first integral and $z(t)=\frac{t-s\cos\alpha}{s\sin\alpha}$ 
in the second one yields
\begin{align*}
H_s(\alpha)=\textstyle
\frac{1}{(s\sin\alpha)^{q-1}} 
\left( \int_{\cot\alpha}^{\cot\frac{\alpha}{2}}
\big( \frac{1}{z^2+1}\big)^q dz 
+ \int_{-\cot\alpha}^{\tan\frac{\alpha}{2}}
\big( \frac{1}{z^2+1}\big)^q dz\right).
\end{align*}
Similarly, one computes
\begin{align*}
G_s(\alpha) 
={}&\textstyle
\frac{1}{(s\sin\alpha)^{q-1}}
\left(\int_{\cot\alpha}^{\cot\frac{\alpha}{2}} 
\big(\frac{z\sin\alpha-\cos\alpha}{z^2+1}\big)^q dz 
+ \int_{-\cot\alpha}^{\tan\frac{\alpha}{2}} 
\big(\frac{z\sin\alpha+\cos\alpha}{z^2+1}\big)^q dz \right).
\end{align*}
Together with the formula of $F$ from \eqref{def_F_angle} this yields 
\begin{align*}
  A_{q,\theta,\delta}(\gamma)\textstyle
  =2^{q+2}\sum_{i=1}^n \int_\delta^\theta 
  s^{1-q}F(\alpha_i(\gamma))ds
  =\frac{2^{q+2}}{q-2}\left( \delta^{2-q}-\theta^{2-q} \right)
  \sum_{i=1}^{n}F(\alpha_i(\gamma)),
\end{align*}
which proves \eqref{eq:annular-blow-up}.
\hfill$\Box$

\subsection{
Optimally immersed curves with prescribed Arnold invariants.}
\label{sec:5.2}

It remains to be investigated whether a sequence of minimizers 
$(\gamma_\delta^\eta)_{\delta>0}$ of the truncated
energies $\TP_{q,\delta}$ converges in $C^1$ if $\delta$ is sent to zero.
So far, we do not have a uniform bound on the energy values 
$\TP_{q,\delta}(\gamma_\delta^\eta)$ independent of $\delta$. 
Hence, our compactness result \cref{compactness} is not directly 
applicable. However, we will be able to deduce a uniform bound on 
$\TP_{q,\theta}( \gamma_\delta^\eta)$ independent of 
$\delta\in
(0,\frac{\theta}2]$ and of $\eta$ as long as $\theta \le\eta$, 
see \cref{energy_uniform_bound}. This will allow us to obtain 
subconvergence of the  
minimizers $\g_\delta^\eta$ to  a limit curve 
$\Gamma^{\eta}$ as $\delta\to 0$. 
Properties of the curve $\Gamma^{\eta}$ are then 
investigated at the end of this section. 

To establish 
the uniform bound on $\TP_{q,\theta}(\gamma_\delta^\eta)$ 
we use identity 
\eqref{relation_energies}. Since the curves $\gamma_\delta^\eta$ are 
minimizers, the energy value
$\TP_{q,\delta}(\gamma_c)$ of any admissible comparison curve
$\g_c\in\HF(\eta,j_\pm,s,\omega)$ gives an upper bound for the first 
summand on 
 the right-hand side of  \eqref{relation_energies}. 
 In addition, Corollary \ref{cor:minimal-right-angle} implies
 that any comparison
 curve $\gamma_c\in\HF(\eta,j_\pm,s,\omega)\subset
 \HF_{n,\eta}$ whose intersection angles $\alpha_i(\g_c)$ equal
 $\frac{\pi}2$ actually \emph{minimizes} the second 
 summand on the right-hand side of \eqref{relation_energies}.
 Here is the precise statement.
\begin{corollary}
\label{energy_uniform_bound}
For $j_+,j_-,s,\omega \in \Z$ consider
the compartment $\HC=\HC(j_\pm,s,\omega)$, and
fix any $\eta\in (0,\eta_1(\HC)]$. Then for any $q>2$ and
$\theta\in (0,\eta)$ there is a constant
$c=c(\theta,q,j_+,j_-,s,\omega)>0$ such that
\begin{align*}
\TP_{q,\theta}(\gamma_\delta^\eta)\leq c
\quad\Foa 0<\delta\le\frac{\theta}2,
\end{align*}
 where $\eta_1(\HC)$ is the constant of
 Theorem \ref{thm:admissible-curves}, and
 the curves $\gamma_\delta^\eta$ are the minimizers obtained in
 \cref{ex_min}.
\begin{proof}
\cref{thm:admissible-curves}  for $\eta:=\eta_1(\HC)$
guarantees the existence of a comparison curve 
$\gamma_c \in C^{1,1}(\Sc,\R^2)\cap\HF(\eta_1,j_\pm,s,\omega)
\cap\HC$ whose
intersection angles equal $\frac{\pi}2$, and
which is also contained in $\HF(\eta,j_\pm,s,\omega)$
by means of \eqref{eq:nested-F}. Set
\begin{align*}
c=  c(\theta,q,j_+,j_-,s,\omega):=\TP_{q,\theta}(\gamma_c),
\end{align*}
which by \cref{Sobolev_endliche_Energie} is finite, 
since $C^{1,1}(\Sc,\R^2)\subset W^{2-\frac1{q},q}(\Sc,\R^2)$. (Notice
that $\delta\le\frac{\theta}2 <
\frac{\eta}2\le\frac{\Lambda(\g_c)}4$
by definition of the class $\HF(\eta,j_\pm,s,\omega)$ so that
\cref{Sobolev_endliche_Energie} is applicable.)
By virtue of \eqref{relation_energies}, the fact that $\g_\delta^\eta$
minimizes the truncated tangent-point energy $\TP_{q,\delta}$
in $\HF(\eta,j_\pm,s,\omega)\cap\HC$, and
by
\eqref{eq:annular-blow-up}  in combination with 
\cref{cor:minimal-right-angle}
we infer
\begin{align*}
\TP_{q,\theta}(\gamma_\delta^\eta)
&\overset{\eqref{relation_energies}}{=}\TP_{q,\delta}(\gamma_\delta^\eta)
\!-\!A_{q,\theta,\delta}(\gamma_\delta^\eta)
\overset{\eqref{eq:annular-blow-up}}{\leq}
\TP_{q,\delta}(\gamma_c)\!-\!A_{q,\theta,\delta}(\gamma_c)
\overset{\eqref{relation_energies}}{=}
c. \qedhere
\end{align*}
\end{proof}
\end{corollary}
With this uniform energy bound, we are now able to pass to the limit of a sequence of minimizers. 

\textit{Proof of \cref{limit_delta}:}
Fix $\eta\in (0,\eta_1(\HC)]$, set
$\delta_k:=\frac{1}{k}$ and consider the subsequence
of $\TP_{q,\delta_k}$-minimizing
curves $\gamma_k^\eta:=\gamma_{\delta_k}^\eta\in\HF(\eta,j_\pm,s,\omega)
\cap\HC$
for $k \in \N$ such that $k\ge\frac{4}{\eta}$, 
whose
existence is guaranteed by \cref{ex_min}. By translational
invariance of the energy we may assume that $0\in\g_k^\eta(\R/\Z)$
for all $k\ge\frac4{\eta}$.
By \cref{energy_uniform_bound} for $\theta:=\frac{\eta}2$
there exists a constant $c=c(\theta,q, j_+,j_-,s,\omega)>0$
such that
$
\TP_{q,\frac{\eta}2}(\gamma_k^\eta)
\le c$ for all $k\ge
\frac{4}{\eta}. 
$
Hence, we can apply the compactness result of \cref{compactness} 
for $\delta:=\frac{\eta}2$
to obtain a subsequence, again denoted by 
$(\gamma_k^\eta)_k$, and some curve
$\Gamma^\eta \in \HF(\eta,j_\pm,s,\omega)$,
such that $\gamma_k^\eta \to \Gamma^\eta$ in $C^1$ 
as $k \to \infty$. Since $\G^\eta$ is a generic
immersion in some open $(j_\pm,s,\omega)$-compartment $\tilde{\HC}$
we find $\gamma_k^\eta\in\tilde{\HC}$ for $k\ll 1$ implying
$\tilde{\HC}=\HC$.
The lower semi-continuity result \cref{lower_sc} implies
that $\TP_{q,\frac{\eta}2}(\G^\eta)\le c$, which by means of
\cref{Sobolev_endliche_Energie} implies that $\G^\eta\in
W^{2-\frac1{q},q}(\R/\Z,\R^2).$  
The minimizers $\gamma_k^\eta$ also minimize the scaled energies 
${\delta_k}^{q-2}\TP_{q,\delta_k}$ in $\HF(\eta,j_\pm,s,\omega)
\cap\HC$. 
By \cite[Corollary 7.20]{dal_maso_gamma_convergence},
 the curve 
$\Gamma^\eta$ is then a minimizer of the limit energy $R_q(\cdot)=
\lim_{k\to\infty}\delta_k^{q-2}\TP_{q,\delta_k}(\cdot)$ introduced in
\cref{def_angle_energy}, and $\G^\eta$ satisfies
\eqref{eq:limit-energy-value} in addition.
According to \cref{limit_energy},
the energy $R_q(\cdot)$ is minimized
if and only if every 
intersection angle is equal to $\frac{\pi}{2}$. Hence, 
$\alpha_i\left(\Gamma^\eta\right)=\frac{\pi}{2}$ for every $i=1,\dots,n$. 

It remains to show \eqref{eq:almost-minimal}.
Fix $\varepsilon>0$. As before, we use the relation 
\eqref{relation_energies} for $\theta:=\frac{\eta}2>\delta$
to write
\begin{align*}
0\leq \TP_{q,\delta}(\Gamma^\eta)-\TP_{q,\delta}(\gamma_\delta^\eta)
=\TP_{q,\frac{\eta}2}(\Gamma^\eta)+A_{q,\frac{\eta}2,\delta}(\Gamma^\eta)
-\TP_{q,\frac{\eta}2}(\gamma_\delta^\eta)-A_{q,\frac{\eta}2,\delta}
(\gamma_\delta^\eta).
\end{align*}
Every self-intersection angle of $\G^\eta$ equals $\frac{\pi}2$ so
that $\Gamma^\eta$ minimizes $A_{q,\frac{\eta}2,\delta}$
according to
 \cref{cor:minimal-right-angle}, so that $A_{q,\frac{\eta}2,\delta}(\Gamma^\eta)
 \leq A_{q,\frac{\eta}2,\delta}(\gamma_\delta^\eta)$. We deduce
\begin{align}\label{eq:TP-difference}
0\leq \TP_{q,\delta}(\Gamma^\eta)-\TP_{q,\delta}(\gamma_\delta^\eta) 
\leq \TP_{q,\frac{\eta}2}(\Gamma^\eta) 
-\TP_{q,\frac{\eta}2}(\gamma_\delta^\eta),
\end{align}
and hence
\begin{align*}
0&\leq \liminf_{\delta \to 0} \left( \TP_{q,\delta}
(\Gamma^\eta)-\TP_{q,\delta}(\gamma_\delta^\eta)\right)
\leq \limsup_{\delta \to 0} \left( \TP_{q,\delta}(\Gamma^\eta)
-\TP_{q,\delta}(\gamma_\delta^\eta)\right)\\
&\overset{\eqref{eq:TP-difference}}{\leq} 
\limsup_{\delta \to 0} \big( \TP_{q,\frac{\eta}2}(\Gamma^\eta) 
- \TP_{q,\frac{\eta}2}(\gamma_\delta^\eta)\big)
= \TP_{q,\frac{\eta}2}(\Gamma^\eta)-\liminf_{\delta \to 0} 
\TP_{q,\frac{\eta}2}(\gamma_\delta^\eta)\leq 0. 
\end{align*}
The last inequality holds by the lower semi-continuity 
with respect to $C^1$-convergence, as shown in \cref{lower_sc}. 
Hence, $\lim_{\delta \to 0} (\TP_{q,\delta}(\Gamma^\eta)
-\TP_{q,\delta}(\gamma_\delta^\eta) )=0$, so there exists $\hat{\delta}
=\hat{\delta}(\varepsilon)>0$ such that for all $0<\delta<\hat{\delta}$
\begin{align*}
\inf_{\HF(\eta,j_\pm,s,\omega)\cap\HC} 
\TP_{q,\delta}\le\TP_{q,\delta}(\Gamma^\eta)
&<\TP_{q,\delta}(\gamma_\delta^\eta)
+\varepsilon=\inf_{\HF(\eta,j_\pm,s,\omega)\cap\HC} \TP_{q,\delta} 
+\varepsilon. 
\end{align*}
\qed

If we initially restrict the minimizing process to curves in 
$\HF(\eta,j_\pm,s,\omega)\cap\HC$ that have a right angle at 
each 
self-intersection, the limit curves $\Gamma^\eta$ indeed minimize 
the truncated energy $\TP_{q,\delta}$ 
for every $\delta\le\frac{\eta}2$
among those curves as stated in \cref{cor:optimal-among-right-angels}.

{\it Proof of \cref{cor:optimal-among-right-angels}.}\,
Fix any $0<\delta \le\frac{\eta}{2}$ and take
$\sigma\in (0,\delta)$. Let $\gamma$
be a curve in $\HF(\eta,j_\pm,s,\omega)\cap\HC$ whose
intersection angles equal $\frac{\pi}2$. 
Denote by $\gamma_\sigma^\eta$ a minimizer of $\TP_{q,\sigma}$ 
in $\HF(\eta,j_\pm,s,\omega)\cap\HC$ whose existence is guaranteed by
Theorem \ref{ex_min}. 
Combining the fact that $\gamma_\sigma^\eta$ is a minimizer with 
\cref{cor:minimal-right-angle} 
yields by means of \eqref{relation_energies}
\begin{align}
\label{minimierer_vergleich}
\TP_{q,\delta}( \gamma_\sigma^\eta)
\overset{\eqref{relation_energies}}{=}
\TP_{q,\sigma}(\gamma_\sigma^\eta)
-A_{q,\sigma,\delta}(\gamma_\sigma^\eta)
\leq\TP_{q,\sigma}(\gamma)-
A_{q,\sigma,\delta}(\gamma)
\overset{\eqref{relation_energies}}{=}\TP_{q,\delta}(\gamma).
\end{align}
Furthermore, \cref{limit_delta} guarantees the existence of
a subsequence $(\gamma_{\sigma_k}^\eta)_k$ that converges in 
$C^1$ to $\Gamma^\eta$ for $\sigma_k \to 0$ as $k\to\infty$. 
Using the lower semi-continuity in \cref{lower_sc} and
\eqref{minimierer_vergleich} then yields 
$
\TP_{q,\delta}(\Gamma^\eta)\leq 
\liminf_{k\to\infty} \TP_{q,\delta}(\gamma_{\sigma_k}^\eta)
\le\TP_{q,\delta}(\gamma).
$
\hfill\qed

The same line of arguments also yields the monotonicity of the energy values
$\TP_{q,\theta}(\G^\eta)$ in  the $\eta$-variable.
\begin{corollary}\label{cor:eta-monoton}
The map 
 $\eta\mapsto
\TP_{q,\theta}(\G^\eta)$ is non-decreasing on the interval $(2\theta,\eta_1(\HC)]$.
In particular,  $\lim_{\eta\searrow 2\theta}\TP_{q,\theta}(\G^\eta)\le
\lim_{\eta\nearrow\eta_1(\HC)}\TP_{q,\theta}(\G^\eta)$.
\begin{proof}
 The truncated energies
can be split into two terms according to \eqref{relation_energies}:
\begin{align}
\TP_{q,\theta}(\g^\eta_\delta) &\overset{\eqref{relation_energies}}{=}
\TP_{q,\delta}(\g_\delta^\eta)-A_{q,\theta,\delta}(\g_\delta^\eta)
\label{eq:monotone-in-eta}\\
&\le \TP_{q,\delta}(\G^{\eta_*})-A_{q,\theta,\delta}(\G^{\eta_*})
\overset{\eqref{relation_energies}}{=}\TP_{q,\theta}(\G^{\eta_*})\Fo 2\theta<\eta\le\eta_*\le
\eta_1(\HC),\notag
\end{align}
where we used the minimality of $\g_\delta^\eta$ in $\HF(\eta,j_\pm,s,\omega)\cap\HC$
containing the set $\HF(\eta_*,j_\pm,s,\omega)\cap\HC$ by means of 
\eqref{eq:nested-F},
and \cref{cor:minimal-right-angle}, since $\G^{\eta_*}$ intersects
exclusively in right angles according to \cref{limit_delta}. Combining  \eqref{eq:monotone-in-eta} with the 
lower semicontinuity of $\TP_{q,\theta}$ (\cref{lower_sc})  yields for 
$2\theta<\eta\le\eta_*\le\eta_1(\HC)$ 
$$
\TP_{q,\theta}(\G^\eta)\le
\liminf_{\delta\to 0}\TP_{q,\theta}(\g^\eta_\delta)
\overset{\eqref{eq:monotone-in-eta}}{\le}\liminf_{\delta\to 0}\TP_{q,\theta}(\G^{\eta_*})=\TP_{q,\theta}(\G^{\eta_*}).\qedhere
$$
\end{proof}
\end{corollary}

\appendix
\section{Auxiliary statements}
\label{appendix}
\begin{lemma}[Global parametrization with local graphs]\label{lem:bauser}
Let $k\ge 1$ and  $\g\in C^k(\R/\Z,\R^2)$ with $|\g'|=1$ on $\R/\Z$.
 Then for any
mutually
distinct arclength parameters $x_1,\ldots,x_N\in\R/\Z$ there exists $r_0>0$ and a 
reparametrization $\tilde{\g}\in C^k(\R/\Z,\R^2)$ of $\g$ with the 
same
orientation as $\g$ satisfying
\begin{equation}\label{eq:reparametrized-graphlike-curve}
\tilde{\g}(x)=\begin{cases}
\g(x) & \Fo x\not\in\bigcup_{i=1}^NB_{2r_0}(x_i)\\
\g(x_i)+(x-x_i)\g'(x_i)+u_i(x-x_i)(\g'(x_i))^\perp & \Fo x\in B_{r_0}(x_i),i=1,\ldots,N,
\end{cases}
\end{equation}
where $u_i\in C^k(\R)$ satisfies $0=u_i(0)=u_i'(0)$ for $i=1,\ldots,N.$ In addition, 
$|\tilde{\g}'(x)|\ge\frac14 $ for all $x\in \bigcup_{i=1}^NB_{2r_0}(x_i),$ and
the balls $B_{3r_0}(x_i)$ are mutually disjoint.
\begin{proof}
Choose $r_0>0$ so small that the subintervals $B_{3r_0}(x_i)\subset\R/\Z$ are
mutually disjoint for $i=1,\ldots,N$ and such that,  by uniform continuity of the unit
tangent,
\begin{equation}\label{eq:tangent-deviation}
|\g'(x)-\g'(x_i)|<\frac19\quad\Foa x\in B_{\frac94 r_0}(x_i)\AND i=1,\ldots,N.
\end{equation}
Set $\epsilon:=\frac94 r_0$ for brevity.

{\it Step 1.}\,
For  $i=1$ take the rotation $R\in SO(2)$ satisfying $R\g'(x_1)=
e_1=(1,0)^T\in\R^2$. The curve $\zeta(x):=R(\g(x+x_1)-\g(x_1))$ is
of class $C^k(\R/\Z,\R^2)$ and satisfies $\zeta(0)=0$,  $\zeta'(0)=e_1$ and
$(\zeta'(0))^\perp=e_2=(0,1)^T\in\R^2.$ In addition, $|\zeta'|=1$ on $\R/\Z$ and the derivative of the first component $\zeta_1$ satisfies
\begin{equation}\label{eq:eta-abl}
1\ge\zeta_1'(x)>\frac89\quad\Foa x\in B_{\epsilon}(0)
\end{equation}
by means of \eqref{eq:tangent-deviation},
so that $\zeta_1|_{B_\epsilon(0)}$ is strictly monotonically increasing with value
$\zeta_1(0)=0$. Hence its inverse 
$g:=(\zeta_1|_{B_\epsilon(0)})^{-1}:V_1:=\zeta_1(B_\epsilon(0))\to 
B_\epsilon(0)$ with $\zeta_1(0)=0\in V_1$ exists and satisfies
$g\in C^k(V_1)$, $g(0)=0$, and $(g\circ\zeta_1)(x)=x$ for all $x\in B_\epsilon(0)$ and
$(\zeta_1\circ g)(y)=y$ for all $y\in V_1$. So, $\zeta(B_\epsilon(0))$ can be reparametrized
in an orientation preserving way as a graph of the function 
$u_1:=\zeta_2\circ g\in C^k(V_1)$
according to
\begin{equation}\label{eq:graph}
V_1\ni y\mapsto \left(\begin{array}{c}(\zeta_1\circ g)(y)\\
                                  (\zeta_2\circ g)(y)\end{array}\right)
= \left(\begin{array}{c}y\\
                                  (\zeta_2\circ g)(y)\end{array}\right)
=\left(\begin{array}{c} y\\
u_1(y)\end{array}\right).
\end{equation}
Notice that $B_{\frac89 \epsilon}(0)=B_{2r_0}(0)\subset V_1$ by virtue of \eqref{eq:eta-abl}.

{\it Step 2.}\,
Choose
a cut-off function $\varphi\in C_0^\infty(B_{2r_0}(0))$ (extended to all of $\R$ by zero), such that $0\le\varphi\le 1$ on $\R$, $\varphi\equiv 1$ on $B_{r_0}(0)$ and
$0\le\varphi'\le\frac2{r_0}$ on $[-2r_0,-r_0]$ and $0\ge\varphi'\ge -\frac2{r_0}$ on $[r_0,2r_0]$.
Now define the global orientation preserving reparametrization 
\begin{equation}\label{eq:eta-definition}
\tilde{\zeta}(x):=\begin{cases}
\zeta(x) & \Fo x\in\R/\Z\setminus B_{2r_0}(0)\\
\zeta\big((1-\varphi(x))x+\varphi(x)g(x)\big)& \Fo x\in B_{2r_0}(0),
\end{cases}
\end{equation}
so that $\tilde{\zeta}'$ coincides with $\zeta'$ on $\R/\Z\setminus B_{2r_0}(0)$. Note that $\tilde{\zeta}$ is of class $C^k(\R/\Z,\R^2)$ 
since $\zeta,g$ are in $C^k$, and
$\varphi$ is smooth and has compact support in $B_{2r_0}(0)$.
For $x\in B_{2r_0}(0)$ compute 
\begin{equation}\label{eq:eta-tilde-strich}
\tilde{\zeta}'(x)=\zeta'\big((1-\varphi(x))x+\varphi(x)g(x)\big)\cdot[\varphi'(x)(g(x)-x)+
\varphi(x)(g'(x)-1)+1],
\end{equation}
and notice that \eqref{eq:eta-abl} implies 
\begin{equation}\label{eq:inverse-derivative}
1\le g'(x)=\frac{1}{\zeta_1'(g(x))}<\frac98 \quad\Fo x\in B_{2r_0}(0), 
\end{equation}
and therefore 
$g'-1\ge 0$ on $B_{2r_0}(0)$,   
and $g(x)-x=g(x)-g(0)-x\in [0,\frac18 x)$ for $x\in [0,2r_0)$
as well as $g(x)-x\in (\frac18 x,0]$ for $x\in (-2r_0,0]$.
Hence we can bound the term in square brackets in \eqref{eq:eta-tilde-strich}  
from below by
\begin{align*}
\varphi'(x)(g(x)-x)+1 > -\frac2{r_0}\cdot
\frac18 x+1>\frac12\quad\Fo x\in [0,2r_0)
\end{align*}
and by $\frac2{r_0}\cdot
\frac18 x+1>\frac12$ for $x\in (-2r_0,0]$. This shows that $|\tilde{\zeta}'(x)|>\frac12$ for all $x\in B_{r_0}(0)$. 

{\it Step 3.}\,
Define $\tilde{\g}(x):=\g(x_1)+R^{-1}\tilde{\zeta}(x-x_1)$,
which is of class $C^k(\R/\Z,\R^2)$,
to find by definition
of $\tilde{\zeta}$ in \eqref{eq:eta-definition} $\tilde{\g}(x)=\g(x_1)+
R^{-1}\tilde{\zeta}(x-x_1)=\g(x_1)+R^{-1}\zeta(x-x_1)=\g(x)$ for all $x\not\in 
B_{2r_0}(x_1)$. For $x\in B_{r_0}(x_1)$, on the other hand, we have
\begin{align*}
\tilde{\g}(x) & =\g(x_1)+R^{-1}\big((x-x_1)e_1+u_1(x-x_1)e_2\big)\\
& =  \g(x_1)+(x-x_1)\g'(x_1)+u_1(x-x_1)(\g'(x_1))^\perp
\end{align*}
since $\varphi\equiv 1 $ on $B_{r_0}(0)$ so that the local graph representation 
\eqref{eq:graph} can be used.

{\it Step 4.}\,
Now repeat steps 1--3 for all $i=2,\ldots,N$ to obtain the full statement. Notice that 
the smallness condition on $r_0$ guarantees that the consecutive reparametrizations locally near $x_{i}$ do not affect the previous ones near $x_1,\ldots,x_{i-1}$. 
\end{proof}
\end{lemma}

\begin{lemma}
  \label{angle_energy_monotone}
  The map $F$ defined in
  \eqref{def_F_angle}
is strictly monotone decreasing and minimized for $\alpha=\frac{\pi}{2}$.
\end{lemma}
We decompose the function $F$ given
in \eqref{def_F_angle}, as  $F(\alpha)=H(\alpha)+G(\alpha)$, where
  \begin{align}
    H(\alpha)&:=\textstyle\frac{1}{(\sin\alpha)^{q-1}}\Big[
    \int_{-\cot\alpha}^{\tan\frac{\alpha}{2}} 
    \big(\frac{1}{z^2+1}\big)^q dz + 
    \int_{\cot\alpha}^{\cot\frac{\alpha}{2}} 
    \big(\frac{1}{z^2+1}\big)^q dz\Big], \label{eq:H}\\ 
    G(\alpha)&:=\textstyle\frac{1}{(\sin\alpha)^{q-1}}\Big[
    \int_{-\cot\alpha}^{\tan\frac{\alpha}{2}} 
    \big(\frac{z\sin\alpha+\cos\alpha}{z^2+1} \big)^q dz + 
    \int_{\cot\alpha}^{\cot\frac{\alpha}{2}} 
    \big(\frac{z\sin\alpha-\cos\alpha}{z^2+1} \big)^q dz\Big].\label{eq:G}
  \end{align}
  To show that $F$ is strictly monotone decreasing we analyze the first
  derivatives of $H$ and $G$. 
  By means of the Leibniz rule for parameter integrals we compute 
\begin{align}
\textstyle\frac{d}{d\alpha} H(\alpha) 
&=\textstyle -(q-1)\frac{\cos\alpha}{(\sin\alpha)^q}
\Big[\int_{\cot\alpha}^{\cot\frac{\alpha}{2}}\big( \frac{1}{z^2+1}\big)^q dz 
+ \int_{-\cot\alpha}^{\tan\frac{\alpha}{2}}\big( \frac{1}{z^2+1}\big)^q dz \Big] 
\notag\\
&+\textstyle\frac{1}{(\sin\alpha)^{q-1}}\Big[
\big(\frac{1}{(\cot\frac{\alpha}{2})^2+1}\big)^q 
\frac{d}{d\alpha}\cot\frac{\alpha}{2}
+\big(\frac{1}{(\tan\frac{\alpha}{2})^2+1}\big)^q
\frac{d}{d\alpha}\tan\frac{\alpha}{2} \Big].\label{eq:1111}
\end{align}
With the identities $(\cot\frac{\alpha}2)^2 +1=(\sin\frac{\alpha}2)^{-2}$, 
$\frac{d}{d\alpha}\cot\frac{\alpha}2=-\frac12 (\sin\frac{\alpha}2)^{-2}$, 
$(\tan\frac{\alpha}2)^2+1=(\cos\frac{\alpha}2)^{-2}$, and 
$\frac{d}{d\alpha}\tan\frac{\alpha}2=\frac12 (\cos\frac{\alpha}2)^{-2}$
we can replace the second line in \eqref{eq:1111} by
the expression
$
\textstyle\frac{1}{(\sin\alpha)^{q-1}}\big[\frac{1}{2}(\cos\frac{\alpha}{2})^{2q-2}
- \frac{1}{2}(\sin\frac{\alpha}{2})^{2q-2}\big]
$
to obtain
\begin{align*}
\textstyle\frac{d}{d\alpha} H(\alpha) &=\textstyle
\frac{1}{(\sin\alpha)^{q-1}} \Big[ -(q-1)(\cot\alpha)\Big(
\int_{\cot\alpha}^{\cot\frac{\alpha}{2}}\big( \frac{1}{z^2+1}\big)^q dz 
+ \int_{-\cot\alpha}^{\tan\frac{\alpha}{2}}\big( \frac{1}{z^2+1}\big)^q dz \Big) \\
&\textstyle +\frac{1}{2}(\cos\frac{\alpha}{2})^{2q-2}- 
\frac{1}{2}(\sin\frac{\alpha}{2})^{2q-2} \Big] .
\end{align*}
We use  the function
 $f(z):=\frac{1}{2}(z^2+1)^{1-q}$ with $f'(z)=(1-q)z(z^2+1)^{-q}$ to
 rewrite the last two summands in square brackets as integrals
\begin{align*}
\textstyle\int_{-\cot\alpha}^{\tan\frac{\alpha}{2}}
f'(z) dz - \int_{\cot\alpha}^{\cot\frac{\alpha}{2}} f'(z) dz 
&=
\textstyle (q-1)\left(\int_{\cot\alpha}^{\cot\frac{\alpha}{2}} 
\frac{z}{(z^2+1)^q} dz-\int_{-\cot\alpha}^{\tan\frac{\alpha}{2}}
\frac{z}{(z^2+1)^q} dz \right),
\end{align*}
so that we can deduce
\begin{align*}
\textstyle\frac{d}{d\alpha}H(\alpha)=\frac{(q-1)}{(\sin\alpha)^{q-1}} 
\left(\int_{\cot\alpha}^{\cot\frac{\alpha}{2}}
\frac{z-\cot\alpha}{(z^2+1)^q}dz -\int_{-\cot\alpha}^{\tan\frac{\alpha}{2}}
\frac{z+\cot\alpha}{(z^2+1)^q} dz\right).
\end{align*}
In the second integral we may substitute $y=z+2\cot(\alpha)$, use
the identity $\tan\frac{\alpha}{2}+2\cot\alpha=\cot\frac{\alpha}{2}$, 
and finally replace the new integration variable $y$ by $z$ again to
arrive at
\begin{align}
\textstyle\frac{d}{d\alpha}H(\alpha)
&=\textstyle\frac{(q-1)}{(\sin\alpha)^{q-1}}
\int_{\cot\alpha}^{\cot\frac{\alpha}{2}}(z-\cot\alpha)
\big(\frac{1}{(z^2+1)^q}-\frac{1}{((z-2\cot\alpha)^2+1)^q}\big)dz
\label{eq:ddaH}
\end{align}
and notice that the prefactor of the integral is strictly positive for
$\alpha\in (0,\pi)$.
For $0<\alpha < \frac{\pi}{2}$ we have  
\begin{align*}
  &(z^2+1)^{-q} < \big((z-2\cot(\alpha))^2+1\big)^{-q} 
  \; \Leftrightarrow \; -4z\cot(\alpha) + 4\cot(\alpha)^2 < 0 \;  \Leftrightarrow \;  z > \cot(\alpha), 
\end{align*}
so that the integrand $g(z):=(z-\cot(\alpha))\left([z^2+1]^{-q}
-[(z-2\cot(\alpha))^2+1]^{-q}\right)$ in
\eqref{eq:ddaH} satisfies
$g(z)<0$ for $z > \cot(\alpha)$ and $0<\alpha <\frac{\pi}{2}$,
and $g(\cot(\alpha))=0$. Furthermore, $g(z)=0$ for all $z \in \left[\cot(\alpha),\cot\left(\frac{\alpha}{2}\right)\right]$ if and only if $\alpha=\frac{\pi}{2}$ since $\cot\left(\frac{\pi}{2}  \right)=0$. Therefore, 
\begin{align}
\frac{d}{d\alpha}H(\alpha)\begin{cases}\textstyle
<0, &\text{ if } \alpha \in \left(0,\frac{\pi}{2}\right),\\
=0, &\text{ if } \alpha=\frac{\pi}{2}.
\end{cases}\label{eq:H-mon}
\end{align}
Thus, $H$ is strictly monotone decreasing and has a global minimum at 
$\alpha=\frac{\pi}{2}$. 

Next, we consider the function $G$ defined in \eqref{eq:G} and abbreviate
the sum of integrals in square brackets as $\tilde{G}(\alpha)$ so that
\begin{align}\label{eq:G-tildeG}
\textstyle\frac{d}{d\alpha}G(\alpha)= 
-(q-1)\frac{\cos\alpha}{(\sin\alpha)^q}\tilde{G}(\alpha)
+\frac{1}{(\sin\alpha)^{q-1}}\frac{d}{d\alpha}\tilde{G}(\alpha). 
\end{align}
Since $q>2$ and $\alpha\in (0,\frac{\pi}2]$ we 
have $\tilde{G}(\alpha)> 0$,
and according to \eqref{eq:G-tildeG} it suffices to prove that $\frac{d}{d\alpha}
\tilde{G}(\alpha)\le 0$ to find
\begin{align}
\frac{d}{d\alpha}G(\alpha)\begin{cases}\textstyle
<0, &\text{ if } \alpha \in \left(0,\frac{\pi}{2}\right),\\
\le 0, &\text{ if } \alpha=\frac{\pi}{2},
\end{cases}\label{eq:G-mon}
\end{align}
which together with \eqref{eq:H-mon} proves the lemma.

Applying the Leibniz rule for parameter integrals yields
\begin{align*}
\textstyle\frac{d}{d\alpha}\tilde{G}(\alpha)={}&
\textstyle
\int_{\cot\alpha}^{\cot\frac{\alpha}{2}}q
\frac{(z\sin\alpha-\cos\alpha)^{q-1}}{(z^2+1)^q}(z\cos\alpha+\sin\alpha)dz\\
&
\textstyle+\int_{-\cot\alpha}^{\tan\frac{\alpha}{2}}q
\frac{(z\sin\alpha+\cos\alpha)^{q-1}}{(z^2+1)^q}
(z\cos\alpha-\sin\alpha)dz\\
&\textstyle
+\left(
\frac{\cot\frac{\alpha}{2}\sin\alpha-\cos\alpha}{(\cot\frac{\alpha}{2})^2+1} 
\right)^q\frac{d}{d\alpha}\cot\frac{\alpha}{2} 
+\left(
\frac{\tan\frac{\alpha}{2}\sin\alpha+\cos\alpha}{(\tan\frac{\alpha}{2})^2+1} 
\right)^q\frac{d}{d\alpha}\tan\frac{\alpha}{2}. 
\end{align*}
Together with $
\cot\frac{\alpha}{2}\sin\alpha-\cos\alpha=1$
and  $
\tan\frac{\alpha}{2}\sin\alpha+\cos\alpha =1, $
we therefore obtain
\begin{align}
&\textstyle\frac{d}{d\alpha}\tilde{G}(\alpha)={}
\textstyle
\int_{\cot\alpha}^{\cot\frac{\alpha}{2}}q
\frac{(z\sin\alpha-\cos\alpha)^{q-1}}{(z^2+1)^q}(z\cos\alpha+\sin\alpha)dz
\nonumber \\
&\textstyle+\int_{-\cot\alpha}^{\tan\frac{\alpha}{2}}q
\frac{(z\sin\alpha+\cos\alpha)^{q-1}}{(z^2+1)^q}
(z\cos\alpha-\sin\alpha)dz
+\frac{1}{2}(\cos\frac{\alpha}{2})^{2q-2}- 
\frac{1}{2}(\sin\frac{\alpha}{2})^{2q-2}. \label{derivative_tilde_Gs}
\end{align}
First, we express the last two summands by means of the functions
$
g_1(z)
:=\frac{(z\sin\alpha-\cos\alpha)^q}{(z^2+1)^{q-1}}$
for $z \in [\cot\alpha,\cot\frac{\alpha}{2}]$
and $
g_2(z)
:=\frac{(z\sin\alpha+\cos\alpha)^q}{(z^2+1)^{q-1}}$ for 
$z \in [-\cot(\alpha),\tan\frac{\alpha}{2}]$
as
\begin{align}
\textstyle\frac{1}{2}(\cos\frac{\alpha}{2})^{2q-2}- 
\frac{1}{2}(\sin\frac{\alpha}{2})^{2q-2}
=\frac{1}{2}\big(\int_{-\cot\alpha}^{\tan\frac{\alpha}{2}}
g_2'(z) dz-\int_{\cot\alpha}^{\cot\frac{\alpha}{2}}g_1'(z)dz \big). \label{int_tilde_gs}
\end{align}
The derivatives can be computed as
\begin{align*}
g_1'(z)&=\textstyle
q\sin\alpha\big(\frac{z\sin\alpha-\cos\alpha}{z^2+1} \big)^{q-1}
-2z(q-1)\big(\frac{z\sin\alpha-\cos\alpha}{z^2+1} \big)^{q} 
\\
&\textstyle
=\frac{(z\sin\alpha-\cos\alpha)^{q-1}}{(z^2+1)^q}\big(q\sin\alpha
-(q-2)z^2\sin\alpha+2z(q-1)\cos\alpha \big),\\
g_2'(z)&
\textstyle
=\frac{(z\sin\alpha+\cos\alpha)^{q-1}}{(z^2+1)^q}\big(q\sin\alpha
-(q-2)z^2\sin\alpha-2z(q-1)\cos\alpha \big).
\end{align*}
With \eqref{int_tilde_gs} we can express the derivative in 
\eqref{derivative_tilde_Gs} purely in terms of integrals. For their 
integrands we 
calculate
\begin{align*}
&\textstyle
q\frac{(z\sin\alpha-\cos\alpha)^{q-1}}{(z^2+1)^q}(z\cos\alpha+\sin\alpha)
-\frac{1}{2}g_1'(z)\\
={}&\textstyle
\frac{(z\sin\alpha-\cos\alpha)^{q-1}}{(z^2+1)^q} \big(qz\cos\alpha+q\sin\alpha
-\frac{q}{2}\sin\alpha+\frac{q-2}{2}z^2\sin\alpha
-z(q-1)\cos\alpha \big)\\
={}&\textstyle
\frac{(z\sin\alpha-\cos\alpha)^{q-1}}{(z^2+1)^q}\big(\frac{q}{2}\sin\alpha
+\frac{q-2}{2}z^2\sin\alpha+z\cos\alpha \big)
\end{align*}
and
\begin{align*}
&\textstyle
q\frac{(z\sin\alpha+\cos\alpha)^{q-1}}{(z^2+1)^q}(z\cos\alpha-\sin\alpha)
+\frac{1}{2}g_2'(z)\\
={}&\textstyle
\frac{(z\sin\alpha+\cos\alpha)^{q-1}}{(z^2+1)^q}\big(-\frac{q}{2}\sin\alpha
-\frac{q-2}{2}z^2\sin\alpha+z\cos\alpha \big).
\end{align*}
Combining all calculations yields
\begin{align}
\textstyle\frac{d}{d\alpha}\tilde{G}(\alpha)&
\textstyle
=\int_{\cot\alpha}^{\cot\frac{\alpha}{2}}
\frac{(z\sin\alpha-\cos\alpha)^{q-1}}{(z^2+1)^q}
\big(\frac{q}{2}\sin\alpha+\frac{q-2}{2}z^2\sin\alpha+z\cos\alpha \big) dz
\notag\\
&\textstyle
+\int_{-\cot\alpha}^{\tan\frac{\alpha}{2}}
\frac{(z\sin\alpha+\cos\alpha)^{q-1}}{(z^2+1)^q}
\big(-\frac{q}{2}\sin\alpha-\frac{q-2}{2}z^2\sin\alpha+z\cos\alpha \big)dz.
\label{eq:ddaG}
\end{align}
In the second integral, we substitute again $y=z+2\cot(\alpha)$ 
to obtain
for that integral
\begin{align*}
&\textstyle
\int_{\cot\alpha}^{\cot\frac{\alpha}{2}}
\frac{(y\sin\alpha-\cos\alpha)^{q-1}}{((y-2\cot\alpha)^2+1)^q}
\big(-\frac{q}{2}\sin\alpha-\frac{q-2}{2}(y-2\cot\alpha)^2\sin\alpha
+(y-2\cot\alpha)\cos\alpha \big)dy \\
={}&\textstyle
\int_{\cot\alpha}^{\cot\frac{\alpha}{2}}
\frac{(y\sin\alpha-\cos\alpha)^{q-1}}{((y-2\cot\alpha)^2+1)^q}
\big(-\frac{q}{2}\sin\alpha-\frac{q-2}{2}y^2\sin\alpha-y\cos\alpha
+2(q-1)\cos\alpha(y-\cot\alpha)\big)dy.
\end{align*}
Hence, replacing the $y$-variable by $z$ and inserting this expression
into \eqref{eq:ddaG} yields
\begin{align*}
\textstyle
&\frac{d}{d\alpha}\tilde{G}(\alpha)={}
\textstyle
\int_{\cot\alpha}^{\cot\frac{\alpha}{2}}\Big[
\frac{(z\sin\alpha-\cos\alpha)^{q-1}}{(z^2+1)^q}
\big(\frac{q}{2}\sin\alpha+\frac{q-2}{2}z^2\sin\alpha+z\cos\alpha \big)\\
&\textstyle
+\frac{(z\sin\alpha-\cos\alpha)^{q-1}}{((z-2\cot\alpha)^2+1)^q}
\big(-\frac{q}{2}\sin\alpha-\frac{q-2}{2}z^2\sin\alpha-z\cos\alpha
+2(q-1)\cos\alpha(z-\cot\alpha)\big)\Big] dz.
\end{align*}
To prove that this expression is non-positive we show as before
 that the integrand is non-positive 
for all $z \in [ \cot\alpha,\cot\frac{\alpha}{2}]$, 
$\alpha\in (0,\frac{\pi}2]$.
As $z \sin\alpha-\cos\alpha\geq 0$ for $z  \geq \cot\alpha$ 
this is equivalent to showing
\begin{align}
  \label{long_inequality}
&\textstyle\frac{\frac{q}{2}\sin\alpha
+\frac{q-2}{2}z^2\sin\alpha+z\cos\alpha}{(z^2+1)^q}
\leq{} \frac{\frac{q}{2}\sin\alpha+\frac{q-2}{2}z^2\sin\alpha+z\cos\alpha-2(q-1)\cos\alpha(z-\cot\alpha)}{\left((z-2\cot\alpha)^2+1\right)^q}. 
\end{align}
Note that for $z\geq \cot\alpha$, $0<\alpha\leq \frac{\pi}{2}$,
and $q\ge 2$ one has
\begin{align*}
\textstyle
h(\alpha):=
\frac{q}{2}\sin\alpha+\frac{q-2}{2}z^2\sin\alpha+z\cos\alpha\geq \sin\alpha 
+\cot\alpha\cos\alpha=\frac{1}{\sin\alpha}>0.
\end{align*}
Hence, we can divide inequality \eqref{long_inequality} above by the 
factor $h(\alpha)$ for $\alpha\in (0,\frac{\pi}2]$ and multiply 
with 
$((z-2\cot\alpha)^2+1)^q$ to obtain the equivalent characterization of 
\eqref{long_inequality},
\begin{align}
&\textstyle
\big(\frac{(z-2\cot\alpha)^2+1}{z^2+1}\big)^q  
\le
\frac{q\sin\alpha+(q-2)z^2\sin\alpha+2z\cos\alpha-4(q-1)\cos\alpha(z-\cot\alpha)}{q\sin\alpha+(q-2)z^2\sin\alpha+2z\cos\alpha}\notag \\ 
\Leftrightarrow \quad & 
\textstyle\big( 1-4\cot\alpha\,\frac{z-\cot\alpha}{ z^2+1 } \big)^q 
\leq 1 -4(q-1)\cot\alpha\,\frac{z-\cot\alpha}{q+(q-2)z^2+2z\cot\alpha}. 
\label{eq:long-ineq}
\end{align}
The maximal value of the the function
$f_\alpha(z):=\frac{z-\cot\alpha}{z^2+1}$ for 
$z \in [ \cot\alpha,\cot\frac{\alpha}{2} ]$ 
for fixed $\alpha\in (0,\frac{\pi}2]$ can be determined with the help 
of its derivative
$
f_\alpha'(z)=\frac{z^2+1-2z(z-\cot\alpha)}{( z^2+1 )^2}
=\frac{1-z^2+2z\cot\alpha}{( z^2+1 )^2}
$
whose zeroes are given by 
\begin{align}
\label{roots}
\textstyle
z_1=\cot\alpha-\frac{1}{\sin\alpha}=-\tan\frac{\alpha}{2} \AND 
z_2=\cot\alpha+\frac{1}{\sin\alpha}=\cot\frac{\alpha}{2}.
\end{align}
Only $z_2$ is contained in the interval $
[ \cot\alpha,\cot\frac{\alpha}{2}]$, so that we conclude for
$\alpha\in (0,\frac{\pi}2]$
\begin{align*}
f_\alpha(\cot\alpha)&=0 \AND  \textstyle
f_\alpha(z_2)=\frac{\cot\frac{\alpha}{2}-\cot\alpha}{(\cot\frac{\alpha}{2})^2+1}
=\frac{1}{2}\tan\frac{\alpha}{2}>0. 
\end{align*}
Hence, $f_\alpha$ attains on $[\cot\alpha,\cot\frac{\alpha}2]$ its
global maximum at $z_2$, which yields 
for all $\alpha\in (0,\frac{\pi}2]$ and $z\in [\cot\alpha,
\cot\frac{\alpha}2 ]$
\begin{align*}
\textstyle
0 \leq 4\cot\alpha f_\alpha(z)\le 4\cot\alpha f_\alpha(z_2)
\overset{\eqref{roots}}{=} 4\cot\alpha\,
\frac{\cot\frac{\alpha}{2}-\cot\alpha}{(\cot\frac{\alpha}{2})^2+1}
\overset{\eqref{roots}}{=}\frac{2\cos(\alpha)}{1+\cos\alpha}\le\frac{2\cos\alpha}{2\cos\alpha}=1, 
\end{align*}
so that
\begin{align*}
\textstyle
\big(1-4\cot\alpha\frac{z-\cot\alpha}{z^2+1}\big)^q\leq 1-4\cot\alpha\frac{z-\cot\alpha}{z^2+1}\quad\Fo z \in [ \cot\alpha,\cot
 \frac{\alpha}{2} ],\,\alpha\in (0,\frac{\pi}2],
\end{align*}
since $q>2$. Therefore, in view of \eqref{eq:long-ineq} above, it suffices to
prove
\begin{align}
\textstyle
4\cot\alpha\,\frac{z-\cot\alpha}{z^2+1} \geq  4(q-1)\cot\alpha\,
\frac{z-\cot\alpha}{q+(q-2)z^2+2z\cot\alpha}\quad
\Fo
z\in [\cot\alpha,\cot\frac{\alpha}2],\,\alpha\in (0,\frac{\pi}2]. 
\label{eq:appconclusion}
\end{align}
For $\alpha=\frac{\pi}2$ both sides vanish, so we can restrict to
$\alpha\in (0,\frac{\pi}2)$ from now on. Also for $z=\cot\alpha$ both
sides are zero.
 For $z\in (\cot\alpha,\cot\frac{\alpha}2]$ and 
 $\alpha\in (0,\frac{\pi}2)$
  inequality \eqref{eq:appconclusion} is equivalent to
\begin{align}
&\textstyle
\frac{1}{z^2+1} \geq  \frac{q-1}{q+(q-2)z^2+2z\cot(\alpha)}
\Leftrightarrow 
-z^2+2z\cot(\alpha)+1 \geq 0,\label{eq:app-finalconclusion}
\end{align}
where we used that $q+(q-2)z^2+2z\cot\alpha =(q-1)(1+z^2)+
1-z^2+2z\cot\alpha>1-z^2+2z\cot\alpha\ge 0$ by means of \eqref{roots}.
By \eqref{roots}, the last inequality in \eqref{eq:app-finalconclusion}
holds for all $z \in (\cot\alpha,\cot\frac{\alpha}{2}]$ and we have proven $\frac{d}{d\alpha}\tilde{G}(\alpha)\leq 0$.
\qed

\section*{Acknowledgments}
The first author was 
supported by {the DFG}-Gra\-du\-ier\-ten\-kolleg 
\emph{Energy, Entropy, and Dissipative Dynamics (EDDy)},  
project no. 320021702/GRK2326. 
The second author's work was partially funded by DFG Grant 
no.\@ Mo 966/7-1 \emph{Geometric curvature functionals: 
energy landscape and discrete methods} (project no. 282535003). 
Substantial parts of the content of this paper are contained in the first 
author's Ph.D. thesis \cite{lagemann_diss}. 

\addtocontents{toc}{\SkipTocEntry}
\bibliographystyle{abbrvhref_ba}
\bibliography{Literaturverzeichnis.bib}

\end{document}